\documentclass [leqno] {amsart}
  
\usepackage {amssymb}
\usepackage {amsthm}
\usepackage {amscd}
\usepackage{bm}
\usepackage{pdflscape}

\usepackage{pstricks}
\usepackage{color}
\usepackage{graphicx}

\usepackage[all]{xy}
\usepackage {hyperref} 

\begin{document}

\theoremstyle{plain}
\newtheorem{introthm}{Theorem}
\newtheorem{proposition}[subsubsection]{Proposition}
\newtheorem{lemma}[subsubsection]{Lemma}
\newtheorem{corollary}[subsubsection]{Corollary}
\newtheorem{thm}[subsubsection]{Theorem}
\newtheorem*{thm*}{Theorem}
\newtheorem{conjecture}[subsubsection]{Conjecture}
\newtheorem{question}[subsubsection]{Question}

\theoremstyle{definition}
\newtheorem{definition}[subsubsection]{Definition}
\newtheorem{notation}[subsubsection]{Notation}
\newtheorem{condition}[subsubsection]{Condition}
\newtheorem{example}[subsubsection]{Example}

\theoremstyle{remark}
\newtheorem{remark}[subsubsection]{Remark}

\numberwithin{equation}{subsubsection}

\newcommand{\eq}[2]{\begin{equation}\label{#1}#2 \end{equation}}
\newcommand{\ml}[2]{\begin{multline}\label{#1}#2 \end{multline}}
\newcommand{\mlnl}[1]{\begin{multline*}#1 \end{multline*}}
\newcommand{\ga}[2]{\begin{gather}\label{#1}#2 \end{gather}}

\newcommand{\arir}{\ar@{^{(}->}}
\newcommand{\aril}{\ar@{_{(}->}}
\newcommand{\are}{\ar@{>>}}

\newcommand{\xr}[1] {\xrightarrow{#1}}
\newcommand{\xl}[1] {\xleftarrow{#1}}
\newcommand{\lra}{\longrightarrow}
\newcommand{\inj}{\hookrightarrow}

\newcommand{\mf}[1]{\mathfrak{#1}}
\newcommand{\mc}[1]{\mathcal{#1}}

\newcommand{\CH}{{\rm CH}}
\newcommand{\Gr}{{\rm Gr}}
\newcommand{\codim}{{\rm codim}}
\newcommand{\cd}{{\rm cd}}
\newcommand{\Spec} {{\rm Spec}}
\newcommand{\supp} {{\rm supp}}
\newcommand{\Hom} {{\rm Hom}}
\newcommand{\End} {{\rm End}}
\newcommand{\id}{{\rm id}}
\newcommand{\Aut}{{\rm Aut}}
\newcommand{\sHom}{{\rm \mathcal{H}om}}
\newcommand{\Tr}{{\rm Tr}}
\newcommand{\Res}{{\rm Res}}
\newcommand{\sExt}{\mathcal{E}xt}


\newcommand{\Grass}{\mathbb{G}{\rm r}} 

\newcommand{\sA}{{\mathcal A}}
\newcommand{\sB}{{\mathcal B}}
\newcommand{\sC}{{\mathcal C}}
\newcommand{\sD}{{\mathcal D}}
\newcommand{\sE}{{\mathcal E}}
\newcommand{\sF}{{\mathcal F}}
\newcommand{\sG}{{\mathcal G}}
\newcommand{\sH}{{\mathcal H}}
\newcommand{\sI}{{\mathcal I}}
\newcommand{\sJ}{{\mathcal J}}
\newcommand{\sK}{{\mathcal K}}
\newcommand{\sL}{{\mathcal L}}
\newcommand{\sM}{{\mathcal M}}
\newcommand{\sN}{{\mathcal N}}
\newcommand{\sO}{{\mathcal O}}
\newcommand{\sP}{{\mathcal P}}
\newcommand{\sQ}{{\mathcal Q}}
\newcommand{\sR}{{\mathcal R}}
\newcommand{\sS}{{\mathcal S}}
\newcommand{\sT}{{\mathcal T}}
\newcommand{\sU}{{\mathcal U}}
\newcommand{\sV}{{\mathcal V}}
\newcommand{\sW}{{\mathcal W}}
\newcommand{\sX}{{\mathcal X}}
\newcommand{\sY}{{\mathcal Y}}
\newcommand{\sZ}{{\mathcal Z}}

\newcommand{\A}{{\mathbb A}}
\newcommand{\B}{{\mathbb B}}
\newcommand{\C}{{\mathbb C}}
\newcommand{\D}{{\mathbb D}}
\newcommand{\E}{{\mathbb E}}
\newcommand{\G}{{\mathbb G}}
\renewcommand{\H}{{\mathbb H}}
\newcommand{\J}{{\mathbb J}}
\newcommand{\M}{{\mathbb M}}
\newcommand{\N}{{\mathbb N}}
\renewcommand{\P}{{\mathbb P}}
\newcommand{\Q}{{\mathbb Q}}
\newcommand{\U}{{\mathbb U}}
\newcommand{\X}{{\mathbb X}}
\newcommand{\Y}{{\mathbb Y}}
\newcommand{\Z}{\mathbb{Z}}

\newcommand{\abs}[1]{\lvert#1\rvert}
\newcommand{\OO}{\sO}
\renewcommand{\AA}{\mathcal{A}}
\newcommand{\CA}{\mathcal{C}\mathcal{A}}
\newcommand{\Bl}[2]{B_{#1}{#2}}
\newcommand{\St}[1]{{\rm St}(#1)}
\newcommand{\BL}{{\rm BL}}
\newcommand{\ugr}{\mathfrak{u}^{{\rm gr}}_{\omega}}
\newcommand{\ugrab}{(\ugr)^{{\rm ab}}}
\newcommand{\ti}[1]{\tilde{#1}}
\newcommand{\Rat}[1]{{\rm Rat}_{#1}}
\newcommand{\Zp}[1]{Z_{#1}}
\newcommand{\CHp}[1]{{\rm CH}_{#1}}
\newcommand{\prop}{{\rm prop}}
\newcommand{\Ext}{{\rm Ext}}
\newcommand{\triples}{(F_*,F^*,T,e)}
\newcommand{\rd}{{\rm red}}
\newcommand{\pr}{{\rm pr}}
\newcommand{\GrAb}{\text{\bf{GrAb}}}
\newcommand{\T}{\text{{\bf T}}}
\newcommand{\CatGrAb}{\text{\bf{Cat}}_{\GrAb,\otimes}}
\newcommand{\CatGrAbV}{\text{\bf{Cat}}_{V/\GrAb,\otimes}}
\newcommand{\Set}{\text{{\bf Sets}}}
\newcommand{\Cht}{(\CH_*,\CH^*,\times,1)}
\newcommand{\ul}[1]{\underline{#1}}

\title[Higher direct images of the structure sheaf]{Higher direct images of the structure sheaf in positive characteristic}

\author{Andre Chatzistamatiou and Kay R\"ulling}
\address{Fachbereich Mathematik \\ Universit\"at Duisburg-Essen \\ 45117 Essen, Germany}
\email{a.chatzistamatiou@uni-due.de}
\email{kay.ruelling@uni-due.de}

\thanks{This work has been supported by the SFB/TR 45 ``Periods, moduli spaces and arithmetic of algebraic varieties''}

\begin{abstract}
We prove vanishing of the higher direct images of the structure (and the canonical)
sheaf for a proper birational morphism with source a smooth variety 
and target the quotient of a smooth variety by a finite group of order prime to the characteristic
of the ground field. We also show that for smooth projective varieties the cohomology of the 
structure sheaf is a birational invariant. These results are well-known in characteristic zero. 
\end{abstract}

\maketitle

\tableofcontents

\section*{Introduction}
In characteristic zero it is a well-known and frequently used fact that 
the higher direct images $R^if_*\OO_X$ of a projective birational morphism $f:X\xr{} Y$ between smooth 
schemes vanish for $i>0$. This statement was proved as a corollary of Hironaka's resolution of 
singularities by resolving the indeterminacies of $f^{-1}$ by successively blowing-up smooth subvarieties of $Y$. In this article we consider the situation over an arbitrary field $k$ and prove this and related results.

In the following all schemes are assumed to be separated and of finite type over $k$, all morphisms are assumed to be $k$-morphisms.
The two main results of this paper are:
\begin{introthm}\label{thm1}
Assume $k$ is perfect. Let $S$ be an arbitrary scheme and let $X,Y$ be integral $S$-schemes. Assume that $X,Y$ are smooth over $k$
and properly birational over $S$, i.e. there exists an integral scheme $Z$ and a commutative diagram
\[\xymatrix@=13pt{   &  Z\ar[dl]_{\tau_X}\ar[dr]^{\tau_Y}  &   \\
              X\ar[dr]_f     &               & Y\ar[dl]^g  \\
                 &      S  &     } \]
such that $\tau_X,\tau_Y$ are proper and birational ($f$ and $g$ being the fixed morphisms to $S$).
Then for all $i$,  there are isomorphisms of $\OO_S$-modules 
\[R^if_*\OO_X \cong R^ig_*\OO_Y,\quad R^if_*\omega_X \cong R^ig_*\omega_Y.\]
\end{introthm}

\begin{introthm}\label{thm2}
Consider a diagram   
$$
\xymatrix
{
&
Y\ar[d]^{f}
\\
\ti{X}
\ar[r]^{\pi}
&
X,
}
$$
where $Y$ and $\ti{X}$ are connected smooth schemes, $X$ is integral and normal,  
$f$ is surjective and finite such that $\deg(f)\in k^*$, and $\pi$ is birational and proper.
Then X is Cohen-Macaulay and
$$
R\pi_*\OO_{\ti{X}}=\OO_X, \quad R\pi_*\omega_{\ti{X}}= \omega_{X}, 
$$
where $\omega_X$ is the dualizing sheaf of $X$.
\end{introthm}
By duality, the two identities in Theorem \ref{thm2} imply each other.

Both theorems are known in characteristic zero: Theorem \ref{thm1} follows from Hironaka's resolution of singularities,
for Theorem \ref{thm2} see \cite{Vi} (which also uses resolution of singularities).
If resolution of singularities is available in positive characteristic then it easily yields Theorem 
\ref{thm1}.

Recall from \cite[I, \S3, p.50]{Kempf}, that a rational resolution of an integral normal scheme $X$ is a resolution (i.e.~a proper birational morphism $g:\tilde{X}\xr{} X$ with $\tilde{X}$ smooth) which satisfies 
$R^ig_*(\OO_{\tilde{X}})=0=R^ig_*(\omega_{\tilde{X}})$ for all $i>0$. 
Thus Theorem \ref{thm1} implies that if an integral normal scheme over a perfect field has a rational resolution, 
then any resolution of $X$ is a rational resolution. 
For a smooth scheme $X$ we obtain $R^ig_*(\OO_{\tilde{X}})=0=R^ig_*(\omega_{\tilde{X}})$ for $i>0$ and any resolution $g$ 
(Corollary \ref{corollarystandartapplication}). Theorem \ref{thm2} asserts that 
$\pi:\tilde{X}\xr{} X$ is a rational resolution; this includes the important special 
case where $X$ is the quotient of $Y$ by a finite group of order prime 
to the characteristic of $k$.

Since resolution of singularities is not yet available in positive characteristic, we develop a 
different approach based on algebraic correspondences.
To get an idea of the methods involved let us sketch the proof of Theorem \ref{thm1} for $S=\Spec\, k$ and $X,Y$ projective (see \textsection3 for the details and a rigorous proof). 

For a scheme $X$ we write $H(X)=\bigoplus_{i,j} H^i(X,\Omega^j_X)$ and $\CH(X)=\oplus_i \CH^i(X)$, with $\CH^i(X)$ the Chow group of
codimension $i$ cycles on $X$. Given smooth projective schemes $X, Y, Z$,  
there is a composition of correspondences:
\begin{equation*}
\begin{split}
\CH(X\times Y)\otimes_\Z \CH(Y\times Z)&\lra \CH(X\times Z),\\ 
H(X\times Y)\otimes_{k} H(Y\times Z)&\to H(X\times Z). 
\end{split}
\end{equation*}
Moreover $a\in \CH(X\times Y)$ (resp.~$\in H(X\times Y)$) defines a map 
on the cohomology $\CH(X)\xr{} \CH(Y)$ (resp.~$H(X)\xr{} H(Y)$) by $c\mapsto \pr_{2*}(\pr^*_1(c)\cup a)$, and composition of correspondences corresponds to the composition of maps.  
Furthermore, there is a cycle map $cl: \CH\to H$, which is compatible with composition.

Now, the proof proceeds as follows.
By assumption there exists a closed integral subscheme $Z\subset X\times Y$  projecting birational to $X$ and $Y$. Let $Z'\subset Y\times X$ be its transpose.
By using the refined intersection product of Fulton we will see that
\[[Z']\circ [Z]= \id_X + E \quad \text{in }\CH\left(\pr_{13}\big((Z\times X)\cap (X\times Z')\big)\right),\]
with $E$ a cycle on $X\times X$ which projects on both sides to subsets of codimension at least one in $X$. 
We will show  that the map defined by $cl(E)\in H(X\times X)$ acts as zero on 
$H^*(X,\OO_X)\oplus H^*(X,\omega_X)$. A similar argument applies for $[Z]\circ [Z']$.
Thus the maps defined by $cl([Z])$ and $cl([Z'])$ are inverse to each other (when restricted to $H^*(X,\OO_X)\oplus H^*(X,\omega_X)$). This
proves Theorem \ref{thm1} in the case $S=\Spec\, k$ and $X,Y$ projective. 

It is not hard to deduce the general statement of Theorem \ref{thm1} once we know it in the case $S=\Spec\, k$. Therefore we have to generalize the above
argument to the case of smooth but not necessarily proper $k$-schemes.  
The problem is that in general a push-forward on $\CH$ or $H$ does not exist. 
However, the variety $Z\subset X\times Y$ is proper over $X,Y,$ and  
by working with cohomology (or Chow groups) with support we can conclude as outlined above. 

One of the main points in this paper is the construction of a cycle map, or natural transformation between cohomology theories with support, 
$
\CH \xr{} H.
$
For this, we first give a definition for (weak) cohomology theories with support.   
We introduce two categories $V^*$ and $V_*$. The objects in both categories are 
$(X,\Phi)$, where $X$ is smooth and $\Phi$ is a family of supports on $X$ (see \ref{def:famsupp} for the definition of a family with supports). 
A morphism $f: X\to Y$ induces a morphism $(X,\Phi)\to (Y,\Psi)$
in $V_*$ if and only if $f\mid \Phi$ is proper and $f(\Phi)\subset \Psi$; $f$ induces a morphism in $V^*$ if and only if $f^{-1}(\Psi)\subset \Phi$.

Then we consider the data $(F_*, F^*, T, e)$, where 
\[F_* : V_*\lra \GrAb, \quad F^*: (V^*)^{\rm op}\lra \GrAb,\]
are functors to graded abelian groups with $F_*(X,\Phi)=F^*(X, \Phi)=:F(X,\Phi)$ as abelian groups, 
$T$ gives for all $(X,\Phi)$, $(Y,\Psi)$ a morphism of abelian groups
\[T_{(X,\Phi), (Y,\Psi)}: F(X,\Phi)\otimes_\Z F(Y,\Psi)\lra F(X\times Y, \Phi\times \Psi),\]
and $e:\Z \to F(\Spec\, k)$ is a morphism of abelian groups. 

These data define a weak cohomology theory with support if the following conditions are satisfied.
\begin{enumerate}
\item The map $T$ is functorial for $F_*,F^*$. 
\item For all diagrams
       \[\xymatrix{ (X',\Phi')\ar[r]^{f'}\ar[d]^{g_X} & (Y',\Psi')\ar[d]^{g_Y}\\
                    (X,\Phi)\ar[r]^f   & (Y,\Psi), }\]
     in which the underlying diagram of schemes is cartesian and \emph{transversal}, 
     with $g_X, g_Y\in V^*$, $f,f'\in V_*$,  the following equality holds:
        \[F^*(g_Y)\circ F_*(f)= F_*(f')\circ F^*(g_X).\] 
\item Some more (very natural) conditions.
\end{enumerate}
The conditions allow us to obtain a calculus with correspondences. 
One key example of a weak cohomology theory with supports is the Chow group
$$
\CH(X,\Phi):=\varinjlim_{W\in \Phi} \CH(W),
$$ 
with $\CH_*$ the (proper) push forward for cycles and $\CH^*$ the refined Gysin homomorphism. 
Another example is the Hodge cohomology 
$
H(X,\Phi):=\bigoplus_{i,j} H^i_{\Phi}(X,\Omega^j_X).
$ 
Here, the definition of $H^*$ is straightforward, but for $H_*$ we use Grothendieck duality for singular schemes since smooth compactifications are not available in characteristic $p$. 
That the Hodge cohomology defines a (weak) cohomology theory with supports is a non-trivial 
fact, of which the proof occupies \textsection2.

In Theorem \ref{mainthm1} we give necessary and sufficient conditions for a (weak) cohomology
theory with supports $F$ to be target of a morphism from $\CH$. Unfortunately, we can do this only
with an additional semi-purity assumption on $F$ (see \ref{semipurity}). 
As an application we prove the existence of a cycle map $\CH\xr{} H$.
We hope that Theorem \ref{mainthm1} will turn out to be useful for proving similar results for the Witt vector cohomology.  

Let us give a short overview of the content of each section. 

In \textsection1 we define
weak cohomology theories with supports and prove basic properties. 
We show that $\CH$ is an example and prove Theorem \ref{mainthm1}. Moreover, we explain the
calculus of correspondences attached to a cohomology theory.

In \textsection2 we show that the Hodge cohomology is another example for a cohomology theory 
with supports. The hard part is the definition of push-forward maps. 
We use Grothendieck's duality theory for singular schemes as developed in \cite{Ha}, \cite{Co}, and
make extensive use of the results given in these references. There are also other approaches to
duality theory, which are more elegant (see e.g. \cite{Lipman2}). But since we use at several places the explicit description
of duality theory as developed by Grothendieck and since it is not clear to the authors, how this classical approach compares
to the one e.g. in \cite{Lipman2}, we will solely stick to the references \cite{Ha} and \cite{Co}.

In \textsection3 we show the existence of a cycle map $\CH \xr{} H$. Moreover, a vanishing statement 
Proposition \ref{propositionvanishing} is proved which enables us to prove Theorem \ref{thm1} and \ref{thm2}. 

In \textsection4 we generalize Theorem \ref{thm1} to the case where $X$ and $Y$ are tame quotients (see Theorem \ref{thm-tame-quotients}). 
This theorem also implies Theorem \ref{thm2}. 

We finish with some open questions. 

In the appendix finally we give a well-known description of the trace morphism for closed embeddings between smooth schemes and for finite and surjective morphisms between smooth schemes, 
which we need in \textsection2. 

\subsubsection*{Acknowledgments}
We are deeply grateful to H\'el\`ene Esnault for her encouragement and patience. 
We thank Manuel Blickle, Georg Hein and Eckart Viehweg for useful discussions. 
Finally, we thank the referee for a careful reading and for suggesting many improvements
in the exposition.

\section{Chow groups with support}

\subsection{Cohomology theories with support}
Let $k$ be a field. All schemes are assumed to be of finite type and separated over $k$.
We begin by recalling basic definitions and notations concerning families of supports.

\begin{definition}\label{def:famsupp}
A family of supports $\Phi$ on $X$ is a non-empty set of closed subsets of $X$ such that
the following holds:
\begin{itemize}
\item[(i)] The union of two elements in $\Phi$ is contained in $\Phi$.
\item[(ii)] Every closed subset of an element in $\Phi$ is contained in $\Phi$. 
\end{itemize}  
\end{definition}

Let $A$ be any set of closed subsets of $X$. The smallest family of supports 
$\Phi_A$ which contains $A$ is given by
\begin{equation}\label{Phiset}
\Phi_A:=\{\bigcup_{i=1}^n Z'_i\,;\, Z'_i \underset{\text{closed}}{\subset} Z_i\in A\}.
\end{equation}
For a closed subset $Z\subset X$ we write $\Phi_Z$ for $\Phi_{\{Z\}}$.

\begin{notation}
Let $f:X\xr{} Y$ be a morphism of schemes and $\Phi$ resp. $\Psi$ a family of 
supports of $X$ resp. $Y$.
\begin{enumerate}
\item We denote by $f^{-1}(\Psi)$ the smallest family of supports on $X$ which contains 
      $\{f^{-1}(Z);Z\in \Psi\}$.
\item We say that $f\mid \Phi$ is proper if $f\mid Z$ is proper for every $Z\in \Phi$.
If $f\mid\Phi$ is proper then $f(\Phi)$ is a family of supports on $Y$. 
\item If $\Phi_1,\Phi_2$ are two families of supports then $\Phi_1\cap \Phi_2$ is a family of supports.
\item If $\Phi$ resp. $\Psi$ is a family of supports of $X$ resp. $Y$ then we denote 
by $\Phi\times \Psi$ the smallest family of supports on $X\times_k Y$ which contains
$\{Z_1\times Z_2; Z_1\in \Phi, Z_2\in \Psi\}$.
\end{enumerate}
\end{notation}

When working with cohomology theories with support it is convenient to define the following
two categories $V_*$ and $V^*$, where for the morphisms in $V_*$ a ``push-forward'' map 
can be expected and for the morphisms in $V^*$ a ``pull-back'' map 
can be expected.   

\begin{definition}
We denote by $V_*$ the category with objects $(X,\Phi)$ where $X$ is a smooth scheme and
$\Phi$ is a family of supports of $X$, and morphisms
$$
\Hom_{V_*}((X,\Phi),(Y,\Psi))=\{f\in \Hom_k(X,Y); \text{$f\mid \Phi$ is proper and $f(\Phi)\subset \Psi$.}\}  
$$

We denote by $V^*$  the category with objects $(X,\Phi)$ where $X$ is a smooth scheme and
$\Phi$ is a family of supports of $X$ (${\rm ob}(V_*)={\rm ob}(V^*)$), and morphisms
$$
\Hom_{V^*}((X,\Phi),(Y,\Psi))=\{f\in \Hom_k(X,Y); f^{-1}(\Psi)\subset \Phi\}.  
$$ 

The composition and the identity comes in both cases from the category of schemes (over $k$).
\end{definition}

\subsubsection{}
Let $X$ be a smooth scheme. For a closed subscheme $W\subset X$ we write 
$(X,W):= (X,\Phi_W)$ in $V^*$ and $V_*$, respectively.
We simply write $X$ for $(X,X)$. 

We have forgetful functors 
$V_*\xr{} \text{\bf{Sch}}_k$ resp. $V^*\xr{} \text{\bf{Sch}}_k$ to the category of schemes, and
we often denote the morphism of schemes induced by a morphism in $V_*$ and
$V^*$ respectively by the same letter.   

For a morphism $f$  in $V_*$ we will say that $f$ is an immersion, flat, \dots, 
if the corresponding morphism of schemes has this property, 
and similarly for morphisms in $V^*$.    
We say that a diagram 
\begin{equation}
\label{MC}
\xymatrix
{
(X',\Phi') \ar[r]^{f'} \ar[d]^{g_X}
&
(Y',\Psi')  \ar[d]^{g_Y}
\\
(X,\Phi) \ar[r]^{f} 
&
(Y,\Psi)  
}
\end{equation}
is cartesian if the diagram of the corresponding schemes is cartesian.  

\subsubsection{Coproducts and ``products''}
For both categories $V^*$ and $V_*$ finite coproducts exist:
$$
(X,\Phi)\coprod (Y,\Psi)=(X\coprod Y,\Phi \cup \Psi).
$$  
For $(X,\Phi)$ let $X=\coprod_i X_i$ be the decomposition into connected components, 
then 
$$
(X,\Phi)=\coprod_i (X_i,\Phi\cap \Phi_{X_i}).
$$
In general products don't exist, and we define 
$$
(X,\Phi)\otimes (Y,\Psi):= (X\times Y,\Phi\times \Psi)
$$
which together with the unit object $\mathbf{1}={\rm Spec}(k)$ and the 
obvious isomorphism $(X,\Phi)\otimes (Y,\Psi) \xr{} (Y,\Psi)\otimes (X,\Phi)$ makes
$V_*$ and $V^*$ into a symmetric monoidal category (see \cite[VII.1]{ML}). 

\subsubsection{}\label{triples}
Consider the following data $(F_*,F^*,T,e)$:  
\begin{enumerate}
\item Two functors $F_*:V_*\xr{}\GrAb$ and $F^*:(V^*)^{op} \xr{} \GrAb$ to the symmetric monoidal category of graded abelian groups, such that 
$F_*(X)=F^*(X)$ as ungraded groups for every object $X\in {\rm ob}(V_*)={\rm ob}(V^*)$. 
We will simply write $F(X):=F_*(X)=F^*(X)$. We use lower indexes 
for the grading on $F_*(X)$, i.e. $F_*(X)=\oplus_i F_i(X)$, and upper indexes for $F^*(X)$.  
\item For every two objects $X,Y\in {\rm ob}(V_*)={\rm ob}(V^*)$ a morphism of graded abelian 
groups (for both gradings):
$$
T_{X,Y}: F(X)\otimes_{\Z} F(Y) \xr{} F (X\otimes Y).
$$
\item A morphism of abelian groups $e:\Z\xr{} F(\Spec(k))$. For all smooth schemes $\pi:X\xr{} \Spec(k)$ we denote by $1_X$ the image of 
$1\in \Z$ via the map 
$
\Z\xr{e} F^*(\Spec(k))\xr{F^*(\pi)} F^*(X).
$ 
\end{enumerate}

\subsubsection{} \label{triplescond} 
The data $(F_*,F^*,T,e)$ is called a \emph{weak cohomology theory with supports}
if the following conditions are satisfied:
\begin{enumerate}
\item \label{triplescondsum} The functor $F_*$ preserves coproducts and $F^*$ maps coproducts to products. 
Moreover, for $(X,\Phi_1),(X,\Phi_2)\in {\rm ob}(V_*)$ with $\Phi_1\cap \Phi_2=\{\emptyset\}$
we require that the map  
$$
F^*(\jmath_1)+F^*(\jmath_2):F^*(X,\Phi_1)\oplus F^*(X,\Phi_2) \xr{} F^*(X,\Phi_1\cup \Phi_2),
$$
with $\jmath_1:(X,\Phi_1\cup \Phi_2)\xr{} (X,\Phi_1)$ and 
$\jmath_2:(X,\Phi_1\cup \Phi_2)\xr{} (X,\Phi_2)$ in $V^*$, is an isomorphism.  
\item The data $(F_*,T,e)$ and $(F^*,T,e)$ respectively define a (right-lax) symmetric monoidal functor (see below).  
\item (Grading) For $(X,\Phi)$ such that $X$ is connected the equality   
$$F_i(X,\Phi)=F^{2\dim X-i}(X,\Phi)$$ holds for all $i$.
\item \label{MP} For all cartesian diagrams \ref{MC} with $g_X,g_Y\in V^*$ and  $f,f'\in V_*$ such that
either $g_Y$ is smooth or $g_Y$ is a closed immersion and $f$ is transversal to $g_Y$ the 
equality 
\begin{equation*}
 F^*(g_Y)\circ F_*(f) = F_*(f')\circ F^*(g_X) 
\end{equation*}
holds.
\end{enumerate}

Recall that $f$ is transversal to $g_Y$ if $(f')^*N_{Y'/Y}=N_{X'/X}$ where $N$ denotes
the normal bundle. The case $X'=\emptyset$ is also admissible; in this case the equality
\ref{triplescond}(\ref{MP}) reads:
$$
 F^*(g_Y)\circ F_*(f) = 0.
$$
The condition \ref{triplescond}(\ref{MP}) implies the projection formula (see Proposition \ref{projectionformula}) and will be needed for a calculus with correspondences.

Recall that $(F_*,T,e)$ is called a right-lax symmetric monoidal functor if 
\begin{itemize}
\item $T$ is associative, i.e. for $X,Y,Z\in {\rm ob}(V_*)$ 
the following diagram is commutative: 
$$
\xymatrix
{
F(X)\otimes F(Y)\otimes F(Z) \ar[r]^{id\otimes T}\ar[d]^{T\otimes id}
&
F(X)\otimes F(Y\otimes Z) \ar[d]^{T}
\\
F(X \otimes Y)\otimes F(Z) \ar[r]^{T} 
&
F(X\otimes Y\otimes Z).
}
$$
\item $T$ is commutative, i.e. for $X,Y\in {\rm ob}(V_*)$ the diagram
$$
\xymatrix
{
F(X)\otimes F(Y) \ar[r]^{T} \ar[d]
&
F(X\otimes Y) \ar[d]
\\
F(Y)\otimes F(X) \ar[r]^{T}
&
F(Y\otimes X)
}
$$ 
is commutative. Here, for two graded abelian groups $A,B$ the morphism 
$A\otimes B\xr{} B\otimes A$ maps $a\otimes b\mapsto (-1)^{\deg(a)\deg(b)} b\otimes a$.
\item The map $e:\Z\xr{} F(\Spec(k))$ renders commutative the diagrams
$$
\xymatrix
{
F(X)\otimes_{\Z} \Z \ar[r]^-{id\otimes e} \ar[rd]^{=} 
&
F(X)\otimes F(\Spec(k))\ar[r]^{T}
&
F(X\otimes \Spec(k)) \ar[dl]^{=}
\\
&
F(X),
&
\\
\Z\otimes_{\Z} F(X) \ar[r]^-{e\otimes id} \ar[rd]^{=}
&
F(\Spec(k))\otimes F(X)\ar[r]^{T}
&
F(\Spec(k)\otimes X) \ar[dl]^{=}
\\
&
F(X).
&
}
$$
\item $T$ is a natural transformation 
\begin{equation*}
T:(V_*\times V_* \xr{F_*\times F_*} \GrAb \times \GrAb \xr{\otimes} \GrAb)
\xr{} (V_*\times V_* \xr{\otimes} V_* \xr{F_*} \GrAb).
\end{equation*}
\end{itemize}

\begin{example}
A first example which satisfies these conditions (see Proposition \ref{tripleschowgroups})
is the Chow group
$$
(X,\Phi)\mapsto {\varinjlim}_{W\in \Phi}\CH_{*}(W).
$$   
The push-forward $V_*\xr{} \GrAb$ is defined in the usual way. To define the pull-back 
$(V^*)^{op}\xr{} \GrAb$ we use Fulton's refined Gysin homomorphism.  
However, in order to get a symmetric functor we have to put the Chow group 
$\CH_d(W)$ in degree $=2d$.

It will be shown in \textsection2 that the Hodge cohomology with support
$$
(X,\Phi)\mapsto \bigoplus_{i,j} H^i_{\Phi}(X,\Omega^j_X)
$$
is another example.
The push-forward is an application of Grothendieck's duality theory.

An example which isn't further considered in this paper is the \'etale cohomology
with support
$$
(X,\Phi)\mapsto H^{*}_{\Phi}(X\times \bar{k},\Q_{\ell}). 
$$ 
\end{example}

\begin{definition}\label{definitionT}
Let $\triples, (G_*,G^*,U,\epsilon)$ be as in \ref{triples} and satisfying the conditions 
\ref{triplescond}. 
By a morphism 
\begin{equation}\label{morphismsT}
\triples \xr{} (G_*,G^*,U,\epsilon)
\end{equation}
we understand a morphism of graded abelian groups (for both gradings)
$$
\phi: F(X)\xr{} G(X) \quad \text{for every $X\in {\rm ob}(V_*)= {\rm ob}(V^*)$,}
$$
such that $\phi$ induces a natural transformation of (right-lax) symmetric monoidal functors 
$$
\phi:(F_*,T,e) \xr{} (G_*,U,\epsilon)\; \text{and}\; \phi:(F^*,T,e) \xr{} (G^*,U,\epsilon),
$$ 
i.e. $\phi$ induces natural transformations $F_*\xr{} G_*, F^*\xr{} G^*$, and 
\begin{equation}\label{morphismTsecond}
\phi\circ T = U \circ(\phi \otimes \phi), \quad \phi\circ e=\epsilon. 
\end{equation}
We denote by $\T$ the category of weak cohomology theories with supports, i.e.~it is the category consisting of objects $\triples$ as in \ref{triples}, and 
satisfying the properties \ref{triplescond}, together with  morphisms \ref{morphismsT}.
\end{definition}

\subsubsection{Cup product} Let $(F_*,F^*,T,e)\in \T$.
For all $(X,\Phi)\in {\rm ob}(V^*)$ we obtain a cup product:
$$
\cup : F(X,\Phi_1)\otimes F(X,\Phi_2) \xr{T} F(X\times X,\Phi_1\times \Phi_2) \xr{F^*(\Delta_X)}
F(X,\Phi_1\cap \Phi_2), 
$$
where $\Delta_X:(X,\Phi_1\cap \Phi_2)\xr{} (X\times X,\Phi_1\times \Phi_2)$ is induced
by the diagonal immersion. The cup product is associative and graded commutative.  

By functoriality we obtain 
\begin{equation}\label{pb/cup}
F^*(f_1)(a)\cup F^*(f_2)(b)=F^*(f_3)(a\cup b)
\end{equation}
for all morphisms $f_1:(X',\Phi_1')\xr{} (X,\Phi_1), f_2:(X',\Phi_2')\xr{} (X,\Phi_2)$
in $V^*$ with $f_1=f_2:=f$ as morphisms of schemes; $f_3:(X',\Phi_1'\cap \Phi_2')\xr{} 
(X,\Phi_1\cap \Phi_2)$ in $V^*$ is induced by $f$. 

\begin{proposition}[Projection formula] \label{projectionformula} 
Let $\triples\in \T$ and let $f:X\xr{} Y$ be a morphism between smooth schemes, 
inducing morphisms 
\begin{equation*}
\begin{split}
f_1:&(X,\Phi_1)\xr{} (Y,\Phi_2) \quad \text{in $V_*$},\\
f_2:&(X,f^{-1}(\Psi))\xr{} (Y,\Psi) \quad \text{in $V^*$}.
\end{split}
\end{equation*}
Then $f$ also induces a morphism
$$f_3\,:\, (X,\Phi_1\cap f^{-1}(\Psi))\xr{} (Y,\Phi_2\cap\Psi) \quad \text{in } V_* $$
and for all $a\in F(X,\Phi_1), b\in F(Y,\Psi)$ the following formulas hold in $F(Y,\Phi_2\cap \Psi)$:
\begin{equation*}
\begin{split}
F_*(f_3)(a\cup F^*(f_2)(b))&=F_*(f_1)(a)\cup b, \\
F_*(f_3)(F^*(f_2)(b)\cup a)&=b \cup F_*(f_1)(a).
\end{split}
\end{equation*}

\begin{proof}
We prove the first equality of the statement, the second is proved in the same way. 
The diagram 
$$
\xymatrix
{
(X,\Phi_1\cap f^{-1}(\Psi)) \ar[r]^{f_3} \ar[d]_{\Delta_X}
&
(Y,\Phi_2\cap \Psi) \ar[dd]^{\Delta_Y}
\\
(X\times X,\Phi_1\times f^{-1}(\Psi))\ar[d]_{id\times f_2}
&
\\
(X\times Y,\Phi_1 \times \Psi) \ar[r]^{f_1\times id}
&
(Y\times Y,\Phi_2\times \Psi)
}
$$
is cartesian and $f\times id$ is transversal to $\Delta_Y$. Thus by \ref{triplescond}(\ref{MP}) we get
\begin{multline*}
F_*(f_3)(a\cup F^*(f_2)(b)) = F_*(f_3)F^*(\Delta_X)F^*(id\times f_2)(T(a\otimes b))=\\
F^*(\Delta_Y)F_*(f_1\times id)(T(a\otimes b)) = F_*(f_1)(a)\cup b .
\end{multline*}
\end{proof}
\end{proposition}


The proof of the following Lemma is straightforward.
\begin{lemma}\label{lemmaaboutone}
\begin{enumerate}
\item For all $(X,\Phi)$ and $a\in F(X,\Phi)$ the equality 
$$
1_X\cup a = a =a\cup 1_X 
$$
holds. In particular $F^*(X)$ is a (graded) ring. 
\item For smooth schemes $X,Y$ we have 
$$
T(1_X\otimes 1_Y)=1_{X\times Y}.
$$
\end{enumerate}
\end{lemma}

\subsubsection{Definition of Chow groups with support} 
In the following we define a first example of an object $(\CH_*,\CH^*,\times,e)\in \T$.

\begin{definition}[Chow groups with support]
Let $\Phi$ be a family of supports on $X$. We define: 
\begin{equation*}
\begin{split}
Z_{\Phi}(X):&=\text{abelian group freely generated by irreducible closed subsets $Z\in \Phi$.}\\
{\rm Rat}_{\Phi}(X):&=\text{subgroup of $Z_{\Phi}(X)$ generated by  $div(f)$, where $f\in k(W)^*$} \\  
& \phantom{=}\,\text{ is a non-zero rational function and $W\in \Phi$ is irreducible.}\\
\CH(X,\Phi):&=Z_{\Phi}(X)/\Rat{\Phi}(X).
\end{split}
\end{equation*}
\end{definition}

For $(X,\Phi), (Y,\Psi)$ we obtain 
\begin{equation}\label{coproductchow}
\CH((X\coprod Y,\Phi\cup \Psi))= \CH(X,\Phi)\oplus \CH(Y,\Psi).
\end{equation}

\subsubsection{Grading}
The groups $Z_{\Phi}(X)$ and $\Rat{\Phi}(X)$ can be graded by dimension: 

$$
\CH_{*}(X,\Phi)=\bigoplus_{d\ge 0} \CH_{d}(X,\Phi)[2d],
$$
where the bracket $[2d]$ means that $\CH_{d}(X,\Phi)$ is considered to be in degree $2d$. 

There is also a grading by codimension. Let $X=\coprod_i X_i$ be the decomposition into 
connected components then 
$
\CH^*(X,\Phi)=\bigoplus_i \CH^*(X_i,\Phi\cap \Phi_{X_i})
$ 
and 
$$
\CH^*(X_i,\Phi\cap \Phi_{X_i}) = \bigoplus_{d\geq 0} \CH^{d}(X_i,\Phi\cap \Phi_{X_i})[2d]
$$
where $\CH^{d}(X_i,\Phi)$ is generated by cycles $[Z]$ with $Z\in \Phi\cap \Phi_X$, $Z$ irreducible, and $\codim_{X_i}(Z)=d$.

\subsubsection{Examples}
If $W\subset X$ is a closed subset then we get 
$$
\CH(X,\Phi_W)=\CH(X,W)=\CH(W),
$$ 
the usual Chow group of $W$.

If $X$ is proper, $U$ is affine, and $\Phi:=\{Z'\,;\, Z'\subset U\}$ then 
$$
\CH(X,\Phi)=Z_{\Phi}(X)=\text{freely generated by closed points of $U$.}
$$ 

\subsubsection{Push forward for Chow groups}
Let $\Phi$ be a family of supports on $X$ (\ref{def:famsupp}).
If $W\subset X$ is a closed subscheme with $W\in \Phi$ 
then $Z(W)=Z_{\Phi_W}(X)\subset Z_{\Phi}(X)$, $\Rat{}(W)=\Rat{\Phi_W}(X)\subset \Rat{\Phi}(X)$ ($\Phi_W$ as 
defined in \ref{Phiset}), and we obtain a map
\begin{equation}\label{pushforwardZUX}
\CH(W) = \CH(X,W) \xr{} \CH(X,\Phi).
\end{equation}
Obviously, $\CH(X,\Phi)$ is the largest quotient of $Z_{\Phi}(X)$ such that there 
are push-forward maps \ref{pushforwardZUX} for every $W\in \Phi$. 

\subsubsection{}\label{pushforwardgeneral}
In general let $f:(X,\Phi)\xr{} (Y,\Psi)$ be a morphism in $V_*$.
There is a push-forward of cycles 
$$
f_*:\Zp{\Phi}(X)\xr{} \Zp{\Psi}(Y), \quad f_*([Z])=\deg(Z/f(Z))\cdot [f(Z)],
$$
for $Z\in \Phi$ irreducible ($\deg(Z/f(Z))=0$ if $\dim(f(Z))<\dim(Z)$).
Push-forward is functorial \cite[\textsection1.4]{F}.

\begin{lemma} With the assumption of \ref{pushforwardgeneral} we get 
 $f_*(\Rat{\Phi}(X))\subset \Rat{\Psi}(Y)$. 
\begin{proof}
Indeed, $\Rat{\Phi}(X)$ is generated
by the images of $\Rat{}(W)$ where $W\in \Phi$.  
The restriction $f\mid W$ is proper and \cite[Proposition~1.4]{F} yields 
$$
f_*({\rm Rat}_{}(W))\subset {\rm Rat}_{}(f(W)),
$$ 
which proves the claim. 
\end{proof}
\end{lemma}

Thus we get an induced map
\begin{equation}
f_*:\CH(X,\Phi)\xr{} \CH(Y,\Psi),
\end{equation}
and a functor 
\begin{equation}
\CH_*: V_* \xr{} \GrAb, \quad \CH_*(X,\Phi):=\CH(X,\Phi), \quad \CH_*(f):=f_*.
\end{equation}

\begin{proposition}\label{indlimit}
Let $\Phi$ be a family of supports of $X$. The following map
is an isomorphism 
$$
{\varinjlim}_{W\in \Phi}\CH(X,W) \xr{} \CH(X,\Phi).
$$
\end{proposition}
\begin{proof}
 This is straightforward.
\end{proof}

\subsubsection{Pull-back for Chow groups}  \label{Fulton}
In order to define a functor 
$$
\CH^*:(V^*)^{op} \xr{} \GrAb
$$
we have to recall Fulton's work on refined Gysin morphisms \cite[\textsection6.6]{F}.

Let $f:X\xr{} Y$ be a morphism between smooth schemes and let $V\subset Y$ be a closed
subscheme. There is a morphism of abelian groups 
$$
f^!:\CH(V) \xr{} \CH(f^{-1}(V))
$$
($f^{-1}(V)= X\times_Y V$)
with the following properties:
\begin{enumerate}
\item For a closed subscheme $V'\subset Y$ with $V\subset V'$ (denote the immersion by $\imath$,
and the immersion $f^{-1}(V)\subset f^{-1}(V')$ by $\jmath$), the equality 
$$
f^! \imath_* = \jmath_* f^!  
$$
(as maps $\CH(V)\xr{} \CH(f^{-1}(V'))$) holds. 
\item If $g:Y\xr{} Z$ is another morphism between smooth schemes and $S\subset Z$ a closed
subscheme then 
$$
f^!\circ g^! = (g\circ f)^!
$$ 
as maps $\CH(S)\xr{} \CH((g\circ f)^{-1}(S))$.
\item 
If $f:X\xr{} Y$ is flat then $f^!=f^*$ where $f^*$ is the usual pull-back map for 
flat morphisms.
\item 
Let 
\begin{equation*}
\xymatrix
{
X' \ar[r]^{f'} \ar[d]^{g_X}
&
Y'  \ar[d]^{g_Y}
\\
X \ar[r]^{f} 
&
Y
} 
\end{equation*}
be a cartesian diagram of smooth schemes and $W\subset X$ a closed subscheme such that
$f\mid W$ is proper. Assume that either $g_Y$ is flat or $g_Y$ is a closed immersion and 
$f$ is transversal to $g_Y$. Then 
$$
g_Y^!f_*=f'_* g_X^!
$$
as maps $\CH(W)\xr{} \CH(f'(g_X^{-1}W))=\CH(g_Y^{-1}f(W))$. This statement is proved
in \cite[Proposition~1.7]{F} for flat morphisms and in \cite[Theorem~6.2(a),(c)]{F} for the
case of a closed immersion.
\end{enumerate}

\begin{remark}
Note that $\CH(W)=\CH(W_{red})$ for every scheme $W$. 
\end{remark}

\subsubsection{Definition of the pull-back map}
Let $f:X\xr{} Y$ be a morphism between smooth schemes and let $V\subset Y$ be a closed
subscheme, thus $f:(X,f^{-1}(V))\xr{} (Y,V)$ is a morphism in $V^*$. We define 
$$
\CH^*(f):=f^!:\CH(Y,V)=\CH(V)\xr{} \CH(f^{-1}(V))= \CH(X,f^{-1}(V)). 
$$ 
For the general case let $f:(X,\Phi)\xr{} (Y,\Psi)$ be any morphism in $V^*$.
For every $V\in \Psi$ the map $f$ induces $(X,f^{-1}(V))\xr{} (Y,V)$ in $V^*$.
Because of \ref{Fulton}(1) and Proposition \ref{indlimit} we obtain 
$$
\CH^*(f): \CH(Y,\Psi)={\varinjlim}_{V\in \Psi}\CH(Y,V) \xr{} {\varinjlim}_{W\in \Phi}\CH(X,W)=\CH(X,\Phi).
$$

The assignment 
\begin{equation}
\CH^*:(V^*)^{op} \xr{} \GrAb, \quad \CH^*(X,\Phi)=\CH(X,\Phi), \quad f\mapsto \CH^*(f),
\end{equation}
defines a functor by \ref{Fulton}(1),(2).

\begin{proposition}\label{tripleschowgroups}
Together with the exterior product $\times$ (see \cite[\textsection1.10]{F}) and the obvious 
unit $1:\Z\xr{} \CH(\Spec(k))$, we obtain an object
$(\CH_*,\CH^*,\times,e)\in \T$. 
\end{proposition}
\begin{proof}
The formula
\ref{triplescond}(\ref{MP}) follows from \ref{Fulton}(1),(4).
\end{proof}

\subsection{Chow groups with support as initial object} 
Given $\Cht$ we are interested in objects $\triples \in \T$ which admit a morphism 
$\Cht\xr{} \triples$. Such morphisms should be viewed as a kind of cycle map, which is compatible with push-forward and pull-back.
Unfortunately, we can only give a satisfactory answer under an additional hypotheses on $\triples$, which we call semi-purity. 

\begin{definition}[Semi-purity]\label{semipurity}
We say that $\triples$ satisfies the semi-purity condition if the following holds:
\begin{itemize}
\item For all smooth schemes $X$ and irreducible closed subsets $W\subset X$ the 
groups 
$
F_i(X,W) 
$
vanish if $i>2\dim W$.
\item For all smooth schemes $X$, closed subsets $W\subset X$, and open
sets $U\subset X$ such that $U$ contains the generic point of every irreducible 
component of $W$, we require the map 
$$
F^*(\jmath):F_{2\dim W}(X,W)\xr{} F_{2\dim W}(U,W\cap U), 
$$ 
induced by $\jmath:(U,W\cap U)\xr{} (X,W)$ in $V^*$, to be \emph{injective}.
\end{itemize}
\end{definition}  

\begin{remark}
For $\Cht$ the condition is satisfied since $\CH_{2\dim W}(X,W)=\Z\cdot [W]$ and $\CH_{i}(X,W)=0$
for $i>2\dim W$.

Let $c$ be the codimension of $W$ in $X$, so that $F_{2\dim W}=F^{2c}$. 
Whenever there are exact sequences 
$$
F^{2c}(X,W\backslash U) \xr{} F^{2c}(X,W)\xr{} F^{2c}(U,U\cap W),
$$
the conditions in \ref{semipurity} follow from $F^{i}(X,W)=0$ for $i<2c$ (and all pairs $(X,W)$), 
which is known as semi-purity in the literature. 
\end{remark}

\begin{thm} \label{mainthm1}
Let $k$ be a perfect field and assume $\triples\in \T$ satisfies the semi-purity condition \ref{semipurity}.
Then $\Hom_{\T}(\Cht,\triples)$ is either empty or contains only one element.
The set $\Hom_{\T}(\Cht,\triples)$ is non-empty if and only if the following 
conditions hold:
\begin{enumerate}
\item \label{traceofone} If $f:X\xr{} Y$ is a finite morphism between smooth connected schemes of equal dimension 
then
$$
F_*(f)(1_X)=\deg(f)\cdot 1_Y.
$$ 
\item \label{rationaleq} For the $0$-point $\imath_0:\Spec(k)\xr{}\P^1$ and the $\infty$-point 
$\imath_{\infty}:\Spec(k)\xr{}\P^1$ the following equality holds:
$$F_*(\imath_0)\circ e=F_*(\imath_{\infty})\circ e.$$ 
\item \label{multiplicities} 
For a closed immersion $\imath:X\xr{} Y$ between smooth schemes
and an effective smooth divisor $D\subset Y$ such that 
\begin{itemize}
\item $D$ meets $X$ properly, thus $D\cap X:=D\times_{Y} X$ is a divisor on X,
\item $D':=(D\cap X)_{red}$ is smooth and connected, and thus $D\cap X=n\cdot D'$
as divisors (for some $n\in \Z, n\geq 1$), 
\end{itemize}
we denote by $\imath_X:X\xr{}(Y,X), \imath_{D'}:D'\xr{} (D,D')$ 
the morphisms in $V_*$ induced by $\imath$, and we define $g:(D,D')\xr{} (Y,X)$ in $V^*$ by  
the inclusion $D\subset Y$. Then the following equality is required to hold:
$$
F^*(g)(F_*(\imath_X)(1_X))=n\cdot F_*(\imath_{D'})(1_{D'}).
$$    
\item \label{class} If $W\subset X$ is an irreducible closed subset then there is an element 
$cl_{(X,W)}\in F_{2\dim W}(X,W)$ with 
$$
F^*(\jmath)(cl_{(X,W)})=F_*(\imath)(1_{U \cap W})
$$
for all open sets $U\subset X$ such that $U\cap W\neq \emptyset$ is smooth, and where 
$\jmath:(U,W\cap U)\xr{} (X,W)$ in $V^*$, $\imath:W\cap U\xr{} (U,W\cap U)$ in $V_*$. 
\end{enumerate}
\end{thm}

We will give the proof after the proof of the following Proposition.

\begin{proposition}\label{propositionmainthm1}
Let $k$ be a perfect field and let 
$F:=\triples\in \T$ satisfy the semi-purity condition \ref{semipurity}. Furthermore, we assume
that the conditions of Theorem (\ref{mainthm1}),(\ref{traceofone})-(\ref{class}) hold for $F$.
Then there is a unique natural transformation of (right-lax) symmetric monoidal functors 
$$
\phi: (\CH_*,\times,1) \xr{} (F_*,T,e)
$$
such that $\phi(1_X)=1_X$ for every smooth scheme $X$.
\begin{proof}
\emph{Uniqueness:} In view of the semi-purity condition \ref{semipurity}, 
\begin{equation}\label{cycleclassa}
\phi([W])=cl_{(X,W)} 
\end{equation} 
is the only choice for an irreducible closed subset $W$ of $X$, $[W]\in \CH_*(X,W)$. 
For a general family of supports $\Phi$ of $X$ the group $\CH_*(X,\Phi)$ is generated
by the images of $[W]$ via $\CH_*(X,W)\xr{} \CH_*(X,\Phi)$ where $W$ runs through all 
irreducible closed subsets $W\in \Phi$. 

\emph{Existence:}
For every smooth scheme $X$ and a family of supports $\Phi$ of $X$ we define
a homomorphism of abelian groups 
\begin{equation}\label{cycleclassmap2}
\phi'_{(X,\Phi)}: Z_{\Phi}(X)\xr{} F(X,\Phi) 
\end{equation}
by $\phi'_{(X,\Phi)}([W])=F_*(\imath_W)(cl_{(X,W)})$ 
for every irreducible closed subset $W\in \Phi$
and $\imath_W:(X,W)\xr{} (X,\Phi)$ in $V_*$ induced by $id_X$. 

\emph{1st Step:} For every morphism $f:(X,\Phi)\xr{} (Y,\Psi)$ in $V_*$ the push-forward $f_*:Z_{\Phi}(X)\xr{}
Z_{\Psi}(Y)$ is well-defined \ref{pushforwardgeneral}.
We claim that for every $f:(X,\Phi)\xr{} (Y,\Psi)$ in $V_*$ the following equality holds:
\begin{equation}\label{propositionmainthm1claim1}
\phi'_{(Y,\Psi)}\circ f_* = F_*(f)\circ \phi'_{(X,\Phi)}.
\end{equation}

Let $W\in \Phi$ be irreducible. If $\dim(f(W))<\dim(W)$ then 
$F_{2\dim W}(Y,f(W))=0$ by semi-purity \ref{semipurity}, thus \ref{propositionmainthm1claim1}
holds in this case.

In the case $\dim f(W) =\dim W=:d$ the map $W\xr{} f(W)$ is generically finite, so that we may find
an open $U\subset Y$ with $U\cap f(W)\neq \emptyset$, $U\cap f(W)$ is smooth,  
$f^{-1}(U)\cap W$ is smooth, and $f':f^{-1}(U)\cap W\xr{} U\cap f(W)$ induced by $f$ is finite. 
Consider the commutative diagram 
$$
\xymatrix
{ 
F_{2d}(X,W) \ar[d]^{F_*(f)} \ar[r]^-{F^*(\jmath')}
&
F_{2d}(f^{-1}(U),W\cap f^{-1}(U)) \ar[d]^{F_*(f)}
&
F_{2d}(W\cap f^{-1}(U)) \ar[l]_-{F_*(\imath')} \ar[d]^{F_*(f')}
\\
F_{2d}(Y,f(W)) \ar[r]^-{F^*(\jmath)}
&
F_{2d}(U,f(W)\cap U) 
&
F_{2d}(f(W)\cap U) \ar[l]_-{F_*(\imath)},
}
$$
where $\jmath:(U,f(W)\cap U)\xr{} (Y,f(W))$ resp.  $\jmath':(f^{-1}(U),W\cap f^{-1}(U))
\xr{} (X,W)$ in $V^*$ are induced by the obvious open immersions, and 
$\imath: f(W)\cap U \xr{} (U,f(W)\cap U)$ resp. 
$\imath': W\cap f^{-1}(U) \xr{} (f^{-1}(U),W\cap f^{-1}(U))$ in $V_*$ are induced by the obvious
closed immersions. From the diagram and condition (\ref{mainthm1}),(\ref{traceofone}) 
we obtain   
$$
F^*(\jmath)F_*(f)(cl_{(X,W)})= \deg(W/f(W))\cdot F_*(\imath)(1_{f(W)\cap U}).
$$
Now, semi-purity \ref{semipurity} implies 
$$
F_*(f)(cl_{(X,W)})= \deg(W/f(W)) \cdot cl_{(Y,f(W))},
$$
which proves the claim \ref{propositionmainthm1claim1}.

\emph{2nd Step:} Let $X$ be a smooth scheme, $W\subset X$ an irreducible closed subset, and 
$D\subset X$ a smooth divisor intersecting $W$ properly, so that $W\cap D:=W\times_X D$ is an effective 
Cartier divisor on $W$. We denote by $[W\cap D]$ the associated Weil divisor and claim that 
\begin{equation}\label{propositionmainthm1claim2}
F^*(\imath_D)(\phi'_{(X,W)}([W]))= \phi'_{(D,W\cap D)}([W\cap D]),
\end{equation} 
where $\imath_D:(D,W\cap D)\xr{} (X,W)$ is induced by $D\subset X$. 

Note that by semi-purity we may replace $X$ by an open subset which contains the generic
points of $(W\cap D)_{red}$. In particular we may assume that the irreducible components 
of $(W\cap D)_{red}$ are disjoint. Letting $V_1,\dots,V_r$ be the irreducible components 
of $(W\cap D)_{red}$ we obtain 
$$
\bigoplus_{i=1}^r F(D,V_i) \xr{\cong} F(D,W\cap D)
$$ 
from (\ref{triplescond})(\ref{triplescondsum}); thus we may assume that $r=1$. If $W$ is
regular (=smooth) in codimension one (e.g. $W$ is normal) then we can find an open $U\subset X$
such that $W\cap U$ and $V_1\cap U\neq \emptyset$ is smooth; thus \ref{propositionmainthm1claim2} follows from (\ref{mainthm1})(\ref{multiplicities}).

Now, let $W$ be not necessarily normal. Since we may assume that $X$ is affine we can find
a closed immersion $\ti{W}\xr{} W\times \P^n$ (over $W$) of the normalisation $\ti{W}$ of $W$.
Setting 
$$
\tilde{X}:=X\times \P^n,\, \tilde{D}:=D\times \P^n,\, \tilde{\imath}:(\tilde{D},\tilde{V}\cap \tilde{D})\xr{} (\tilde{X},\tilde{V}),
$$
we obtain
\begin{equation*}
\begin{split}
  F^*(\imath)(\phi'_{(X,W)}[W])&= F^*(\imath)F_*({\rm pr}_1)(\phi'_{(\ti{X},\ti{W})}([\ti{W}])) \quad \text{(\ref{propositionmainthm1claim1})} \\
                    &=F_*({\rm pr}_{1\mid \tilde{D}}) F^*(\tilde{\imath}) (\phi'_{(\ti{X},\ti{W})}([\ti{W}])) \quad (\ref{triplescond})(\ref{MP}) \\
                    &=F_*({\rm pr}_{1\mid \tilde{D}})(\phi'_{(\ti{D},\ti{W}\cap \ti{D})}([\tilde{W}\cap \tilde{D}]))  \\
                    &=\phi'_{(D,W\cap D)}({\rm pr}_{1*}([\tilde{W}\cap \tilde{D}]))) \quad \text{(\ref{propositionmainthm1claim1})} \\
                    &=\phi'_{(D,W\cap D)}([W\cap D]). 
\end{split}
\end{equation*} 

\emph{3rd Step:}
For all $(X,\Phi)$ we claim that the map $\phi'_{(X,\Phi)}$ satisfies
\begin{equation}\label{propositionmainthm1claim3}
\phi'_{(X,\Phi)}({\rm Rat}_{\Phi}(X))=0;
\end{equation}
and thus induces a natural transformation
$$
\phi:\CH_* \xr{} F_*.
$$ 

Let $W\subset X\times \P^1$ be irreducible such that ${\rm pr}_1(W)\in \Phi$ and $W\xr{} \P^1$
is dominant. By using the 2nd Step (\ref{propositionmainthm1claim2}) we obtain 
$$
F^*(\imath_{\epsilon})(\phi'_{(X\times \P^1,W)}([W]))=\phi'_{(X,{\rm pr}_1(W))}([W\cap (X\times \{\epsilon\})])
$$
for $\epsilon \in \{0,\infty\}$, $\imath_{\epsilon}:(X\times\{\epsilon\},{\rm pr}_1(W)) \xr{}
(X\times \P^1,W)$.

Thus $F^*(\imath_0)=F^*(\imath_{\infty})$ will prove the claim (\ref{propositionmainthm1claim3}). 
It is not difficult to see that this follows from the projection formula and 
\begin{equation}\label{0=infty}
F_*(\imath'_0)(1_X)=F_*(\imath'_{\infty})(1_X)
\end{equation}
in $F(X\times \P^1)$, where $\imath'_{\epsilon}:X\times \{\epsilon\}\xr{\subset} X\times \P^1$.

In view of (\ref{triplescond})(\ref{MP}) the equality \ref{0=infty} is implied by (\ref{mainthm1})(\ref{rationaleq}).

\emph{4th Step:} The only assertion left to prove is 
$$
\phi \circ \times = T \circ (\phi\otimes \phi), \quad \phi\circ 1 =e.
$$  
The second equality holds by definition. For the first equality it is sufficient 
to show 
$$
\phi'_{(X\times Y,W\times V)}([W]\times [V])=T(\phi'_{(X,W)}([W]) \otimes \phi'_{(Y,V)}([V]))
$$ 
for smooth schemes $X,Y$ and irreducible closed subsets $W\subset X,V\subset Y$. Again by
semi-purity we may assume that $W$ and $V$ are smooth in which case the statement follows
from Lemma \ref{lemmaaboutone}.
\end{proof}
\end{proposition}

\begin{proof}[Proof of Theorem \ref{mainthm1}] \label{mainthm1proof}
Set $\CH:=\Cht$ and $F:=\triples$. 
For $\phi\in \Hom_{\T}(\CH,F)$ we get $$\phi(1_X)=1_X$$ for all smooth 
schemes $X$; thus Proposition \ref{propositionmainthm1} implies that 
$\Hom_{\T}(\CH,F)$ is either empty or contains only one element.

Obviously, the conditions (\ref{mainthm1})(\ref{traceofone})-(\ref{class}) are necessary 
for $\Hom_{\T}(\CH,F)$ to be non-empty. So let us assume that the conditions 
are satisfied. Proposition \ref{propositionmainthm1} yields a natural transformation of
right-lax symmetric monoidal functors 
$$
\phi:(\CH_*,\times,1) \xr{} (F_*,T,e).
$$
We need to prove that $\phi$ induces a natural transformation 
$$
\phi:\CH^*\xr{} F^*.
$$

\emph{1st Step:} Assume that $f:(X,\Phi)\xr{} (Y,\Psi)$ in $V^*$ is smooth. We claim the
commutativity of the following diagram:
$$
\xymatrix
{
\CH(Y,\Psi) \ar[r]^{\CH^*(f)} \ar[d]^{\phi}
&
\CH(X,\Phi) \ar[d]^{\phi}
\\
F(Y,\Psi) \ar[r]^{F^*(f)}
&
F(X,\Phi). 
}
$$ 
It is sufficient to prove 
$$
F^*(f)(\phi_{(Y,V)}([V]))=\phi_{(X,f^{-1}(V))}(f^*[V])
$$
for all irreducible closed subsets $V\subset Y$. By using semi-purity we may replace
$Y$ by an open set and thus assume that $V$ is smooth. We obtain
\begin{equation*}
\begin{split}
F^*(f)(\phi_{(Y,V)}([V]))&= F^*(f)F_*(\imath_V)(1_V) \\
                       &= F_*(\imath_{f^{-1}(V)})F^*(f_{\mid f^{-1}(V)})(1_V) \quad (\ref{triplescond})(\ref{MP})\\
                       &= F_*(\imath_{f^{-1}(V)})(1_{f^{-1}(V)})\\
                       &= \phi_{(X,f^{-1}(V))}([f^{-1}(V)]),
\end{split}
\end{equation*} 
where $\imath_V:V\xr{} (Y,V), \imath_{f^{-1}(V)}:f^{-1}(V)\xr{} (X,f^{-1}(V))$.

\emph{2nd Step:} Let $p:E\xr{}X$ be a vector bundle and let $s:X\xr{} E$ be the zero section.
We claim that for every closed subscheme $W\subset X$ the following diagram is commutative:
$$
\xymatrix
{
\CH(E,p^{-1}(W)) \ar[d]^{\phi} \ar[r]^-{\CH^*(s)}
&
\CH(X,W) \ar[d]^{\phi}
\\
F(E,p^{-1}(W)) \ar[r]^-{F^*(s)}
&
F(X,W).
}
$$
Indeed, by homotopy invariance we may write any  $a\in \CH(p^{-1}(W))$ as $a=\CH^*(p)(b)$ with
$b\in \CH(W)$. Thus by the 1st Step:
$$F^*(s)(\phi(a))= F^*(s)F^*(p)(\phi(b))=\phi(b)=\phi(\CH^*(s)(a)).$$ 

\emph{3rd Step:} For every closed subscheme $W\subset X$ denote by 
$
\imath_0:(X,W)\xr{} (X\times \P^1,W\times \P^1),
$
resp.~$
\imath_{\infty}:(X,W)\xr{} (X\times \P^1,W\times \P^1),
$  
the morphisms in $V^*$ and $V_*$ induced by the inclusion $X\times\{0\}\subset X\times \P^1$,
resp.~$X\times\{\infty\}\subset X\times \P^1$. We claim that 
$$
F^*(\imath_0)=F^*(\imath_{\infty}).
$$ 
Indeed, if $p:(X\times \P^1, W\times \P^1)\xr{} (X, W)$ is the first projection, then  
$$
F^*(\imath_{\epsilon})(a)=F_*(p)F_*(\imath_\epsilon)(\phi([X])\cup F^*(\imath_\epsilon)(a))=F_*(p)(F_*(\imath_{\epsilon})(\phi([X]))\cup a)
$$
for $\epsilon\in \{0,\infty\}$. Since 
$F_*(\imath_{\epsilon})(\phi([X]))=\phi([X\times \{\epsilon\}])$ the claim follows from 
$[X\times \{0\}]=[X\times \{\infty\}]$ in $\CH^1(X\times \P^1)$.

\emph{4th Step:}  Let $f:X\xr{} Y$ be a closed immersion and $V\subset Y$ a closed 
subscheme; set $W:=f^{-1}(V)=V\times_Y X$. Then $f$ induces $f:(X,W)\xr{} (Y,V)$
in $V^*$ and we claim that  
$$
\xymatrix
{
\CH(Y,V) \ar[r]^{\CH^*(f)} \ar[d]^{\phi}
&
\CH(X,W) \ar[d]^{\phi}\\
F(Y,V) \ar[r]^{F^*(f)} 
&
F(X,W)
}
$$ 
is a commutative diagram.

Again, it is sufficient to prove 
$$F^*(f)(\phi([V]))=\phi(\CH^*(f))$$
for $V$ integral.

For the proof we use deformation to the normal cone (see \cite[\textsection5]{F}). Let 
\begin{equation*}
\begin{split}
M^0&:={\rm Bl}_{X\times \{\infty\}}(Y\times \P^1) \backslash {\rm Bl}_{X\times \{\infty\}}(Y\times \{\infty\}), \\ 
\tilde{M}^0&:={\rm Bl}_{W\times \{\infty\}}(V\times \P^1) \backslash {\rm Bl}_{W\times \{\infty\}}(V\times \{\infty\}),
\end{split}
\end{equation*}
then $\tilde{M}^0\subset M^0$ is closed, and $M^0, \tilde{M}^0$ are flat over $\P^1$. 
We have a closed immersion 
$\imath_X:X\times \P^1\xr{} M^0$ resp. $\imath_W:W\times \P^1 \xr{} \tilde{M}^0$ 
which deforms the immersion $X\subset Y$ resp. $W\subset V$ over $\P^1\backslash \{\infty\}$
to the zero section of the normal cone over $\infty$.

Since $W\times \P^1=\tilde{M}^0\cap (X\times \P^1)$ we obtain morphisms 
$$
\imath_{\epsilon}:(X\times\{\epsilon\},W\times \{\epsilon\}) \xr{} (M^0,\tilde{M}^0) 
$$
in $V^*$ for $\epsilon\in \{0,\infty\}$. By the 3rd Step we know that  
$F^*(\imath_0)=F^*(\imath_{\infty}).$

Consider the projection  
$
p:(Y\times (\P^1\backslash \{\infty\}),V\times (\P^1\backslash \{\infty\})) 
\xr{} (Y,V)$ 
in  $V^*$. Note that  $\tilde{M}^0$ is the closure of $V\times(\P^1\setminus\{\infty\})$ in $M^0$, and thus 
$$\CH^*(p)([V])=\CH^*(\jmath)([\tilde{M}^0])$$  
with 
$\jmath:(Y\times (\P^1\backslash \{\infty\}),V\times (\P^1\backslash \{\infty\}))\xr{} (M^0,\tilde{M}^0)$ 
the open immersion. By using the 1st Step we get
$
F^*(p)(\phi([V]))=F^*(\jmath)(\phi([\tilde{M}^0]))
$
and thus
$$
F^*(f)(\phi([V]))=F^*(\imath_0)(\phi([\tilde{M}^0]))=F^*(\imath_\infty)(\phi([\tilde{M}^0])).
$$
 
Now, let us compute $F^*(\imath_{\infty})$. The morphism $\imath_{\infty}$ has a factorization
$$
\imath_{\infty}: (X,W) \xr{s} (N_{Y/X},C_{V/W}) \xr{t} (M^0,\tilde{M}^0),
$$ 
where $N_{Y/X}$ is the normal bundle and $C_{V/W}$ is the normal cone. Note that 
$N_{Y/X}$ is a smooth divisor in $M^0$, which intersects $\tilde{M}^0$ properly (being the fiber of $M^0\xr{} \P^{1}$ over $\infty$), so that we may apply 
(\ref{propositionmainthm1claim2}) to $t$. Moreover $s$ is the zero section of the normal bundle. The zero section
also induces a morphism
$$
s': (X,W) \xr{} (N_{Y/X},N_{Y/X}\times_X W)
$$
in $V^*$; denote by
$\tau: (N_{Y/X},C_{V/W}) \xr{} (N_{Y/X},N_{Y/X}\times_X W)$ the morphism in $V_*$ induced
by the identity map. Then \ref{triplescond}(\ref{MP}) yields
$$
F^*(s)=F^*(s')\circ F_*(\tau).
$$
Thus we get  
\begin{equation*}
\begin{split}
F^*(\imath_{\infty})(\phi([\tilde{M}^0]))    &=F^*(s')F_*(\tau)F^*(t)(\phi([\tilde{M}^0])) \\
                                   &=F^*(s')F_*(\tau)(\phi(\CH^*(t)([\tilde{M}^0]))) \quad \text{(\ref{propositionmainthm1claim2})}\\
                                   &=\phi(\CH^*(\imath_{\infty})([\tilde{M}^0])) \quad (\text{2nd Step})\\
                                   &=\phi(\CH^*(\imath_{0})([\tilde{M}^0]))=\phi(\CH^*(f)([V])).
\end{split}
\end{equation*}

\emph{5th Step:} Let $f:(X,\Phi)\xr{} (Y,\Psi)$ be any morphism in $V^*$. We have to prove 
that 
$$\phi\circ \CH^*(f)=F^*(f)\circ \phi.$$ 
Indeed, $f$ factors through 
$$
f:(X,\Phi)\xr{(id,f)} (X\times Y,{\rm pr}_2^{-1}(\Psi)) \xr{{\rm pr}_2} (Y,\Psi).
$$
By the 1st Step we may reduce to the case of the closed immersion $(id,f)$, and 
by using Proposition \ref{indlimit} the statement follows from the 4th Step.
\end{proof}

\subsection{Correspondences} \label{notation:correspondences}
Let $\triples\in \T$.
Let $X_i,i=1,2,3,$ be smooth varieties and $\Phi_{ij}$, for $ij=12,23,13,$ be families
of supports on $X_i\times X_j$. Denote by $p_{ij}:X_1\times X_2\times X_3\xr{} X_i\times X_j$ the 
projection. Suppose that 
\begin{equation}\label{conditioncomposition}
\begin{cases} \text{$p_{13} \mid p_{12}^{-1}(\Phi_{12})\cap p_{23}^{-1}(\Phi_{23})$ is proper,} 
\\
p_{13}(p_{12}^{-1}(\Phi_{12})\cap p_{23}^{-1}(\Phi_{23}))\subset \Phi_{13}. \end{cases}
\end{equation}
Then we  define
$$ F(X_1\times X_2,\Phi_{12})\otimes F(X_2\times X_3,{\Phi_{23}})\xr{} F(X_1\times X_3,\Phi_{13}), \quad a\otimes b\mapsto b\circ a$$
to be the composition 
\begin{multline}\label{equation-def-circ}
F(X_1\times X_2,\Phi_{12})\otimes F(X_2\times X_3,{\Phi_{23}}) 
\xr{F^*(p_{12})\otimes F^*(p_{23})}\\
F(X_1\times X_2\times X_3, p_{12}^{-1}(\Phi_{12}))\otimes 
F(X_1\times X_2\times X_3,p_{23}^{-1}(\Phi_{23})) \xr{\cup} \\
 F(X_1\times X_2\times X_3,p_{12}^{-1}(\Phi_{12})\cap p_{23}^{-1}(\Phi_{23})) \xr{F_*(p_{13})}
 F(X_1\times X_3,\Phi_{13}).
\end{multline}

\subsubsection{}
Let $\Phi_{ij}'$ for $ij=12,23,13,$ be families
of supports on $X_i\times X_j$. Suppose that 
$$
\begin{cases} \text{$p_{13} \mid p_{12}^{-1}(\Phi_{12}')\cap p_{23}^{-1}(\Phi_{23}')$ is proper,} 
\\
p_{13}(p_{12}^{-1}(\Phi_{12}')\cap p_{23}^{-1}(\Phi_{23}'))\subset \Phi_{13}', \end{cases}
$$
and $\Phi_{ij}'\subset \Phi_{ij}$ for $ij=12,23,13$.
Obviously, the following diagram is commutative 
$$
\xymatrix
{
F(X_1\times X_2,{\Phi_{12}'})\otimes F(X_2\times X_3,\Phi_{23}') \ar[d] \ar[r]^-{\circ} 
&
F(X_1\times X_3,\Phi_{13}')\ar[d]
\\
F(X_1\times X_2,{\Phi_{12}})\otimes F(X_2\times X_3,\Phi_{23}) \ar[r]^-{\circ}
&
F(X_1\times X_3,{\Phi_{13}}).
}
$$

The most important case for us will be $(\CH_*,\CH^*,\times,1)$. For later use we record the following particular case of the above discussion.

\begin{lemma}\label{remark:support}
Let $X_i,i=1,2,3,$ be smooth schemes and $\Phi_{ij}$, for $ij=12,23,13,$ be families of supports on $X_i\times X_j$, which satisfy \eqref{conditioncomposition}.
For $a\in Z_{\Phi_{12}}(X_1\times X_2)$ and $b\in Z_{\Phi_{23}}(X_2\times X_3)$ we define
\begin{equation}\label{equation-def-supp(a,b)}
\supp(a,b):=p_{13}(p_{12}^{-1}(\supp(a))\cap p_{23}^{-1}(\supp(b))),
\end{equation}
which is a closed subset contained in $\Phi_{13}$. 
The families of supports 
$$
\Phi_{12}'=\Phi_{\supp(a)}, \quad \Phi_{23}'=\Phi_{\supp(b)}, \quad 
\Phi_{13}'=\Phi_{\supp(a,b)}
$$ 
satisfy \eqref{conditioncomposition}. 
The cycles $a, b$ define in the obvious way classes 
\begin{align*}
\tilde{a}\in \CH(\supp(a)), \quad &\tilde{b}\in \CH(\supp(b))\\
a\in \CH(X_1\times X_2,\Phi_{12}), \quad &b\in \CH(X_2\times X_3,\Phi_{23}).
\end{align*} 
Then $b\circ a$ is the image of $\tilde{b}\circ \tilde{a}$ via the map 
$$
\CH(\supp(a,b))\xr{} \CH(X_1\times X_3,\Phi_{13}).
$$ 
\end{lemma}

The Lemma \ref{remark:support} helps to understand the composition 
of two cycles $a,b$ via the purely set-theoretic computation of 
$\supp(a,b)$. Frequently we are able to compute the composition 
over suitable good open subsets; this is the motivation for the
next Lemma.

\begin{lemma}\label{lemma-composition-localization}
Let $X_i,i=1,2,3,$ be smooth schemes. Let $a\in Z(X_1\times X_2), b\in Z(X_2\times X_3)$ be algebraic cycles such that
$$
p_{13}\mid p^{-1}_{12}\supp(a)\cap p^{-1}_{23}\supp(b) \quad \text{is proper.}
$$
Let $X'_1\subset X_1, X'_3\subset X_3$ be open subsets; define $a'\in Z(X'_1\times X_2), b'\in Z(X_2\times X'_3)$ as the restrictions of $a,b$. We denote
by $p'_{ij}$ the projections from $X'_1\times X_2\times X'_3$.
\begin{enumerate}
\item The restriction of $p'_{13}$ to $p'^{-1}_{12}\supp(a')\cap p'^{-1}_{23}\supp(b')$ is proper.
\item The equality 
$$
\supp(a',b')=\supp(a,b)\cap (X'_1\times X'_3)
$$
holds, where $\supp(a,b)$ is defined in \eqref{equation-def-supp(a,b)}.
\item The composition $b'\circ a'$ is the image of $b\circ a$ via the localization 
map
$$
\CH(\supp(a,b))\xr{} \CH(\supp(a',b')).
$$
(Here $\supp(a',b')\subset \supp(a,b)$ is an open subset by (2)).
\end{enumerate}
\begin{proof}
By definition we obtain
$$
\supp(a')=\supp(a)\cap (X'_1\times X_2), \quad
\supp(b')=\supp(b)\cap (X_2\times X'_3).$$ 

For (1). Let $Z_{12}\subset X_1\times X_2, Z_{23}\subset X_2\times X_3$ be 
closed subsets such that 
$$
p_{13}\mid p^{-1}_{12}Z_{12}\cap p^{-1}_{23}Z_{23} \quad \text{is proper.}
$$
Set $Z_{12}'=Z_{12}\cap (X'_1\times X_2), Z'_{23}=Z_{23}\cap (X_2\times X'_{3})$.
Obviously, 
$$
p'^{-1}_{12}Z'_{12}\cap p'^{-1}_{23}Z'_{23}=(p^{-1}_{12}Z_{12}\cap p^{-1}_{23}Z_{23})\cap (X'_1\times X_2\times X'_3).
$$
Thus, if $p^{-1}_{12}Z_{12}\cap p^{-1}_{23}Z_{23}$ is proper
over $X_1\times X_3$ then $p'^{-1}_{12}Z'_{12}\cap p'^{-1}_{23}Z'_{23}$ is 
proper over $X'_1\times X'_{3}$.

Statement (2) is a straightforward computation. For (3): By using 
the definition of $\circ$ in \eqref{equation-def-circ} it is 
straightforward to show that the diagram 
$$
\xymatrix
{ 
\CH(X_1\times X_2,\supp(a))\otimes \CH(X_2\times X_3,\supp(b))\ar[d]\ar[r]^-{\circ} 
&
\CH(X_1\times X_3,\supp(a,b))\ar[d]\\
\CH(X'_1\times X_2,\supp(a'))\otimes \CH(X_2\times X'_3,\supp(b'))\ar[r]^-{\circ} 
&
\CH(X'_1\times X'_3,\supp(a',b'))\\
}
$$
is commutative.
\end{proof}
\end{lemma}

\subsubsection{} For two smooth schemes $X,Y$ and a family $\Phi$ (resp.~$\Psi$) 
of supports of $X$ (resp.~$Y$) we define a family of supports $P(\Phi,\Psi)$
on the product by 
\begin{multline}\label{definitionP}
P(\Phi,\Psi):=\{Z\subset X\times Y; \text{$Z$ is closed, ${\rm pr}_2 \mid Z$ is proper,} \\ 
\text{$Z\cap {\rm pr}_1^{-1}(W)\in {\rm pr}_2^{-1}(\Psi)$ for every $W\in \Phi$.} \}
\end{multline}

Let $X_i,i=1,2,3,$ be smooth schemes and let $\Phi_{i}$ be a family
of supports on $X_i$ for $i=1,2,3$. It is easy to see that $\Phi_{ij}:=P(\Phi_i,\Phi_j)$
satisfy the condition \ref{conditioncomposition} and therefore  
\begin{equation} \label{compositioncor}
F(X_1\times X_2,P(\Phi_1,\Phi_2))\otimes F(X_2\times X_3,P(\Phi_2,\Phi_3)) \xr{} 
F(X_1\times X_3,P(\Phi_1,\Phi_3)),  
\end{equation}
where $a\otimes b\mapsto b\circ a$, is well-defined.

\begin{proposition}
\begin{enumerate}
\item Let $X_i,i=1,\dots,4,$ be a smooth scheme and let $\Phi_{i},i=1,\dots,4,$ be a family 
of supports of $X_i$.
For all $a_{ij}\in F(X_i\times X_j,P(\Phi_i,\Phi_j))$ the following equality holds
$$
a_{34}\circ (a_{23}\circ a_{12})=(a_{34}\circ a_{23})\circ a_{12}.
$$
\item For any $(X,\Phi)$ the diagonal immersion induces a morphism 
$\imath: X\xr{} (X\times X,P(\Phi,\Phi))$ in $V_*$. We set 
$$
\Delta_{(X,\Phi)}:=F_*(\imath)(1_X).
$$
The equality $\Delta_{(X,\Phi)} \circ g=g$ resp.  $g \circ \Delta_{(X,\Phi)}=g$ holds for all 
$(Y,\Psi)$ and 
$g\in F(Y\times X,P(\Psi,\Phi))$ resp. $g\in  F(X\times Y,P(\Phi,\Psi))$.   
\end{enumerate}
\begin{proof}
The proof of the first statement is as in \cite[Proposition~16.1.1]{F} but one has to keep track of the supports, which is straightforward.

The second statement is an easy computation.
\end{proof}
\end{proposition}

\subsubsection{Grading}
For $(X,\Phi)$ and $(Y,\Psi)$ there are two different gradings on
$
F(X\times Y,P(\Phi,\Psi)),
$
coming from $F_*$ and $F^*$. Unfortunately, neither are compatible with $\circ$ 
from (\ref{compositioncor}). We define a new grading by 
$$
F(X\times Y,P(\Phi,\Psi))^i=\bigoplus_{X'} F^{2\dim(X')+i}(X'\times Y,P(\Phi,\Psi)),
$$
where $X'$ runs through the connected components of $X$. With this grading $\circ$
becomes a morphism of graded abelian groups. 

By the definition of the grading there are choices. We could also define a grading   
$$
F(X\times Y,P(\Phi,\Psi))_i=\bigoplus_{X'} F_{2\dim(X')+i}(X'\times Y,P(\Phi,\Psi)).
$$

\begin{definition}
To an object $F=\triples\in \T$ we attach the graded additive symmetric monoidal  category
$Cor_F$ with objects 
${\rm ob}(Cor_F)={\rm ob}(V_*)={\rm ob}(V^*)$ and morphisms 
$$
\Hom_{Cor_F}((X,\Phi),(Y,\Psi))=F(X\times Y,P(\Phi,\Psi))
$$
with composition law $a\otimes b\mapsto{} b\circ a$ (\ref{compositioncor}). The identity
is $\Delta_{(X,\Phi)}$. 

The product $\otimes$ on $Cor_F$ is defined by 
$$
(X,\Phi)\otimes (Y,\Psi) :=(X\times Y,\Phi\times \Psi),
$$ 
and for two morphisms $f\in F(X\times X',P(\Phi,\Phi')), g\in F(Y\times Y',P(\Psi,\Psi')),$
we define 
\begin{equation*}
\begin{split}
f\otimes g&\in \Hom_{Cor_F}((X,\Phi)\otimes (Y,\Psi),(X',\Phi')\otimes (Y',\Psi'))\\
f\otimes g&:=F_*(id_X\times \mu_{X',Y}\times id_{Y'})(T(f\otimes g)),
\end{split}
\end{equation*}
where $\mu_{X',Y}$ is the permutation of the factors $(X',\Phi'),(Y,\Psi)$. 
\end{definition}

\subsubsection{}
Given two objects $F,G\in \T$ and a morphism $\phi:F\xr{} G$ in $\T$ we obtain 
a functor of graded additive symmetric monoidal categories
$$
Cor(\phi):Cor_F \xr{} Cor_G
$$
which is given by 
$$
\phi: F(X\times Y,P(\Phi,\Psi)) \xr{} G(X\times Y,P(\Phi,\Psi))
$$
for all $(X,\Phi), (Y,\Psi)$. This provides a functor 
$$
Cor: \T \xr{} \CatGrAb, \quad F \mapsto Cor_F, \quad \phi\mapsto Cor(\phi).
$$
Here, $\CatGrAb$ is the category of graded additive symmetric monoidal categories.

\subsubsection{} 
In order to state the properties of $Cor$ it is convenient to introduce the 
category $V$ with objects ${\rm ob}(V)={\rm ob}(V_*)={\rm ob}(V^*)$ and only
morphisms the identity $id_X$ (for every $X\in {\rm ob}(V)$). There are 
obvious functors $V\xr{} V_*, V\xr{} V^*,$ and $V\xr{} Cor_F$
for all $F\in \T$. We define $\CatGrAbV$ to be the category with functors 
$V\xr{} X$ as objects ($X\in \CatGrAb$) and commutative diagrams
\begin{equation}\label{underV}
\xymatrix
{
X\ar[rr]^{f}
&
&
Y,
\\
&
V \ar[ur]\ar[ul]
&
}
\end{equation}
with $f\in \Hom_{\CatGrAb}(X,Y)$, as morphisms. In general a functor 
$f$ is called \emph{under} $V$ if the diagram \ref{underV} is commutative.

\begin{proposition}\label{propositionTtoCor}
The functor $Cor:\T\xr{} \CatGrAbV$ is fully faithful. 
\begin{proof}
Given $F,G\in \T$ and $\phi:F\xr{} G$ we can recover 
$\phi:F(X)\xr{} G(X)$, for $X\in {\rm ob}(V)$, from the map $Cor(\phi)$:
\begin{equation}
\label{CortoT}
\Hom_{Cor_F}(\Spec(k),X) \xr{} \Hom_{Cor_G}(\Spec(k),X). 
\end{equation}
On the other hand given $\psi:Cor_F\xr{} Cor_G$ in $\CatGrAbV$ then \ref{CortoT}
defines a morphism $F\xr{} G$ in $\T$. 
\end{proof}
\end{proposition}

\subsubsection{}\label{VtoCor}
For all $F\in \T$ there is a functor 
$$\rho_F:Cor_F\xr{} \GrAb$$ 
defined by 
\begin{equation*}
\begin{split}
  \quad &\rho_F(X,\Phi)=F(X,\Phi) \\
&\rho_F(\gamma)=(a\mapsto F_*({\rm pr}_2)(F^*({\rm pr}_1)(a)\cup \gamma)) \quad \text{for $\gamma\in F(X\times Y,P(\Phi,\Psi))$.}
\end{split}
\end{equation*}
The map $\rho_F(\gamma):F(X,\Phi)\xr{} F(Y,\Psi)$ is well-defined since 
$
{\rm pr}_2 \mid {\rm pr}_1^{-1}(\Phi)\cap P(\Phi,\Psi) 
$
is proper and ${\rm pr}_1^{-1}(\Phi)\cap P(\Phi,\Psi) \subset {\rm pr}_2^{-1}(\Psi)$
by definition of $P(\Phi,\Psi)$. Functoriality is again a straightforward computation.

Moreover, there are functors 
$$\tau_*^F:V_*\xr{} Cor_F, \quad \tau^*_F:(V^*)^{op}\xr{} Cor_F,$$ 
(under $V$) such that 
$$
\rho_F \circ \tau^F_*=F_*, \quad \rho_F \circ \tau^*_F=F^*.
$$ 
The functor $\tau_*^F:V_*\xr{} Cor_F$ is defined by mapping a morphism $f:(X,\Phi)\xr{} (Y,\Psi)$ 
to $F_*(id,f)(1_X)$ where $(id,f):X\xr{} (X\times Y,P(\Phi,\Psi))$ is in $V_*$. Similarly, the
functor $\tau^*_F:V^*\xr{} Cor_F$ is defined by mapping a morphism $f:(X,\Phi)\xr{} (Y,\Psi)$
to $F_*(f,id)(1_X)$ with  $(f,id):X\xr{} (Y\times X,P(\Psi,\Phi))$ in $V_*$. 
The equalities $\rho_F \circ \tau^F_*=F_*$ and $\rho_F \circ \tau^*_F=F^*$ follow easily
from the projection formula.

\begin{lemma}\label{lemmaCortau}
If $\phi:F\xr{} G$ is a morphism in $\T$ then 
\begin{equation*}
\begin{split}
Cor(\phi)\circ \tau_*^F&=\tau_*^G, \\
Cor(\phi)\circ \tau^*_F&=\tau^*_G.
\end{split}
\end{equation*}
\begin{proof}
For the first equality let $f:(X,\Phi)\xr{} (Y,\Psi)$ be a morphism in $V_*$. 
We get 
\begin{multline*}
Cor(\phi)(\tau^F_*(f)) = Cor(\phi)(F_*(id,f)(1_X)) =\phi(F_*(id,f)(1_X))=\\
G_*(id,f)(\phi(1_X))=G_*(id,f)(1_X)=\tau^G_*(f).
\end{multline*}
The second equality is proved in the same way.
\end{proof}
\end{lemma}

\section{Hodge cohomology with support}
For a smooth scheme $X$ and a family of supports $\Phi$ of $X$, we define 
$$
H(X,\Phi):=\bigoplus_{i,j} H^i_{\Phi}(X,\Omega_X^j),
$$
and call this $k$-vector space the Hodge cohomology of $X$ with support in $\Phi$.
We denote by $H^*(X,\Phi)$ the graded abelian group, which in degree $n$ equals
\eq{Hodgeuppergrading}{H^n(X,\Phi)= \bigoplus_{i+j=n} H^i_\Phi(X, \Omega^j_X).}
We denote by $H_*(X,\Phi)$ the graded abelian group, which in degree $n$ equals
\eq{Hodgelowergrading}{H_n(X,\Phi)=\bigoplus_r H^{2\dim X_r-n}(X_r, \Phi),}
where $X=\coprod_r X_r$ is the decomposition into connected components. 
We define 
\eq{e}{e:  \Z\to H(\Spec\, k)=k}
to be the natural map sending 1 to 1.

The goal of this section is to provide the object functions $H_*$ and $H^*$ with the structure of functors 
$$H_*:V_*\xr{} \GrAb,\quad H^*:(V^*)^{op} \xr{} \GrAb$$ 
and to define for each $(X,\Phi), (Y,\Psi)\in{\rm ob}(V_*)={\rm ob}(V^*)$ a morphism of graded abelian groups (for both gradings)   
$$
T_{(X,\Phi), (Y,\Psi)}: H(X,\Phi)\otimes H(Y,\Psi)\xr{} H(X\times Y,\Phi\times \Psi)
$$
such that
$
(H_*,H^*,T, e)
$
is an object in $\T$, i.e. it is a datum as in \ref{triples} and satisfies the properties \ref{triplescond}.

\subsection{Pullback}
\subsubsection{}
We work in the bounded derived category of quasi-coherent sheaves $D^b(X)$ on a scheme $X$. 
(The bounded derived category of coherent sheaves will be denoted by $D^b_c(X)$.)
Let $f:X\xr{} Y$ be a morphism of schemes, $\Phi$ resp. $\Psi$ a family of supports of $X$
resp. $Y$. There is an isomorphism of functors 
\begin{equation}\label{pfsupport}
R\Gamma_{f^{-1}(\Psi)} \xr{\cong} R\Gamma_{\Psi}Rf_*.
\end{equation}
If $\Psi\subset \Psi'$ for another family of supports $\Psi'$ then the diagram 
\begin{equation}\label{refinementpfsupport}
\xymatrix
{
R\Gamma_{f^{-1}(\Psi')} \ar[r]^{\text{\ref{pfsupport}}}
&
R\Gamma_{\Psi'}Rf_*
\\
R\Gamma_{f^{-1}(\Psi)} \ar[r]^{\text{\ref{pfsupport}}} \ar[u]
&
R\Gamma_{\Psi}Rf_* \ar[u] 
}
\end{equation}
is commutative.
Moreover, if $g:Z\xr{} X$ is another morphism of schemes
then the following diagram is commutative 
\begin{equation}\label{compositionpfsupport}
\xymatrix
{
R\Gamma_{(f\circ g)^{-1}(\Psi)} \ar[rr]^{\text{\ref{pfsupport} for $f\circ g$}}
\ar[d]_{\text{\ref{pfsupport} for $g$}}
&
&
R\Gamma_{\Psi}R(f\circ g)_*
\\
R\Gamma_{f^{-1}(\Psi)}Rg_* \ar[urr]_{\text{\ref{pfsupport} for $f$}}.
}
\end{equation}

\subsubsection{}
For a morphism $f:X\xr{} Y$ of schemes we have 
$$
id \xr{} Rf_*Lf^*,
$$
and thus we obtain a morphism of functors
\begin{equation}\label{pbsupport}
R\Gamma_{\Psi}\xr{} R\Gamma_{f^{-1}(\Psi)} Lf^*;
\end{equation}
it easily follows from \ref{refinementpfsupport} that the diagram
\begin{equation}\label{refinementpbsupport}
\xymatrix
{
R\Gamma_{\Psi'} \ar[r]^-{\text{\ref{pbsupport}}}
&
R\Gamma_{f^{-1}(\Psi')}Lf^*
\\
R\Gamma_{\Psi} \ar[r]^-{\text{\ref{pbsupport}}} \ar[u]
&
R\Gamma_{f^{-1}(\Psi)}Lf^* \ar[u] 
}
\end{equation}
commutes for $\Psi\subset \Psi'$.
From \ref{refinementpbsupport} and \ref{compositionpfsupport} 
it follows that for another morphism $g:Z\xr{} X$ of schemes the
following diagram is commutative  
\begin{equation}\label{compositionpbsupport}
\xymatrix
{
R\Gamma_{\Psi} \ar[r]^{\text{\ref{pbsupport} for $f$}} \ar[d]_{\text{\ref{pbsupport} for $f\circ g$}}
&
R\Gamma_{f^{-1}(\Psi)} Lf^* \ar[dl]^{\text{\ref{pbsupport} for $g$}} 
\\
R\Gamma_{(f\circ g)^{-1}(\Psi)} L(f\circ g)^*.
}
\end{equation}

For a morphism  $f:(X,\Phi)\xr{} (Y,\Psi)$ in $V^*$ (i.e. $f^{-1}(\Psi)\subset \Phi$) 
the morphism 
\begin{multline*}
R\Gamma_{\Psi}\Omega_Y^d\xr{} R\Gamma_{f^{-1}(\Psi)} Lf^*\Omega_Y^d =  R\Gamma_{f^{-1}(\Psi)} f^*\Omega_Y^d\xr{}  
                                                   R\Gamma_{f^{-1}(\Psi)} \Omega_X^d \xr{}  R\Gamma_{\Phi} \Omega_X^d 
\end{multline*}
(for $d\geq 0$) gives a morphism 
\begin{equation}\label{pullback}
H^*(f):H(Y,\Psi)\xr{} H(X,\Phi).
\end{equation}
It is a straightforward computation that 
$f\mapsto H^*(f)$ defines a functor $(V^*)^{op}\xr{} \GrAb$.

\subsection{Push-forward in the derived category}
We recall the following notations from duality theory (see \cite{Ha}, \cite{Co}): 
Let $X$ be a separated $k$-scheme of finite type with structure map $\pi: X\to \Spec\, k$. 
We have $\pi^!k\in D^b_c(X)$. (In fact if $X$ has dimension $d$ then $\pi^!k$ has non-zero cohomology only in the interval $[-d, 0]$.
This follows from \cite[V, Prop.7.3 and its proof]{Ha}  and \cite[(3.1.25)]{Co}.) 
We denote 
\eq{dual}{D_X:=R\sHom_X(-, \pi^!k): D^b_c(X)\to D^b_c(X).}
If $f: X\to Y$ is  a {\em proper} morphism between $k$-schemes, we have the trace map
\eq{trace}{\Tr_f : Rf_* f^!\xr{} \id,}
which is a natural transformation of functors on $D^+_c(Y)$. 
For maps $f: X\to Y$ and $g: Y\to Z$, we have the canonical isomorphisms
\eq{c}{c_{f,g}: (gf)^!\xr{\simeq} f^!g^! \quad \text{in } D^+_c(X).}


\begin{notation}\label{dual-notation}
Given a bounded complex $C$ in $D(X)$ and a morphism of complexes $\varphi : A\to B$ in $D(X)$ we will often denote
the morphism $R\sHom_X(C,\varphi): R\sHom_X(C, A)\to R\sHom_X(C,B)$ simply by $\varphi$ and
 the morphism $R\sHom_X(\varphi, C): R\sHom_X(B,C)\to R\sHom_X(A,C)$ by $\varphi^\vee$. It will always be clear from the context what $C$ is in the particular situation. 
\end{notation}

\begin{definition}\label{defpush-forward}
 Let $f: X\to Y$ be a proper $k$-morphism. Let $\pi_X$ and $\pi_Y$ denote the structure maps of $X$ and $Y$
 respectively. Then we define
\[f_*: Rf_*D_X(\Omega^q_X)\to D_Y(\Omega^q_Y),\quad q\ge 0,\] 
to be the composition
\mlnl{Rf_*R\sHom_X(\Omega^q_X,\pi_X^!k)\xr{c_{f,\pi_Y}} Rf_*R\sHom_X(\Omega^q_X, f^!\pi_Y^!k)\\ 
        \xr{nat.} R\sHom_Y(Rf_*\Omega^q_X, Rf_*f^!\pi_Y^!k) \xr{\Tr_f} R\sHom_Y(Rf_*\Omega^q_X, \pi_Y^!k)  \\  
                             \xr{(f^*)^\vee} R\sHom_Y(\Omega^q_Y, \pi_Y^!k).}
\end{definition}

\begin{remark}\label{alt-def-pf}
\begin{enumerate}
 \item Notice that the composition of the middle two arrows in the above composition is just the standard Grothendieck duality isomorphism (see e.g. \cite[(3.4.10)]{Co})
                  $ Rf_*R\sHom_X(-, f^!(-))\xr{\simeq} R\sHom_Y(Rf_*(-), -).$
\item It is straightforward to check that the above push-forward also equals the following composition:
       \mlnl{Rf_*D_X(\Omega^q_X)\xr{({^af^*})^\vee} Rf_*D_X(Lf^*\Omega^q_Y)  \xr{c_{f,\pi_Y}}
               Rf_*R\sHom_X(Lf^*\Omega^q_Y, f^!\pi_Y^!k) \\ \xr{{\rm adj.}} R\sHom_Y(\Omega^q_Y, Rf_*f^!\pi_Y^!k) \xr{\Tr_f} D_Y(\Omega^q_Y).}
Here ``adj.'' denotes the isomorphism (see \cite[II, Prop. 5.10]{Ha})
\[Rf_*R\sHom_X(Lf^*(-),-)\cong R\sHom_Y(-,Rf_*(-)) \quad \text{on } D^-_{\rm c}(Y)\times D^+(X)\] 
and ${^af^*}: Lf^*\Omega^q_Y\to\Omega^q_X$ is the morphism corresponding to $\Omega^q_Y\to Rf_*\Omega^q_X$ under $H^0(Y,-)$ applied to the above isomorphism.
\end{enumerate} 
\end{remark}

\begin{proposition}\label{properties-pf}
 \begin{enumerate}
  \item $\id_*=\id$.
   \item Let $f: X\to Y$ and $g: Y\to Z$ be two proper maps with $X$ and $Y$ of pure dimension $d_X$ and $d_Y$ respectively. Then 
            \[(g\circ f)_* = g_*\circ Rg_*(f_*): Rg_*Rf_*D_X(\Omega^q_X)\to D_Z(\Omega^q_Z).\]
   \item Let  
           \[\xymatrix{X' \ar[r]^{u'}\ar[d]_{f'} & X\ar[d]^f \\ Y'\ar[r]^u & Y }\]
          be a cartesian diagram with $f$ proper, $u$ \'etale and $X$ of pure dimension $d$, then the following diagram commutes
         \eq{etalebasechange}{\xymatrix{u^*Rf_*D_X(\Omega^q_X)\ar[d]_\simeq\ar[r]^{u^*(f_*)} &   u^*D_Y(\Omega^q_Y)\ar[d]^\simeq \\    
                     R{f'}_*D_X(\Omega^q_{X'})\ar[r]^{f'_*} & D_{Y'}(\Omega^q_{Y'}), }}
          where the vertical maps are the natural isomorphisms (in the proof we will make these isomorphisms precise).
\end{enumerate}
\end{proposition}

\begin{proof}
 (1) is clear. By \cite[Lem 3.4.3, (TRA1)]{Co} and \cite[p.139, (VAR1)]{Co} we have
\eq{Trcomp}{\Tr_{g\circ f}= \Tr_g\circ Rg_*(\Tr_f)\circ R(g\circ f)_*(c_{f,g}): Rg_*Rf_*(g\circ f)^!\xr{\simeq} \id.}
and  
\eq{ccomp}{c_{f,g}\circ c_{g\circ f, h}= f^!(c_{g,h})\circ c_{f, h\circ g}: (h\circ g\circ f)^!\to f^!g^!h^!,}
where $h: Z\to W$ is a third map. This implies (2).

Now to make the vertical maps in (3) precise we need some further notations:
Let $\beta_u: u^*R\sHom(-,-)\xr{\simeq} R\sHom(u^*(-), u^*(-))$, $\alpha: u^*Rf_*\xr{\simeq} Rf_*{u'}^*$ and
$e_u: u^*\xr{\simeq} u^!$ be the natural isomorphisms. Then the vertical map on the left of \eqref{etalebasechange} is given by 
$c_{u',\pi_X}^{-1}\circ e_{u'}\circ \beta_{u'}\circ \alpha$ and the vertical map on the right of \eqref{etalebasechange} is given by
$c_{u,\pi_Y}^{-1}\circ e_u\circ\beta_u$. Thus we have to prove
\eq{formulaetalebasechange}{c_{u,\pi_Y}^{-1}\circ e_u\circ\beta_u\circ u^*(f_*)= f'_*\circ c_{u',\pi_X}^{-1}\circ e_{u'}\circ \beta_{u'}\circ \alpha.}
Denote by 
\[b_{u,f}: {u'}^*f^!\xr{\simeq} {f'}^!u^* \]
the isomorphism of \cite[VII, Cor 3.4, (a), 5)]{Ha} (see also \cite[(3.3.24)]{Co}).
Then it is easy (but tedious) to check, that \eqref{formulaetalebasechange}  follows from 
\[u^*(\Tr_f)= \Tr_{f'}\circ R{f'}_*(b_{u,f})\circ\alpha:  u^*Rf_*f^!\lra R{f'}_*{f'}^!u^* \]
(see \cite[Lem 3.4.3, (TRA4)]{Co}) and the following Lemma.
\end{proof}

\begin{lemma}\label{basechangecompatibility}
Let 
 \[\xymatrix{X' \ar[r]^{u'}\ar[d]_{f'}\ar[dr]^h & X\ar[d]^f \\ Y'\ar[r]^u & Y }\]
be a cartesian diagram with $u$ \'etale.
Then the following diagram commutes (notations as above)
\[\xymatrix{{u'}^*f^!\ar[rr]^{b_{u,f}}\ar[d]_{e_{u'}} &  & {f'}^!u^*\ar[d]^{{f'}^!(e_u)}\\
             {u'}^!f^!  & h^!\ar[l]_-{c_{u',f}}\ar[r]^-{c_{f',u}} & {f'}^!u^! .}\]

\end{lemma}

\begin{proof}[Proof of Lemma.]
We will make extensive use of the notations in \cite{Ha} and \cite{Co}. 
All maps and functors involved in the statement are defined, for example, in \cite[(3.3.6),(3.3.15), (3.3.21), (3.3.25)]{Co}.
Using these definitions for the residual complex $K=\pi_Y^\Delta k$  on $Y$ together with the identity $u^*K=u^\Delta K$ and the commutativity of
\[\xymatrix{u^*D_KD_{u^*K}D_{u^*K} \ar[r]^-\eta\ar[d]_{\beta_u} & u^*D_K\ar[d]^{\beta_u}\\ 
            D_{u^*K}u^*D_KD_K\ar[r]^-\eta  & D_{u^*K}u^*,}\]
one checks that one is reduced to prove the commutativity of the following diagram:
\eq{resetalebasechange}{\xymatrix{{u'}^*f^\Delta K\ar[d]_{\varphi_{u'}^{-1}} &  & {f'}^\Delta u^* K\ar[d]^{{f'}^!(\varphi_u^{-1})}\ar[ll]_{d_{u,f}}\\
             {u'}^\Delta f^\Delta K & h^\Delta K \ar[l]_-{c_{u',f}}\ar[r]^-{c_{f',u}} & {f'}^\Delta u^\Delta K .}}
Here the maps are the analogs in the category of residual complexes of the maps in the lemma (see \cite[IV, Thm 3.1, Thm 5.5]{Ha}).
Since we work with actual complexes now, the commutativity of the above diagram is a local question. 
Thus take $U\subset X$ open such that $f_{|U}$ factors as $U\xr{i} P\xr{p} Y$ with $i$ a closed immersion and $p$ smooth.
Then $f'_{|U'}$ also factors as $U'\xr{i'} P'\xr{p'} Y'$. By the construction of $f^\Delta$ in the proof of \cite[VI, Thm 3.1, a)]{Ha} we have
$f^\Delta= i^yp^z$ and also ${f'}^\Delta= {i'}^y{p'}^z$. Now by \cite[VI, Thm 5.5, 2)]{Ha} we have 
\[c_{p,i}\circ d_{u,f}= d_{u_P,i}\circ {i'}^\Delta (d_{u,p})\circ c_{p',i'},\]
with $u_P$ being the base change of $u$ by $p:P\to Y$.
This equality implies that it is sufficient to prove the commutativity of \eqref{resetalebasechange} for $i$ and $p$ separately. Thus we are reduced to consider the two cases $f$ is finite or smooth. The case $f$ smooth is immediate, the case $f$ finite follows from \cite[Thm 3.3.1, 2. (VAR4)]{Co}. 
\end{proof}

\begin{remark}\label{restrictingtosmoothlocus}
\begin{enumerate}
 \item Let $\pi:X\to \Spec\, k$ be smooth of pure dimension $d$.
Then there is a canonical isomorphism $\pi^! k\cong \Omega^d_X[d]=:\omega_X[d]$.
More general for any $j\ge 0$ and $n\in \Z$ we have the isomorphism
\eq{isoD(Omega)}{\Omega^j_X[n]\xr{\simeq} D_X(\Omega^{d-j}_X)[n-d],}
which is defined to be the composition of 
\[\Omega^j_X[n]\xr{\simeq}\sHom_X(\Omega^{d-j}_X, \Omega^d_X)[n], \quad \alpha\mapsto (\beta\mapsto \alpha\wedge\beta)\]
(notice that we make the choice of a sign here) with
\[\sHom_X(\Omega^{d-j}_X, \Omega^d_X)[n]=\sHom^\bullet(\Omega^{d-j}_X,\Omega^d_X[d])[n-d]\cong R\sHom(\Omega^{d-j}_X, \pi^! k)[n-d].\]
\item Let $X$ be a $k$-scheme of pure dimension $d$ and $U\subset X$ a smooth open subscheme, then we have for all $j\ge 0$
\[\Omega^j_U\simeq \sHom^\bullet_U(\Omega^{d-j}_U, \Omega^d_U[d])[-d]\simeq D_X(\Omega^{d-j}_X)_{|U}[-d],\]
where the first isomorphism is as in (1) and the second is given by restriction
(or to be more precise, first use the isomorphism $\Omega^d_U[d]\simeq \pi_U^! k$ and then the vertical isomorphism on the right in \eqref{etalebasechange} with $U\inj X$ instead of $u: Y'\to Y$). 
\end{enumerate}

 \end{remark}
 
\begin{lemma}\label{pf-for-cartesianproduct}
Let $\pi_X: X\to \Spec\, k$ be proper of pure dimension $d_X$ and $\pi_Y: Y\to \Spec \, k$ smooth of pure dimension $d_Y$. 
We denote by $\pr_2:X\times Y \to Y$ the projection (it is proper) and set $d:=\dim (X\times Y)$. Then for all $j\ge 0$ there is a morphism in $D^+_{\rm c}(X\times Y)$
\[\gamma: \pr_2^!(\OO_Y)\otimes (\pr_2^*\Omega^{j-d_X}_Y[d_Y])\lra D_{X\times Y}(\pr_2^*\Omega^{d-j}_Y)\]
satisfying the following conditions:
\begin{enumerate}
\item For $U\subset X$ open and smooth denote by $p_2:U\times X \to Y$ the restriction of $\pr_2$. Then $\gamma\mid U\times Y$ is the composition
  \mlnl{(\pr_2^!(\OO_Y)\otimes \pr_2^*\Omega^{j-d_X}_Y[d_Y])_{|U\times Y}\xr{\simeq} \Omega^{d_X}_{U\times Y/Y}[d_X]\otimes p_2^*\Omega^{j-d_X}_Y[d_Y]\\
         \xr{\simeq}  \Omega^{d_X}_{U\times Y/Y}[d_X]\otimes p_2^*R\sHom_Y(\Omega^{d-j}_Y, \Omega_Y^{d_Y}[d_Y])
                                                                 \xr{\simeq\,\text{nat.}} D_{U\times Y}(p_2^*\Omega^{d-j}_Y).  }
Here the last isomorphism is induced by the composition of the canonical isomorphisms 
$\Omega_{U\times Y/Y}^{d_X}[d_X]\otimes p_2^*\Omega^{d_{Y}}_Y[d_Y]\cong\Omega^d_{U\times Y}[d]\cong\pi_{U\times Y}^!(k)$.
\item The following diagram commutes:
    \[\xymatrix@C=6pt{ R\pr_{2*}(\pr_2^!(\OO_Y)\otimes \pr_2^*\Omega^{j-d_X}_Y[d_Y])\ar[rr]^\gamma\ar[d]_{\text{proj. formula}} & &
                R\pr_{2*}D_{X\times Y}(\pr_2^*\Omega^{d-j}_Y)\ar[d]\\
                R\pr_{2*}(\pr_2^!(\OO_Y))\otimes\Omega^{j-d_X}_Y[d_Y]\ar[dr]_{\Tr_{\pr_2}\otimes \id}   &  & 
                 R\sHom_Y(\Omega^{d-j}_Y, \pi_Y^!k) \\
                 & \Omega^{j-d_X}_Y[d_Y],\ar[ur]_{\eqref{isoD(Omega)}}^\simeq   }\]
                 where the vertical map on the right is $\Tr_{\pr_2}\circ{\rm adjunction}\circ c_{\pr_2, \pi_Y}$. 
\end{enumerate}
\end{lemma}

\begin{proof}
In \cite[(4.3.12)]{Co} is defined a map 
\[e_{\pr_2}: \pr_2^!(\OO_Y)\otimes^L\pr_2^*\pi_Y^! k \lra \pr_2^!\pi_Y^! k\] 
such that 
\eq{e-restr}{\xymatrix{ (\pr_2^!(\OO_Y)\otimes^L\pr_2^*\pi_Y^! k)_{|U\times Y}\ar[r]^-{e_{\pr_2}\mid U\times Y} & (\pr_2^!\pi_Y^! k)_{|U\times Y}\\
                  \Omega^{d_X}_{U\times Y/Y}[d_X]\otimes p_2^*\Omega^{d_Y}_Y[d_Y]\ar[u]_\simeq\ar[r]^-\simeq & \Omega^d_{U\times Y}[d]\ar[u]^\simeq}}
commutes, where the vertical map on the left is the composition of the canonical isomorphism $\Omega^d_{U\times Y}[d]\cong \pi_{U\times Y}^!(k)$ with $c_{p_2, \pi_Y}:\pi_{U\times Y}^!\cong p_2^!\pi_Y^!$.
Furthermore by \cite[Thm 4.4.1]{Co} the following diagram commutes:

\eq{e-trace}{\xymatrix{ R\pr_{2*}(\pr_2^!\OO_Y\otimes^L\pr_2^*\pi_Y^!k)\ar[r]^-{e_{\pr_2}}\ar[d]_{\text{ proj. formula}} &                                                                                     R\pr_{2*}\pr_2^!(\pi_Y^!k)\ar[d]^{\Tr_{\pr_2}}\\    
              R\pr_{2*}(\pr_2^!\OO_Y)\otimes^L\pi_Y^!k\ar[r]_-{\Tr_{\pr_2}\otimes\id} & \pi_Y^!k.  }}
We define $\gamma$ to be the composition
\begin{eqnarray*}
\pr_2^!(\OO_Y)\otimes\pr_2^*\Omega^{j-d_X}_Y[d_Y] & \xr{\id\otimes\eqref{isoD(Omega)}}      &  \pr_2^!(\OO_Y)\otimes\pr_2^*R\sHom_Y(\Omega^{d-j}_Y,\pi_Y^!k)\\
                                                  & \xr{nat.}                               &  R\sHom(\pr_2^*\Omega^{d-j}_Y, \pr_2^!(\OO_Y)\otimes^L\pr_2^*\pi_Y^!k)\\
                                                  &\xr{c_{\pr_2,\pi_Y}^{-1}\circ e_{\pr_2}} & D_{X\times Y}(\pr_2^*\Omega^{d-j}_Y).
\end{eqnarray*}
It follows from  \eqref{e-restr} and \eqref{e-trace}, that $\gamma$ satisfies (1) and (2).
\end{proof}

\begin{proposition}\label{1.7.9}
Let $i: X\inj Y$ be a closed immersion of pure codimension $c$ between smooth $k$-schemes of pure dimension $d_X$ and $d_Y$, respectively.
Then for all $q\ge 0$
\[R\ul{\Gamma}_X \Omega_Y^q[c]\cong \sH^c_X(\Omega^q_Y) \quad \text{in } D^b_{\rm qc}(\sO_Y).\]
Suppose further the ideal sheaf of $X$ in $\sO_Y$ is generated by a sequence $t=t_1,\ldots, t_c$  of global sections of $\sO_Y$.
Define a morphism $\imath_X^q$ by
\[\imath_X^q: i_*\Omega^q_X\to \sH^c_X(\Omega^{c+q}_Y),\quad \alpha\mapsto (-1)^c\genfrac{[}{]}{0pt}{}{dt\tilde{\alpha}}{t},\]
where $\tilde{\alpha}\in \Omega^q_Y$ is any lift of $\alpha$ and $dt=dt_1\wedge\ldots\wedge dt_c$. (Here we use the notation of \ref{A1}.) Then the following diagram commutes in $D^b_{\rm qc}(\sO_Y)$
\eq{1.7.9.1}{\xymatrix{i_*\Omega^q_X\ar[r]^-{\eqref{isoD(Omega)}}\ar[dr]_{\imath_X^q} & i_*D_X(\Omega_X^{d_X-q})[-d_X]\ar[r]^{i_*} & D_Y(\Omega^{d_X-q}_Y)[-d_X]\ar[r]^-{\eqref{isoD(Omega)}} & \Omega^{c+q}_Y[c]\\
                                                                                 &  \sH^c_X(\Omega_Y^{c+q})\ar[r]^\simeq          & R\ul{\Gamma}_X(\Omega_Y^{c+q})[c].\ar[ur]       } } 

\end{proposition}
\begin{proof}
The first statement is well-known (see also Lemma \ref{trace-for-lci-appendix}). It remains to prove the commutativity of \eqref{1.7.9.1}.  Let 
$\pi_X: X\to \Spec k$ and $\pi_Y: Y\to \Spec k$ be the structure maps.  By Definition \ref{defpush-forward} the top row in \eqref{1.7.9.1} is given by the 
following composition in $D^b_{\rm qc}(\sO_Y)$ 
\begin{eqnarray*}
i_{*}\Omega^q_X    & \xr{\eqref{isoD(Omega)}}      &  i_* R\sHom(\Omega_X^{d_X-q}, \pi_X^!k)[-d_X] \\
                   & \xr{c_{i,\pi_Y}}              &  i_*R\sHom(\Omega_X^{d_X-q}, i^!\pi_Y^!k)[-d_X] \\
                   & \xr{\text{nat.}}              &  R\sHom(i_*\Omega_X^{d_X-q}, i_*i^!\pi_Y^!k)[-d_X] \\
                   & \xr{\Tr_i}                    &  R\sHom(i_*\Omega_X^{d_X-q}, \pi_Y^!k)[-d_X] \\
                   & \xr{ (i^*)^\vee}              &  R\sHom(\Omega_Y^{d_X-q}, \pi_Y^!k )[-d_X] \\
                   & \xr{\eqref{isoD(Omega)}^{-1}} &  \Omega_Y^{c+q}[c]. 
\end{eqnarray*}
We set $\imath_X:=\imath^{d_X}_X$. Then it follows from Lemma \ref{1.7.11} and the definition of \eqref{isoD(Omega)}, that the above composition equals 
\begin{eqnarray}\label{1.7.9.2}
i_*\Omega_X^q      & \xr{\text{multipl.}}          &  i_*\sHom(\Omega_X^{d_X-q}, \Omega_X^{d_X})\\
                   & \xr{\text{nat.}}              &  \sHom(i_*\Omega_X^{d_X-q}, i_*\Omega_X^{d_X}) \nonumber\\
                   & \xr{\imath_X}                 &  \sHom(i_*\Omega_X^{d_X-q},\sH^c_X(\Omega_Y^{d_Y})) \nonumber\\
                   & \xr{(i^*)^\vee}                   &  \sHom(\Omega_Y^{d_X-q}, \sH^c_X(\Omega_Y^{d_Y})) \nonumber\\
                   & \xr{(*)}                      &  \sHom^\bullet(\Omega_Y^{d_X-q}, \Omega_Y^{d_Y}[c]) \nonumber\\
                   & \xr{\text{multipl.}^{-1}}     &  \Omega^{c+q}_Y[c],     \nonumber
\end{eqnarray}
where $(*)$ is induced by $\sH^c_X(\Omega_Y^{d_Y})\cong R\ul{\Gamma}_X(\Omega_Y^{d_Y})[c]\to \Omega_Y^{d_Y}[c]$.
There is a natural isomorphism $\varphi: \sH^c_X(\Omega^{c+q}_Y)\xr{\simeq} \sHom(\Omega^{d_X-q}_Y,\sH^c_X(\Omega^{d_Y}_Y) )$ coming from the isomorphisms
\begin{align*}
\sHom(\Omega^{d_X-q}_Y,\sH^c_X(\Omega^{d_Y}_Y) ) & \cong R\sHom(\Omega^{d_X-q}_Y,R\ul{\Gamma}_X(\Omega^{d_Y}_Y)[c])\\ 
                                                 & \cong R\ul{\Gamma}_X(R\sHom(\Omega^{d_X-q}_Y,\Omega^{d_Y}_Y))[c] \cong R\ul{\Gamma}_X(\Omega^{c+q}_Y)[c]\\
                                                 & \cong \sH^c_X(\Omega_Y^{c+q}).
 \end{align*}
This isomorphism is explicitly given by
\[\varphi: \sH^c_X(\Omega^{c+q}_Y)\xr{\simeq} \sHom(\Omega^{d_X-q}_Y,\sH^c_X(\Omega^{d_Y}_Y) ),\quad \genfrac{[}{]}{0pt}{}{\alpha}{t^n}\mapsto \left(\beta\mapsto \genfrac{[}{]}{0pt}{}{\alpha\beta}{t^n}\right).\]
The composition \eqref{1.7.9.2} equals
\[i_*\Omega_X^q \xr{\varphi^{-1}\circ (i^*)^\vee\circ \imath_X\circ (\text{nat.})\circ (\text{multipl.})} \sH^c_X(\Omega^{c+q}_Y)\cong R\ul{\Gamma}_X(\Omega^{c+q}_Y)[c]\to \Omega^{c+q}_Y[c].\]
It is straightforward to check that $\imath^q_X= \varphi^{-1}\circ (i^*)^\vee\circ \imath_X\circ (\text{nat.})\circ (\text{multipl.})$ and this implies the commutativity of \eqref{1.7.9.1}.
\end{proof}

\begin{corollary}\label{pf-pb-cl-emb}
Assume we have a cartesian square
\[\xymatrix{X'\ar@{^(->}[r]^{i'}\ar[d]_{g_X} & Y'\ar[d]^{g_Y}\\
            X\ar@{^(->}[r]^{i}&  Y,   }\]
in which $X, X', Y, Y' $ are smooth of pure dimension $d_X, d_{X'}, d_Y, d_{Y'}$, $i$ is a closed immersion and we have
$c:= d_Y-d_X=d_{Y'}-d_{X'}$. Then for all $q\ge 0$ the following diagram commutes in $D^b_{\rm qc}(Y)$
\[\xymatrix{ i_*Rg_{X*}\Omega^q_{X'}=Rg_{Y*}i'_*\Omega^q_{X'}\ar[r]     & Rg_{Y*}\Omega^{c+q}_{Y'}[c]\\
             i_*\Omega^q_X\ar[r]\ar[u]^{g_X^*}                          &  \Omega^{c+q}_Y[c],\ar[u]^{g_Y^*}      }\]
where the lower horizontal morphism is given by the composition
\[i_*\Omega^q_X \xr{\eqref{isoD(Omega)}} i_*D_X(\Omega_X^{d_X-q})[-d_X]\xr{i_*} D_Y(\Omega^{d_X-q}_Y)[-d_X]\xr{\eqref{isoD(Omega)}} \Omega^{c+q}_Y[c]\]
and the upper horizontal morphism by $Rg_{Y*}$ applied to the analogous map for $i'$.
\end{corollary}

\begin{proof}
Since $R\ul{\Gamma}_X Rg_{Y*}= Rg_{Y*}R\ul{\Gamma}_{X'}$ we naturally have a commutative diagram
\[\xymatrix{Rg_{Y*}\sH^c_{X'}(\Omega^{c+q}_{Y'})\ar[r]^\simeq & Rg_{Y*}R\ul{\Gamma}_{X'}(\Omega^{c+q}_{Y'})[c]\ar[r] &  Rg_{Y*}(\Omega^{c+q}_{Y'})[c]\\
            \sH^c_{X}(\Omega^{c+q}_Y)\ar[u]^{g_Y^*}\ar[r]^\simeq & R\ul{\Gamma}_X(\Omega^{c+q}_Y)[c]\ar[u]^{g_Y^*}\ar[r] &  \Omega^{c+q}_Y[c]\ar[u]^{g_Y^*},     }\]
where the $g_Y^*$ on the very left is defined in such a way that the left square commutes. By Proposition \ref{1.7.9} it thus suffices to prove the commutativity of 
\[\xymatrix{  i_*g_{X*}\Omega^q_{X'}= g_{Y*}i'_*\Omega^q_{X'}\ar[r]^-{\imath^q_{X'}} & g_{Y*}\sH^c_{X'}(\Omega^{c+q}_{Y'}) \\
                i_*\Omega^q_X\ar[u]^{g_X^*}\ar[r]^{\imath^q_X} & \sH^c_X(\Omega^{c+q}_Y).\ar[u]^{g_Y^*}  }\]
This is a local question, we may therefore assume that the ideal of $X$ in $Y$ is generated by a sequence $t_1,\ldots, t_c$ of global sections of $\sO_Y$.
Then $g_Y^*t_1,\ldots , g_Y^*t_c$ is a sequence of global sections of $\sO_{Y'}$, which generate the ideal sheaf of $X'$ in $Y'$. Hence the assumption
follows from the explicit description of $\imath^q_X$ and $\imath^q_{X'}$ in Proposition \ref{1.7.9}.
\end{proof}

\begin{proposition}\label{finite-pf}
Let $f:X\to Y$ be a finite and surjective morphism between smooth schemes, which are both of pure dimension $n$. We denote by
\[\tau_f: \oplus_q f_*\Omega^q_X\to \oplus_q\Omega^q_Y\]
the composition
\[\oplus_q f_*\Omega^q_X\xr{\eqref{isoD(Omega)}} \oplus_q f_*D_X(\Omega^{n-q}_X)\xr{f_*} \oplus_q D_Y(\Omega^{n-q}_Y)\xr{\eqref{isoD(Omega)}} \oplus_q \Omega^q_Y.\]
Then:
\begin{enumerate}
 \item In degree 0 the map $\tau_f$ equals the usual trace on the finite and locally free $\sO_Y$-module $f_*\sO_X$,   $Tr_{X/Y}: f_*\sO_X\to \sO_Y$.
 \item For $\alpha\in f_*\Omega^a_X$ and $\beta\in \Omega^b_Y$ we have
           \[\tau_f(\alpha f^*\beta)=\tau_f(\alpha)\beta\]
\item The composition $\tau_f\circ f^*: \oplus_q\Omega^q_Y\to \oplus_q\Omega^q_Y$ equals multiplication with the degree of $f$.
\end{enumerate}

\end{proposition}
\begin{proof}
All statements are local in $Y$. We may therefore assume, that $f$ factors as $X\xr{i} P\xr{\pi} Y$, where $i$ is a regular closed immersion of pure codimension $d$ and $\pi$ is smooth of relative dimension $d$; further
we may assume that the ideal sheaf of $X$ in $P$ is generated by $d$ global sections $t_1,\ldots, t_d$ of $\sO_P$. Then in degree $n$ the map $\tau_f$
equals the trace map $\tau_f^n: f_*\omega_X\to \omega_Y$ from \ref{A3.1} and in degree $q$ the map $\tau_f$ thus equals the composition
\mlnl{f_*\Omega^q_X\cong f_*\sHom_X(\Omega^{n-q}_X,\omega_X)\xr{{\rm nat.}}\sHom_Y(f_*\Omega^{n-q}_X, f_*\omega_X)\\
            \xr{\tau_f^n} \sHom(f_*\Omega^{n-q}_X,\omega_Y ) \xr{\circ f^*}\sHom(\Omega^{n-q}_Y, \omega_Y)\cong \Omega^q_Y. }
Thus  for $\alpha\in f_*\Omega_X^q$ Lemma \ref{A3.2} gives the following formula for $\tau_f(\alpha)$: 
In 
           \[i^*\Omega_P^{d+q}=\bigoplus_{r+s=d+q} i^*(\Omega^r_{P/Y})\otimes f^*\Omega^s_Y\]
         write
\[i^*(dt_d\wedge\ldots\wedge dt_1\wedge \tilde{\alpha})=\sum_{r+s=d+q} \sum_j i^*\gamma_{j,r}\otimes f^*\beta_{j,s} ,\quad \gamma_{j,r}\in \Omega_{P/Y}^r, \, \beta_{j,s}\in\Omega^s_Y,\]
where $\tilde{\alpha}\in \Omega^q_P$ is a lift of $\alpha$.
Then  
\eq{formula-finite-pf}{\tau_f(\alpha) = (-1)^{d(d-1)/2} \sum_j  \Res_{P/Y}\genfrac{[}{]}{0pt}{}{\gamma_{j,d}}{t_1,\ldots, t_d} \beta_{j,q} \in \Omega_Y^q. }
This formula immediately implies (2). For any  $a\in \sO_X$ we have 
\[\tau_f(a)=(-1)^{d(d-1)/2} \Res_{P/Y}\genfrac{[}{]}{0pt}{}{\tilde{a}dt_d\wedge\ldots\wedge dt_1}{t_1,\ldots, t_d}=\Res_{P/Y}\genfrac{[}{]}{0pt}{}{\tilde{a}dt_1\wedge\ldots\wedge dt_d}{t_1,\ldots, t_d},\]
which equals $\Tr_{X/Y}(a)$ by \cite[p.240, (R6)]{Co}, hence (1). Finally (3) is a direct consequence of (1) and (2).
\end{proof}

\begin{remark}
The trace map from Proposition \ref{finite-pf} and its properties are well-known, see e.g. \cite[\S 16]{Ku}, where the trace is considered in a much greater generality.
But there the construction is done via an ad hoc method not using the duality formalism. Therefore the connection to the trace map above is not a priori clear.
\end{remark}

\subsection{Push-forward for Hodge cohomology with support}
\begin{definition}\label{compactification}
Let $f: (X, \Phi)\to (Y,\Psi)$ be a morphism in $V_*$ with $X$ equidimensional. We define a compactification of $f$ to be a factorization 
\[f=\bar{f}\circ j: (X,\Phi)\inj (\bar{X}, \Phi)\lra (Y,\Psi),\]
where  $\bar{X}$ is equidimensional (but possibly singular), $j$ is an open immersion and $\bar{f}$ is proper.
Notice that since $f\mid \Phi$ is proper, $\Phi$ is also a family of supports on $\bar{X}$. 
The compactification will be denoted by $(j, \bar{f})$.

By Nagata's compactification theorem (see, e.g. \cite{Co2}) any $f$ in $V_*$ admits  a compactification.
\end{definition}

\begin{definition}[Push-forward]\label{pfwithsupport}
Let $f: (X, \Phi)\to (Y,\Psi)$ be a morphism in $V_*$ and assume that $X$ and $Y$ are of pure dimension $d_X$ and $d_Y$ respectively
and set $r:=d_X-d_Y$. 
Let 
\[(X, \Phi)\xr{j} (\bar{X}, \Phi)\xr{\bar{f}} (Y,\Psi),\]
be a compactification of $f$.  We define the push-forward
\[H_*(f): H(X,\Phi)\lra H(Y,\Psi)\]
as the following composition:
\mlnl{H(X,\Phi)\simeq \bigoplus_{i,j} H^i_{\Phi}(\bar{X}, D_{\bar{X}}(\Omega^{d_X-j}_{\bar{X}})[-d_X])\xr{\text{nat.}} 
      \bigoplus_{i,j} H^{i-d_X}_{\bar{f}^{-1}(\Psi)}(\bar{X}, D_{\bar{X}}(\Omega^{d_X-j}_{\bar{X}}))\\
 \xr{\oplus \bar{f}^j_*} \bigoplus_{i,j} H^{i-d_X}_{\Psi}(Y, D_{Y}(\Omega^{d_X-j}_{Y}))
  \xr{\simeq\, \eqref{isoD(Omega)}} \bigoplus_{i,j} H^{i-r}_{\Psi}(Y, \Omega^{j-r}_{Y})= H(Y,\Psi),}
where the first isomorphism is the composition of \eqref{isoD(Omega)} for $n=0$ with the excision isomorphism. 
Notice that we obtain a morphism of graded abelian groups $H_*(f): H_*(X,\Phi)\to H_*(Y,\Psi)$, see \eqref{Hodgelowergrading}.

This definition is independent of the chosen compactification. 

We extend the definition to the case of non-equidimensional $X$ and $Y$ additively. 
\end{definition}
\begin{proof}
We have to prove the independence of $H_*(f)$ from the chosen compactification. 
Let 
\[\xymatrix{ X\ar[r]^{j_2}\ar[d]_{j_1}   & X_2 \ar[d]^{f_2}   \\
               X_1\ar[r]_{f_1}\ar[ur]^g & Y          }\]
be a commutative diagram with $d:=\dim X_1=\dim X_2=d_X$, $j_1$ and $j_2$ open and $f_1$, $f_2$ proper. Notice that $g$ is automatically proper.
Then the following diagram commutes:
\[\xymatrix@=11pt{  &  H^{i-d}_{\Phi}(D_{X_2}(\Omega^{d-j}_{X_2}))\ar[r] &
                                          H^{i-d}_{f_2^{-1}(\Psi)}(D_{X_2}(\Omega^{d-j}_{X_2}))\ar[dr]^-{f_{2*}} & \\
             H^i_\Phi(\Omega^j_X)\ar[ur]^-\simeq\ar[dr]_-\simeq  & &  & H^{i-d}_\Psi(D_Y(\Omega^{d-j}_Y)). \\ 
             &  H^{i-d}_{\Phi}(D_{X_1}(\Omega^{d-j}_{X_1}))\ar[r]\ar[uu]_{g_*} & 
                H^{i-d}_{f_1^{-1}(\Psi)}(D_{X_1}(\Omega^{d-j}_{X_1}))\ar[ur]_-{f_{1*}}\ar[uu]^{g_*}  &
            }\]
Indeed, the left triangle commutes since $g_*\mid X= \id_*$ by \ref{properties-pf}, (3), the square in the middle obviously commutes and
the triangle on the right commutes by \ref{properties-pf}, (2). 

Two arbitrary compactifications of $f$ always receive a map from a third one and thus the general case follows from the case above.
 
\end{proof}

\begin{proposition}\label{properties-pfwithsupports}
 \begin{enumerate}
   \item $H_*(\id)=\id$.
  \item Let $f:(X,\Phi)\to (Y,\Psi)$ and $g:(Y,\Psi)\to (Z,\Xi)$ be two morphisms in $V_*$. Then 
           \[H_*(g\circ f) = H_*(g)\circ H_*(f) : H(X,\Phi)\lra H(Z,\Xi).\]
  \item If $f:(X,\Phi)\to (Y,\Psi)$ in $V_*$ is finite, then $H_*(f)$ is induced by the trace map $\tau_f$ from Proposition \ref{finite-pf}.
 \end{enumerate}

\end{proposition}
\begin{proof}
(1) follows from \ref{properties-pf}, (1).  Now for (2) we may assume that $X,Y,Z$ are connected. 
Let $(j_X, f_1)$ and $(j_Y, g_1)$ be  compactifications of $f$ and $g$, respectively.
Let $(j_{X_1}, f_2)$ be a compactification of $j_{Y}\circ f_1$. Thus we have a commutative diagram
\[\xymatrix@=15pt{X_2 \ar[dr]^-{f_2}      & &    \\ 
            X_1\ar[u]^-{j_{X_1}}\ar[dr]^(.4){f_1} & Y_1\ar[dr]^-{g_1}\\
            X\ar[u]^-{j_X}\ar[r]_-f & Y\ar[u]_(.4){j_Y}\ar[r]_-g & Z,  
}\]
with vertical arrows open immersions and diagonal arrows proper. Replacing $X_1$ by $f_2^{-1}(Y)$, we may assume that the parallelogram is cartesian. 
Now by \ref{properties-pf}, (3) the following diagram commutes
\[\xymatrix{ H^i_{f_2^{-1}(\Psi)}(D_{X_2}(\Omega^j_{X_2}))\ar[r]^-{f_{2 *}}\ar[d]^\simeq & H^i_{\Psi}(D_{Y_1}(\Omega^j_{Y_1}))\ar[d]^\simeq\\
             H^i_{f_1^{-1}(\Psi)}(D_{X_1}(\Omega^j_{X_1}))\ar[r]^-{f_{1 *}} & H^i_{\Psi}(D_Y(\Omega^j_Y)) .}\]
Thus (2) follows from \ref{properties-pf}, (2). (3) follows immediately from the definitions.
\end{proof}

\begin{lemma}\label{Hodgepf-cartesianproduct}
Consider a cartesian diagram
\[\xymatrix{(X\times Y', \Phi')\ar[r]^-{f'}\ar[d]_{g_{X\times Y}} & (Y', \Psi')\ar[d]^{g_Y}\\
             (X\times Y, \Phi)\ar[r]^-f    & (Y,\Psi),   }\]
            such that $f$ is induced by the projection to $Y$, $f,f'\in V_*$ and $g_{X\times Y}, g_Y\in V^*$. Then
        \[H^*(g_Y)\circ H_*(f)=H_*(f')\circ H^*(g_{X\times Y}).\]
 Furthermore, $H_*(f): H(X\times Y,\Phi)\to H(Y,\Psi)$ factors over the projection  
 $H(X\times Y, \Phi)\to \oplus_{i,j}H^i_\Phi(X\times Y, \pr_1^*\Omega^{d_X}_X\otimes \pr_2^*\Omega^j_Y )$.
\end{lemma}

\begin{proof}
We may assume $X$ and $Y$ of pure dimension $d_X$ and $d_Y$, respectively and we set $d:=d_X+d_Y$. 
We embed $X$ as an open in a proper $k$-scheme $\bar{X}$ of pure dimension $d_X$. Then 
$(X\times Y,\Phi)\xr{j} (\bar{X}\times Y, \Phi)\xr{\pr_2} (Y,\Psi)$ is a compactification of $f$, where $j$ is the open embedding and $\pr_2$ is induced by the projection to $Y$. Similar we obtain a compactification for $f'$, in which case we write $\pr_2'$ for the projection to $Y'$.
The second statement of the lemma follows from Definition \ref{pfwithsupport}, Remark \ref{alt-def-pf}, (2)  and the following commutative diagram:
\[\xymatrix{\Omega^j_{X\times Y}\ar[r]^\simeq\ar[d]_{\text{projection}} & 
        D_{X\times Y}(\Omega^{d-j}_{X\times Y} )[-d]\ar[d]^{(\pr_2^* )^\vee } \\
         \pr_1^*\Omega^{d_X}_X\otimes \pr_2^*\Omega^{j-d_X}_Y\ar[r]^-\simeq & D_{X\times Y}(\pr_2^*\Omega^{d-j}_Y)[-d].   }\]

Now we come to the first statement of the lemma. Consider the following diagram (we use a shortened notation):
 \[\xymatrix{  &   H^{i-d}_{\pr_2^{-1}(\Psi)}(D_{\bar{X}\times Y}(\Omega^{j-d}_{\bar{X}\times Y}))\ar[d]^{(\pr_2^*)^\vee}\ar[dr]^{{\pr_2}_*} &\\
           H^i_\Phi(\Omega^j_{X\times Y})\ar[ur]\ar[r]\ar[dr] & H^{i-d}_{\pr_2^{-1}(\Psi)}(D_{\bar{X}\times Y}(\pr_2^*\Omega^{j-d}_Y))\ar[r]  &
                   H^{i-d_X}_\Psi(\Omega^{j-d_X}_Y), \\
      &  H^{i-d}_{\pr_2^{-1}(\Psi)}(\pr_2^!(\OO_Y)\otimes \pr_2^*\Omega^{j-d_X}_Y[d_Y])\ar[u]^{\gamma}\ar[ur] & 
    }\] 
 here we use the notation  of Lemma  \ref{pf-for-cartesianproduct}, furthermore the upper map on the left is induced by excision, the middle and the lower map on the left are induced
 by projection and excision and the middle and the lower map on the right are induced
 by the corresponding maps from Lemma \ref{pf-for-cartesianproduct}, (2). It follows from Lemma \ref{pf-for-cartesianproduct} and 
 Remark \ref{alt-def-pf}, (2) that all the triangles in this diagram commute. Replacing $Y$ by $Y'$ and $\pr_2$ by $\pr_2'$  we obtain a similar 
 commutative diagram. Thus it remains to show that the following diagram is commutative:
 \eq{Hodgepf-cartesianproduct-diagram}{
   \xymatrix{ H^i_\Phi(\Omega^j_{X\times Y})\ar[d]^{H^*(g_{X\times Y})}\ar[r]^-{\text{proj.}} & 
                               H^{i-d}_{\pr_2^{-1}(\Psi)}(\pr_2^!(\OO_Y)\otimes \pr_2^*\Omega^{j-d_X}_Y[d_Y])\ar[r]^-{\Tr_{\pr_2}\otimes \id}&   H^{i-d_X}_\Psi(\Omega^{j-d_X}_Y)\ar[d]^{H^*(g_Y)}\\
              H^i_{\Phi'}(\Omega^j_{X\times Y'})\ar[r]^-{\text{proj.}} & 
                     H^{i-d'}_{{\pr_2'}^{-1}(\Psi')}({\pr_2'}^!(\OO_{Y'})\otimes {\pr_2'}^*\Omega^{j-d_X}_{Y'}[d_{Y'}])\ar[r]^-{\Tr_{\pr_2'}\otimes\id}&  H^{i-d_X}_{\Psi'}(\Omega^{j-d_X}_{Y'}),
            }                           }
where $d'=d_X+d_{Y'}$. To this end we define the map 
\[\tau_f: Rf_*R\ul{\Gamma}_\Phi(\omega_{X\times Y/Y}[d_X])\to R\ul{\Gamma}_{\Psi}\sO_Y,\]
to be the composition
\mlnl{Rf_*R\ul{\Gamma}_\Phi(\omega_{X\times Y/Y}[d_X])\xr{\text{excision }\simeq} R\pr_{2*}R\ul{\Gamma}_\Phi(\pr_2^!\sO_Y)\\ 
                                 \xr{\rm nat.}R\pr_{2*}R\ul{\Gamma}_{\pr_2^{-1}(\Psi)}(\pr_2^!\sO_Y) \xr{\simeq} R\ul{\Gamma}_\Psi R\pr_{2*}\pr_2^!\sO_Y\xr{\Tr_{\pr_2}} R\ul{\Gamma}_\Psi\sO_Y.}
Then the upper horizontal line in diagram \eqref{Hodgepf-cartesianproduct-diagram} equals $H^{i-d}(Y,-)$ applied to the following composition:
\mlnl{Rf_*R\ul{\Gamma}_\Phi\Omega^j_{X\times Y}[d]\xr{\rm projection} Rf_*R\ul{\Gamma}_\Phi(\omega_{X\times Y/Y}[d_X]\otimes f^*\Omega^{j-d_X}_Y[d_Y] )\\
\xr{\simeq} Rf_*R\ul{\Gamma}_\Phi(\omega_{X\times Y/Y}[d_X])\otimes \Omega^{j-d_X}_Y[d_Y] \xr{\tau_f\otimes\id} R\ul{\Gamma}_\Psi(\Omega^{j-{d_X}}_Y)[d_Y].}
(There is no intervention of signs in the definition of the projection map, this is compatible with the fact, that the isomorphism 
$\omega_{X\times Y}[d]\cong \omega_{X\times Y/Y}[d_X]\otimes f^*\omega_{Y}[d_Y]$ is defined without a sign, see \cite[(2.2.6)]{Co}. )
The lower horizontal line in the diagram \eqref{Hodgepf-cartesianproduct-diagram} equals $H^{i-d'}(Y',-)$ applied to the analog composition for $f'$. Then 
it is straightforward to check that the commutativity of diagram \eqref{Hodgepf-cartesianproduct-diagram} is implied by the commutativity of
\eq{diagram2}{\xymatrix{Rf_*R\ul{\Gamma}_\Phi(\omega_{X\times Y/Y}[d_X])\ar[r]^-{\tau_f}\ar[d]_{g_{X\times Y}^*} & R\ul{\Gamma}_{\Psi} \sO_Y \ar[d]^{g_Y^*}\\
            Rg_{Y*}Rf'_*R\ul{\Gamma}_{\Phi'}(\omega_{X\times Y'/Y'}[d_X])\ar[r]^-{\tau_{f'}} & Rg_{Y*}R\ul{\Gamma}_{\Psi'} \sO_{Y'}.    } }

To prove the commutativity of this last diagram, we can clearly assume (by definition of the pull-back and $\tau_f$), that $\Phi'=g_{X\times Y}^{-1}(\Phi)$ and $\Psi'=g_{Y}^{-1}(\Psi)$.
We define the map 
\[\alpha:R\pr_{2*}\pr_2^!\sO_Y\to Rg_{Y*}R\pr'_{2*}(\pr_2')^!\sO_{Y'} \]  to be following composition
\begin{eqnarray}
R\pr_{2*}\pr_2^!(\pi_Y^*k)  & \xr{b_{\pi_Y, \pi_{\bar{X}}}^{-1}}& R\pr_{2*}\pr_1^*(\pi_{\bar{X}}^!k)\nonumber\\
                            &\xr{}&  Rg_{Y*}R\pr'_{2*}(\pr_1')^*(\pi_{\bar{X}}^!k)\nonumber\\
                            & \xr{b_{\pi_{Y'},\pi_{\bar{X}}}}&  Rg_{Y*} R\pr'_{2*}(\pr_2')^!(\pi_{Y'}^*k)=Rg_{Y*}R\pr'_{2*}(\pr_2')^!\sO_{Y'}\nonumber,
\end{eqnarray}
where $b_{\pi_Y, \pi_{\bar{X}}}: \pr_1^*\pi_{\bar{X}}^!\simeq \pr_2^!\pi_Y^*$ is the isomorphism from \cite[VII, Cor 3.4, 5)]{Ha} and the middle map is the composition of the natural maps
\[R\pr_{2*}\pr_1^*\to Rg_{Y*}Lg_{Y}^*R\pr_{2*}\pr_1^*\to Rg_{Y*}R\pr'_{2*}Lg_{\bar{X}\times Y}^*\pr_1^*\cong Rg_{Y*} R\pr'_{2*}(\pr_1')^*.\]
Now the commutativity of diagram \eqref{diagram2} follows from the commutativity of 
\[\xymatrix{Rf_*R\ul{\Gamma}_\Phi(\omega_{X\times Y/Y}[d_X])\ar[r]\ar[d]_{g_{X\times Y}^*} &  R\ul{\Gamma}_{\Psi}R\pr_{2*}\pr_2^!\sO_Y\ar[d]_{\alpha} \\
            Rg_{Y*}Rf'_*R\ul{\Gamma}_{\Phi'}(\omega_{X\times Y'/Y'}[d_X])\ar[r]  &  R\ul{\Gamma}_{\Psi}Rg_{Y*}R\pr'_{2*}(\pr_2')^!\sO_{Y'}, }\]
which is clear by the explicit description of the isomorphisms $b_{(-,-)}$ in the smooth case (see \cite[VII, Cor. 3.4, (a), Var 6]{Ha});
and from the commutativity of the diagram
\[\xymatrix{R\ul{\Gamma}_{\Psi}R\pr_{2*}\pr_2^!\sO_Y\ar[d]_{\alpha}\ar[r]^-{\Tr_{\pr_2}} & R\ul{\Gamma}_\Psi(\sO_Y)\ar[d]^{g_Y^*}\\
            R\ul{\Gamma}_{\Psi}Rg_{Y*}R\pr'_{2*}(\pr_2')^!\sO_{Y'}\ar[r]^-{\Tr_{\pr_2'}} &  R\ul{\Gamma}_\Psi Rg_{Y*}(\sO_{Y'}),     }\]
which follows from \cite[VII, Cor. 3.4, (b), TRA 4]{Ha}. Hence the statement.

\end{proof}

\begin{proposition}\label{pf-pb-general}
Let 
\[\xymatrix{(X',\Phi')\ar[r]^{f'}\ar[d]_{g_X} & (Y',\Psi')\ar[d]^{g_Y}\\
            (X,\Phi)\ar[r]^{f}&  (Y,\Psi),   }\]
be a cartesian square with $f,f'\in V_*$ and $g_X, g_Y\in V^*$. Assume either that $g_Y$ is flat or
$g_Y$ is a closed immersion and $f$ is transversal to $Y'$. Then
\[H^*(g_Y)\circ H_*(f)=H_*(f')\circ H^*(g_X).\]
\end{proposition}
\begin{proof}
Embedding $X$ in $X\times Y$ via the graph morphism, the above diagram splits as 
\[\xymatrix{ (X',\Phi')\ar@{^(->}[r]\ar[d]^{g_X} & (X\times Y', \Phi')\ar[r]^{\pr_2}\ar[d]^{\id\times g_Y} & (Y',\Psi')\ar[d]^{g_Y}\\
              (X,\Phi)\ar@{^(->}[r]   & (X\times Y, \Phi)\ar[r]^{\pr_2} & (Y,\Psi).}\]
Both squares are cartesian, the projections $\pr_2$ are smooth and the inclusions are closed. If $g_Y$ is a closed immersion and $f$ is 
transversal to $Y'$, then $\id\times g_Y: X\times Y'\to X\times Y$ is transversal to $X\inj X\times Y$. Thus the statement follows from
Proposition \ref{properties-pfwithsupports}, (2),  Corollary \ref{pf-pb-cl-emb} and Lemma \ref{Hodgepf-cartesianproduct}.
\end{proof}

\begin{lemma}\label{pf-divisor}
 Let $X$ be smooth and $\imath: D\inj X$ the inclusion of a smooth divisor. Let $\Phi$ be a family of supports on $D$ and denote by $\imath_1: (D,\Phi)\to (X,\Phi)$ the map in $V_*$ induced by 
$\imath$. Then
$H_*(\imath_1): H^i_\Phi(D,\Omega^{j}_D)\to H^{i+1}_\Phi(X,\Omega^{j+1}_X)$ is the connecting homomorphism of the 
long exact cohomology sequence associated to the following exact sequence
\eq{Res-exact-seq}{0\to \Omega^{j+1}_X\to \Omega^{j+1}_X(\log D)\xr{\rm Res} \imath_*\Omega^j_D\to 0,}
where $Res(\frac{dt}{t}\alpha)=\imath^*(\alpha)$, for $t\in \sO_X$ a regular element defining $D$ and $\alpha\in \Omega^j_X$. 
In particular if $\Phi\subset X$ is supported in codimension $\ge i+1$ in $X$, then $H_*(\imath_1)$ is injective on $H^i_\Phi$.
\end{lemma}

\begin{proof}
By remark \ref{restrictingtosmoothlocus}, (1), the map $\imath_*$ from Definition \ref{defpush-forward} induces a map, which  we denote by $\imath_*$ again,
\[\imath_* : \imath_*\Omega^j_D\to \Omega^{j+1}_X[1].\]
It suffices to show that this map coincides with the edge homomorphism coming from the distinguished triangle \eqref{Res-exact-seq}, which we denote by $\partial_{\Res}$.
The diagram
\[\xymatrix{\sHom(\Omega^{n-(j+1)}_X, \imath_*\omega_D)\ar[r] & \sHom(\Omega^{n-(j+1)}_X, \omega_X[1])\\
             \imath_*\Omega^{j-1}_D\ar[r]\ar[u] & \Omega^j_X[1],\ar[u]_\simeq   }\]
where $n=\dim X$ and the vertical maps are induced by multiplication from the left, is commutative for both  $\imath_*$ and $\partial_{\Res}$. Thus we only need to consider the case $j=n-1$.

Let $K^\bullet$ be the complex $\OO_X(-D)\to\OO_X$ in degree $[-1,0]$. Then
$K^\bullet\to \imath_*\OO_D$ is a locally free resolution.
We denote by $\Tr_\imath'$ the composition $\imath_*\omega_{D/X}[-1]\xr{\eta_\imath} \imath_*\imath^!\sO_X\xr{\Tr_\imath} \sO_X$, where $\eta_\imath$ is the fundamental local isomorphism (see \eqref{FLI}).
Then $\Tr_\imath'$ is given by 
\[\Tr_\imath': \imath_*\omega_{D/X}[-1]\xl{\simeq}\sHom^\bullet(K^\bullet,\sO_X)\to \sO_X.\]
Here the first map is in degree 1 given by (see \eqref{FLI-explicit})
\[\sHom(\sO_X(-D),\sO_X)=\sO_X(D)\to \imath_*\omega_{D/X} , \quad \frac{1}{t}\mapsto -t^\vee ,\]
where $t$ is a regular parameter defining $D$,
and  the second map (in degree 0) by $\sHom(\sO_X,\sO_X)=\sO_X$. (See the proof of Lemma \ref{trace-for-lci-appendix} and in particular \eqref{compute-Rhom}.)

It thus follows from the commutative diagram \eqref{1.7.11.1} that $\imath_*: \imath_*\omega_D\to \omega_X$ equals the composition
\[\imath_*\omega_D\to \imath_*\omega_{D/X}\otimes \omega_X\xr{\Tr_\imath'[1]\otimes\id} \omega_X[1],\] 
where the first map is given by  $\alpha\mapsto t^\vee\otimes (dt\wedge\tilde{\alpha})$, with $\tilde{\alpha}$ a lift.
Obviously the following diagram commutes
\[\xymatrix{ 0\ar[r] &  \omega_X\ar[r]\ar@{=}[d] & \omega_X(\log D)\ar[d]\ar[r]^{\rm -Res} & \imath_*\omega_D\ar[d]\ar[r] & 0\\ 
           0\ar[r] & \sHom(\OO_X,\OO_X)\otimes\omega_X\ar[r] &  \sO_X(D)\otimes\omega_X\ar[r]^{-1} &  \imath_*\omega_{D/X}\otimes \omega_X\ar[r] & 0. }\]
And by the above (and the sign conventions from \cite[1.3]{Co}) the map $\imath_*\omega_{D/X}\otimes \omega_X\to\sHom(\OO_X,\OO_X)[1]\otimes\omega_X $ induced by the
lower exact sequence equals $-(\Tr_\imath'\otimes \id)$. (Here we need that $\sHom^\bullet(K^\bullet, \sO_X)[1]= \sHom(K^{-(\bullet+1)},\sO_X)$.)
Thus the commutativity of the above diagram yields $i_*= -\partial_{-\Res}=\partial_\Res$. 
\end{proof}

\subsection{The K\"unneth morphism}
For $(X,\Phi)$ and $(Y,\Psi)\in {\rm ob}(V_*)={\rm ob}(V^*)$ the K\"unneth morphism 
\eq{defkuenneth}{\times:H^i_\Phi(X, \Omega^p_X)\times H^j_\Psi(Y, \Omega^q_Y)\to H^{i+j}_{\Phi\times \Psi}(X\times Y, \Omega^{p+q}_{X\times Y}) }
is defined as the composition of the cartesian product with multiplication. 
Choose flasque resolutions $\Omega^p_X\to I^\bullet$, $\Omega^q_Y\to J^\bullet$  and denote  
$K^i_\Phi={\rm Ker} (\Gamma_\Phi I^i \to \Gamma_\Phi I^{i+1})$ and $K^j_\Psi={\rm Ker}(\Gamma_\Psi J^j \to \Gamma_\Psi J^{j+1} )$.  
Then  $\pr_1^{-1}I^\bullet\otimes_k \pr^{-1}_2 J^\bullet$ is a resolution of 
$\pr_1^{-1}\Omega^p_X\otimes_k \pr^{-1}_2 \Omega^q_Y$ and \eqref{defkuenneth} is induced by the composition of the natural maps
\mlnl{K^i_\Phi \otimes_k K^j_\Psi \to H^{i+j}_{\Phi\times \Psi}(\pr_1^{-1}I^\bullet\otimes_k \pr^{-1}_2 J^\bullet)\to \\
        H^{i+j}_{\Phi\times \Psi}(X\times Y, \pr_1^{-1}\Omega^p_X\otimes_k \pr^{-1}_2 \Omega^q_Y) \to 
          H^{i+j}_{\Phi\times \Psi}(X\times Y, \Omega^{p+q}_{X\times Y}).}
We define
\begin{equation}\label{KT}
T: H(X,\Phi)\otimes H(Y,\Psi) \xr{} H(X\times Y,\Phi\times \Psi)
\end{equation}
by the formula
$$
T(\alpha_{i,p}\otimes \beta_{j,q})=(-1)^{(i+p)\cdot j}(\alpha_{i,p}\times \beta_{j,q}),
$$
where $\alpha_{i,p}\in H^i_\Phi(X, \Omega^p_X), \beta_{j,q}\in H^j_\Psi(Y, \Omega^q_Y),$
and $\times$ is the map in \ref{defkuenneth}.

\begin{proposition}\label{kuenneth}
The triples $(H_*, T, e)$ and $(H^*,T, e)$  define right-lax symmetric monoidal functors (see \ref{triplescond}).
\end{proposition}

We will need the following lemma.

\begin{lemma}\label{derivedprojformula}
 Let $f: X\to Y$ be a morphism, assume $Y$ to be smooth and $X$ of pure dimension $d$. Then for any $p, q\ge 0$ there is a morphism
\[\mu: D_X(\Omega^{d-p}_X)\otimes f^*\Omega^q_Y \lra D_X(\Omega^{d-(p+q)}_X),\]
such that: 
\begin{enumerate}
 \item If $U\subset X$ is a smooth open subset, then the following diagram commutes
\[\xymatrix{ D_U(\Omega^{d-p}_U)[-d]\otimes {f_{|U}}^*\Omega^q_Y \ar[r]^-{\mu_{|U}[-d]}\ar[d]^{\simeq}_{\eqref{isoD(Omega)}} & D_U(\Omega^{d-(q+p)}_U)[-d]\ar[d]^{\simeq}_{\eqref{isoD(Omega)}}\\
               \Omega^p_U\otimes {f_{|U}}^*\Omega^q_Y\ar[r] & \Omega^{p+q}_U, }\]
where the lower horizontal map is given by $\alpha\otimes \beta\mapsto \alpha\wedge f^*(\beta)$.
\item If $f$ is proper, then the following diagram commutes 
\[\xymatrix{ Rf_*D_X(\Omega^{d-p}_X)\otimes \Omega^q_Y\ar[r]^-\simeq\ar[dr]_{f_*\otimes \id} & 
                                  Rf_*( D_X(\Omega^{d-p}_X)\otimes f^*\Omega^q_Y)\ar[r]^-\mu &  Rf_*D_X(\Omega^{d-(p+q)}_X)\ar[d]^{f_*}\\
               & D_Y(\Omega^{d-p}_Y)\otimes \Omega^q_Y \ar[r]   & D_Y(\Omega^{d-(p+q)}_Y),
                  }\]
where the lower horizontal map is induced by, $\sHom(\Omega^{d-p}_Y, \Omega^{d_Y}_Y)\otimes\Omega^q_Y\to\sHom(\Omega^{d-(q+p)}_Y, \Omega^{d_Y}_Y),$ $ \varphi\otimes \alpha \mapsto \varphi(\alpha\wedge(-))$. 
\end{enumerate}

\end{lemma}

\begin{proof}
 We denote by $\pi_X$ and $\pi_Y$ the structure maps of $X$ and $Y$, respectively. Since $\pi_X^!k$ and $\pi_Y^!k$ are dualizing complexes, 
 they are represented by bounded complexes of injectives $I_X^\bullet$, $I_Y^\bullet$ and 
$\Tr_f: f_*\pi_X^!k\cong f_*f^!\pi_Y^!k \to \pi_Y^!k$ is thus represented by a morphism of complexes $\Tr_f: f_*I_X^\bullet\to I_Y^\bullet$.
Now the map 
\[\mu : \sHom_X(\Omega^{d-p}_X, I_X^\bullet)\otimes f^*\Omega^q_Y \lra \sHom_X(\Omega^{d-(p+q)}_X, I_Y^\bullet)\]
is in degree $n$ given by 
\[\sHom_X(\Omega^{d-p}_X, I_X^n)\otimes f^*\Omega^q_Y\lra \sHom_X(\Omega^{d-(p+q)}_X, I_Y^n), \quad 
   \theta\otimes \alpha \mapsto \theta(f^*(\alpha)\wedge -).\]
It is immediate that this defines a map of complexes which satisfies (1). For (2) we observe, that it suffices to check the commutativity of
\[\xymatrix{\sHom_X(f_*\Omega^{d-p}_X, f_*I_X^\bullet)\otimes \Omega^q_Y \ar[d]_{\Tr_f\circ(-)\circ f^*}\ar[r]^\mu & 
                                       \sHom_X(f_*\Omega^{d-(p+q)}_X, f_*I_X^\bullet)\ar[d]^{\Tr_f\circ(-)\circ f^*}\\
            \sHom_Y(\Omega^{d-p}_Y, I_Y^\bullet)\otimes \Omega^q_Y\ar[r] &  \sHom(\Omega^{d-(p+q)}_Y, I_Y^\bullet),  }\]
which is straightforward. 

\end{proof}

\begin{proof}[Proof of Proposition \ref{kuenneth}.]
Recall that $H^*(X,\Phi)$ is graded by \ref{Hodgeuppergrading} and 
$H_*(X,\Phi)$ is graded by \ref{Hodgelowergrading}. 
The morphism $T$ respects the grading for both gradings. 
In the following we will work with the upper grading $H^*$. All arguments
will also work for the lower grading $H_*$ because the difference 
between lower and upper grading is an even integer.

By using the associativity of $\times$ \eqref{defkuenneth} it  is 
straightforward to prove the associativity of $T$. Let us prove the commutativity 
of $T$, i.e. that the diagram 
\begin{equation}
\label{proof-Tcomm}
\xymatrix
{
H(X,\Phi)\otimes H(Y,\Psi) \ar[r]^{T} \ar[d]
&
H(X\times Y,\Phi\times \Psi) \ar[d]
\\
H(Y,\Psi)\otimes H(X,\Phi) \ar[r]^{T}
&
H(Y\times X,\Psi\times \Phi)
}
\end{equation}
is commutative. The left vertical map is defined by $a\otimes b\mapsto 
(-1)^{\deg(a)\deg(b)} b\otimes a$, and the right vertical map is given 
by $H^*(\epsilon_1)$ and $H_*(\epsilon_2)$, respectively, with 
\begin{align*}
\epsilon_1&:(Y\times X,\Psi\times \Phi) \xr{} (X\times Y,\Phi\times \Psi) \\
\epsilon_2&: (X\times Y,\Phi\times \Psi) \xr{} (Y\times X,\Psi\times \Phi)  
\end{align*}
the obvious morphisms $\epsilon_1\in V^*$ and $\epsilon_2\in V_*$. Obviously,
$H^*(\epsilon_1) = H_*(\epsilon_2)$, thus we may work with $H^*(\epsilon_1)$
in the following. Note that the diagram 
$$
\xymatrix
{
H^i_\Phi(X, \Omega^p_X)\times H^j_\Psi(Y, \Omega^q_Y)
\ar[r]^{\text{\eqref{defkuenneth}}} \ar[d]
&
H^{i+j}_{\Phi\times \Psi}(X\times Y, \Omega^{p+q}_{X\times Y}) \ar[d]^{(-1)^{p\cdot q}H^*(\epsilon_1)} 
\\
H^j_\Psi(Y, \Omega^q_Y) \times H^i_\Phi(X, \Omega^p_X)
\ar[r]^{\text{\eqref{defkuenneth}}}
&
H^{i+j}_{\Psi\times \Phi}(Y\times X, \Omega^{p+q}_{Y\times X})
}
$$
is commutative, where the left vertical arrow is defined by 
$a\times b\mapsto (-1)^{i\cdot j} (b\times a)$. By using this diagram it 
is a straightforward calculation to prove the commutativity of \eqref{proof-Tcomm}.

We still need to prove the functoriality of $T$ for $H^*$ and $H_*$. For $H^*$
this follows immediately from the definitions. Let us prove the functoriality
for $H_*$.  
We will write $H_\Phi^i(-)$ instead of $H^i_\Phi(X, -)$. By using the 
commutativity of $T$ \eqref{proof-Tcomm} it is enough to prove that the following diagram commutes
\begin{equation}\label{proof-TandH}
\xymatrix{ H^i_\Phi(\Omega^p_X) \times H^j_\Psi(\Omega^q_Y)\ar[d]_{H_*(h)\times \id}\ar[r]^T & 
                                H^{i+j}_{\Phi\times\Psi}(\Omega^{p+q}_{X\times Y})\ar[d]^{H_*(h\times\id)} \\
             H^{i-r}_{\Phi'}(\Omega^{p-r}_{X'}) \times H^j_{\Psi}(\Omega^q_Y)\ar[r]^T & 
                                H^{i+j-r}_{\Phi'\times\Psi}(\Omega^{p+q-r}_{X'\times Y}), }
\end{equation}
 for any $(Y,\Psi)\in V_*$ and $h:(X,\Phi)\to (X',\Phi')$ in $V_*$ and $r=\dim X-\dim X'$ ($X$ and $X'$ are assumed to be equidimensional). Equivalently, 
the diagram as in \eqref{proof-TandH}, but with $\times$ instead of $T$
as horizontal arrows, commutes.
Observe that  $\times$ can be factored as
\eq{altT}{\times: H^i_\Phi(\Omega^p_X) \times H^j_\Psi(\Omega^q_Y)\xr{ H^*(\pr_1)\times\id}
                 H^i_{\Phi\times Y}(\Omega^p_{X\times Y}) \times H^j_{\Psi}(\Omega^q_Y) \to H^{i+j}_{\Phi\times\Psi}(\Omega^{p+q}_{X\times Y}), }
where the map on the right is the composition of the cartesian product with the multiplication map
$\Omega^p_{X\times Y}\otimes_k \pr_2^{-1}\Omega^q_Y\to \Omega^{p+q}_{X\times Y}$, $\alpha\otimes \beta\mapsto \alpha\wedge \pr_2^*\beta$.
By Proposition \ref{pf-pb-general} the diagram  
\eq{pbT}{\xymatrix{ H^i_\Phi(\Omega^p_X)\ar[d]_{H(h_*)}\ar[r]^-{H^*(\pr_1)} & 
                  H^i_{\Phi\times Y}(\Omega^p_{X\times Y})\ar[d]^{H_*(h \times \id)}\\
             H^{i-r}_{\Phi'}(\Omega^{p-r}_{X'})\ar[r]^-{H^*(\pr_1)} &  H^{i-r}_{\Phi'\times Y}(\Omega^{p-r}_{X'\times Y})   }}
commutes. Thus it suffices to prove that
\[\xymatrix{ H^i_{\Phi\times Y}(\Omega^p_{X\times Y})\times H^j_\Psi(\Omega^q_Y)\ar[d]_{H_*(h\times \id )\times\id}\ar[r] & 
                                H^{i+j}_{\Phi\times\Psi}(\Omega^{p+q}_{X\times Y})\ar[d]^{H_*(h\times \id)} \\
             H^{i-r}_{\Phi'\times Y}(\Omega^{p-r}_{X'\times Y})\times H^j_\Psi(\Omega^q_Y)\ar[r] & 
                                H^{i+j-r}_{\Phi'\times\Psi}(\Omega^{p+q-r}_{X'\times Y})}\]
commutes.

 Now let $\bar{h}: \bar{X}\to X'$ be a compactification of $h$ and set $d=\dim X+ \dim Y$. We write 
\[\omega_{\bar{X}\times Y}^p:= D_{\bar{X}\times Y}(\Omega^{d-p}_{\bar{X}\times Y}), \quad 
                                                                     \omega^p_{X'\times Y}:=D_{X'\times Y}(\Omega^{d-p}_{X'\times Y}).\]
Notice that ${\omega_{\bar{X}\times Y}^p}_{|X\times Y}\cong \Omega^p_{X\times Y}[d]$ and
$\omega_{X'\times Y}^p\cong \Omega^{p-r}_{X'\times Y}[d-r]$. With this notation the push-forward is  a morphism 
\[(\bar{h}\times \id)_*: R(\bar{h}\times \id)_*\omega^p_{\bar{X}\times Y}\to \omega^p_{X'\times Y}\]
and we have to show that the following diagram commutes:
\[\xymatrix{ H^i_{\Phi\times Y}(\omega^p_{\bar{X}\times Y})\times H^j_\Psi(\Omega^p_Y)\ar[d]_{(h\times\id)_*\times\id}\ar[r] & 
                                H^{i+j}_{\Phi\times\Psi}(\omega^{p+q}_{\bar{X}\times Y})\ar[d]^{(h\times\id)_*} \\
             H^i_{\Phi'\times Y}(\omega^p_{X'\times Y})\times H^j_\Psi(\Omega^p_Y)\ar[r] & 
                                H^{i+j}_{\Phi'\times\Psi}(\omega^{p+q}_{X'\times Y}),}\]
where the upper map is given by the cartesian product composed with the $\mu$ from Lemma \ref{derivedprojformula}. 
Clearly we may assume $\Phi=\bar{h}^{-1}(\Phi')$; thus 
\eq{psitopsi'}{H^i_{\Phi\times Y}(\omega^p_{\bar{X}\times Y})=H^i_{\Phi'\times Y}(R(\bar{h}\times\id)_*\omega^p_{\bar{X}\times Y}).}
Now it follows from Lemma \ref{derivedprojformula}, (2), that it is enough to prove the commutativity of the following two diagrams
\[\xymatrix{H^i_{\Phi\times Y}(\omega^p_{\bar{X}\times Y})\times H^j_\Psi(\Omega^q_Y)  \ar[d]_{(\bar{h}\times \id)_*\times\id}\ar[r] &
  H^{i+j}_{\Phi'\times\Psi}(R(\bar{h}\times\id)_*(\omega^p_{\bar{X}\times Y})\otimes_k \pr_2^{-1}\Omega^q_{Y})
                                                                                \ar[d]^{(\bar{h}\times\id)_*\otimes \id}\\
  H^i_{\Phi'\times Y}(\omega^p_{X'\times Y})\times H^j_\Psi(\Omega^q_Y)\ar[r] &
  H^{i+j}_{\Phi'\times\Psi}(\omega^p_{X'\times Y}\otimes_k \pr_2^{-1}\Omega^q_Y), }\]
here the upper horizontal map is the composition of \eqref{psitopsi'} with the cartesian product, and
\[\xymatrix{H^{i+j}_{\Phi'\times\Psi}(R(\bar{h}\times\id)_*(\omega^p_{\bar{X}\times Y})\otimes_k \pr_2^{-1}\Omega^q_{Y})
                                                                                     \ar[d]_{(\bar{h}\times\id)_*\otimes\id}\ar[r] &
   H^{i+j}_{\Phi'\times\Psi}(R(\bar{h}\times\id)_*(\omega^p_{\bar{X}\times Y})\otimes_\sO \pr_2^*\Omega^q_{Y}) 
                                                                                            \ar[d]^{(\bar{h}\times\id)_*\otimes\id}\\
     H^{i+j}_{\Phi'\times\Psi}(\omega^p_{X'\times Y}\otimes_k \pr_2^{-1}\Omega^q_Y)\ar[r] &
   H^{i+j}_{\Phi'\times\Psi}(\omega^p_{X'\times Y}\otimes_\sO \pr_2^*\Omega^q_Y). }\]
For this take injective resolutions $\omega^p_{\bar{X}\times Y}\to I^\bullet$ and $\omega^p_{X'\times Y}\to J^\bullet$, then
the push-forward is given by an actual morphism $(\bar{h}\times\id)_* I^\bullet\to J^\bullet$. Now the commutativity of the first diagram
is easily checked by taking an injective resolution of $\Omega^q_Y$. For the commutativity of the second diagram we observe, that 
$(\bar{h}\times\id)_*I^\bullet\otimes_k \pr^{-1}_2\Omega^q_Y$  and $(\bar{h}\times \id)_*I^\bullet\otimes_\sO\pr^*_2\Omega^q_Y$ still
represent $R(\bar{h}\times\id)_*\omega^p_{\bar{X}\times Y}\otimes_k \pr^{-1}_2\Omega^q_Y$  and 
$R(\bar{h}\times\id)_*\omega^p_{\bar{X}\times Y}\otimes_\sO \pr^*_2\Omega^q_Y$, respectively and similar with $(\bar{h}\times \id)_*I^\bullet$ 
replaced by $J^\bullet$ and $(\bar{h}\times\id)_*\omega^p_{\bar{X}\times Y}$ by $\omega^p_{X'\times Y}$. Thus it is enough to check the commutativity
using these complexes, which is obvious.
\end{proof}

\subsection{Summary}

Let $(H_*, H^*, T, e)$ be the datum defined above, i.e. $H_*: V_*\to \GrAb$ is defined on objects by \eqref{Hodgelowergrading} and on morphisms by Definition \ref{pfwithsupport},
$H^*:(V^*)^{op}\to \GrAb$ is defined on objects by \eqref{Hodgeuppergrading} and on morphisms by \eqref{pullback}, $T$ is defined by \eqref{KT} and $e$ by \eqref{e}.

\begin{thm}\label{hodgeistriple}
The datum $(H_*, H^*, T, e)$ is an object in $\T$, i.e. it is a datum as in \ref{triples} and satisfies the properties \ref{triplescond}.
\end{thm}

We denote by $HP$ the pure part of $H$, i.e.
\eq{purehodgewithsupp}{HP(X,\Phi):=\bigoplus_{n\ge 0} H^n_\Phi(X,\Omega^n_X), \quad (X,\Phi)\in V_*,}
 and let $HP^*(X,\Phi)$ be the graded abelian group, which in degree $2n$ equals
\[HP^{2n}(X,\Phi)= H^n_{\Phi}(X, \Omega^n_X)\]
and is zero in odd degrees. The graded abelian group $HP_*(X,\Phi)$ is defined as in \eqref{Hodgelowergrading}.
Then the functors $H_*$, $H^*$ induce functors $HP_*$, $HP^*$ and $T$ and $e$ restrict to $HP$. We obtain:

\begin{corollary}\label{purehogdeistriple}
The datum $(HP_*, HP^*, T, e)$ is an object in $\T$.
Furthermore $HP$ satisfies the semi-purity condition from Definition \ref{semipurity} and the natural map
$(HP_*, HP^*,T, e)\to (H_*, H^*, T, e)$ is a morphism in  $\T$.
\end{corollary}
\begin{proof}
The semi-purity condition follows from \cite[Exp.III, Prop. 3.3]{SGA2}.
\end{proof}

\section{Cycle class map to Hodge cohomology and applications}

{\em In this section $k$ is assumed to be a perfect field (unless stated otherwise).}

\subsection{Cycle class} 

\begin{proposition}[Cycle class]\label{cycleclassmap}
Let $X$ be a smooth scheme and $W\subset X$ an irreducible closed
subset of codimension $c$. There is a class $cl_{(X,W)}=cl(W)\in H^c_W(X,\Omega_X^c)$ with the property  
$$
H^*(\jmath)(cl(W))=H_*(\imath_{U\cap W})(1)
$$
for every open subset $U\subset X$ such that $U\cap W$ is smooth (and non-empty), where
$\jmath:(U,W\cap U) \xr{} (X,W)$ is induced by the open immersion, 
$\imath:W\cap U \xr{} (U,W\cap U)$ is induced by the closed immersion, and $1$ is the identity element of the ring $H^0(X, \sO_X)$.
\end{proposition} 

\begin{remark}
The cycle class in the above proposition is Grothendieck's ``fundamental class'' (see e.g. \cite[p.39, (ii)]{Lipman}).
For the convenience of the reader and to be sure about the compatibility with the push-forward constructed in the previous section,
we give a proof of the proposition, which is standard.
\end{remark}

\begin{proof}
\emph{1st Step:} Let $\eta$ be the generic point of $W$. We define
$$
H_{\eta}^c(X,\Omega^c_X):=\varinjlim_{\eta\in U} H^c_{U\cap W}(U,\Omega_U^c)
$$
where the inductive limit runs over all open sets $U\subset X$ with $\eta\in U$. 
Choose $U$ such that $U\cap W\neq \emptyset$ is smooth. The image of 
$H_*(\imath_{U\cap W})(1)\in H^c_{U\cap W}(U,\Omega_U^c)$ in $H_{\eta}^c(X,\Omega^c_X)$ doesn't 
depend on the choice of $U$ by \ref{triplescond}(\ref{MP}). We denote this class by $cl(W)_{\eta}$.

\emph{2nd Step:} A class $a\in H_{\eta}^c(X,\Omega^c_X)$ is in the image of 
$$
H_W^c(X,\Omega_X^c) \xr{} H_{\eta}^c(X,\Omega^c_X)
$$
(i.e. extends to a global class) if and only if for all  $1$-codimensional points $x$ in $W$ there is an open subset 
$U\subset X$ containing $x$, so that $a$ lies in the image of 
$$
H_{W\cap U}^c(U,\Omega_U^c) \xr{} H_{\eta}^c(X,\Omega^c_X).
$$
Indeed, the Cousin resolution yields an exact sequence 
\begin{equation}\label{Cousin}
0\xr{} H_W^c(X,\Omega_X^c) \xr{} H_{\eta}^c(X,\Omega^c_X) \xr{} \bigoplus_{x\in W, cd(x)=1}
H_x^{c+1}(X,\Omega_X^c),
\end{equation}
and $H_{U\cap W}^c(U,\Omega_U^c)\xr{} H_{\eta}^c(X,\Omega_X^c)\xr{} H_x^{c+1}(X,\Omega_X^c)$ 
vanishes for all $x$ and $U$ as above.

\emph{3rd Step:} If $W$ is normal then $cl(W)_{\eta}$ extends (uniquely)
to a class in $H^c_W(X,\Omega_X^c)$. Indeed,  since $W$ is regular in codimension $=1$
and we assume that $k$ is perfect, we may choose an open $U\subset X$ such that $U\cap W$
is smooth and $U\cap W$ contains all points of codimension $=1$ of $W$. So that 
the class extends by the 2nd Step. Note that the extension is unique because of the 
exact sequence \ref{Cousin}. 

\emph{4th Step:} We claim that the class $cl(W)_{\eta}$ extends to a class in 
$H^c_W(X,\Omega_X^c)$. In view of the 2nd Step it is sufficient to extend the class
at all points $x\in W$ of codimension $=1$. Thus we may assume that $X$ (and therefore $W$)
is affine. The normalisation $\ti{W}\xr{} W$ is a finite morphism and thus projective. 
Choose an embedding $\ti{W}\xr{} W\times_k \P^n_k$ over $W$. The previous step yields 
a class $cl(\ti{W})\in H^{n+c}_{\ti{W}}(X\times \P^n,\Omega_{X\times \P^n}^{n+c})$. Consider
$H_*({\rm pr}_1)(cl(\ti{W}))\in H^{c}_W(X,\Omega_X^c)$; for an open $U\subset X$ such that
$W\cap U\neq \emptyset$ and $U\cap W$ is smooth  we obtain
\begin{multline*}
H^*(\jmath)H_*({\rm pr}_1)(cl(\ti{W}))=H_*({\rm pr}_{1 \mid U\times \P^n})H^*(\jmath')(cl(\ti{W}))
= \\
H_*({\rm pr}_{1 \mid U\times \P^n})H_*(\imath_{(U\times \P^n)\cap \ti{W}})(1) = 
H_*(\imath_{U\cap W})(1),
\end{multline*}
with $\jmath':(U\times \P^n,(U\times \P^n)\cap \ti{W})\xr{} (X\times \P^n,\ti{W})$.
Thus $H_*({\rm pr}_1)(cl(\ti{W}))$ is the desired lift.
\end{proof}

\subsubsection{Explicit description of the cycle class}\label{classexplicit}
Let $X$ be a smooth scheme and $W\subset X$ an irreducible closed subset of codimension $c$ with generic point $\eta\in X$.
Denote $A=\OO_{X,\eta}$. Then
$$ H^c_\eta(X, \Omega^c_{X})=\varinjlim_{f\subset \mathfrak{m}_\eta }\frac{\Omega^c_A}{(f)},$$
where the limit is over all $A$-sequences $f=(f_1,\ldots, f_c)$ of length $c$ which are contained
in $\mathfrak{m}_\eta$ (in particular $\sqrt{(f_1,\ldots, f_c)}=\mathfrak{m}_\eta$). 
The class of $\omega\in \Omega^c_A$ under the composition  in $\Omega^c_A\to\Omega^c_A/(f)\to H^c_\eta(X, \Omega^c_{X})$ is denoted by
$\genfrac[]{0pt}{}{\omega}{f}$. See \ref{A1} for details. 

Now let $U$ be an affine open subset of $X$ such that $U\cap W$ is smooth and the ideal of $W\cap U$ in $\OO_U$ is generated 
by global sections $t_1,\ldots,t_c$ on $U$. Then by Proposition \ref{1.7.9}
$$cl(W)_\eta =(-1)^c\genfrac[]{0pt}{}{dt_1\cdots dt_c}{t_1,\ldots, t_c}.$$

\begin{lemma}\label{ccmpbcd1} 
For a closed immersion $\imath:X\xr{} Y$ between smooth schemes
and an effective smooth divisor $D\subset Y$ such that 
\begin{itemize}
\item $D$ meets $X$ properly, thus $D\cap X:=D\times_{Y} X$ is a divisor on X,
\item $D':=(D\cap X)_{red}$ is smooth and connected, and thus $D\cap X=n\cdot D'$
as divisors (for some $n\in \Z, n\geq 1$), 
\end{itemize}
we denote by $\imath_X:X\xr{}(Y,X), \imath_{D'}:D'\xr{} (D,D')$ 
the morphisms in $V_*$ induced by $\imath$, and we define $g_2:(D,D')\xr{} (Y,X)$ in $V^*$ by  
the inclusion $g:D\inj Y$. Then the following equality holds:
$$
H^*(g_2)(H_*(\imath_X)(1_X))=n\cdot H_*(\imath_{D'})(1_{D'}).
$$
\begin{proof}
Let $c$ be the codimension of $X$ in $Y$ and $g_3:(D,D')\xr{} (Y,D')$ be induced by the inclusion $D\subset Y$. 
Then 
$$
H_*(g_3): H^c_{D'}(D,\Omega_D^c)\xr{} H^{c+1}_{D'}(Y,\Omega_Y^{c+1}) 
$$
is injective (by Lemma \ref{pf-divisor}) and thus we need to prove
\[H_*(g_3)H^*(g_2)H_*(\imath_X)(1_X)=n\cdot H_*(g_3)H_*(\imath_{D'})(1_{D'}).\]
Let  $g_1: D\xr{} (Y,D)$ be induced by $g$, then projection formula \ref{projectionformula} gives
\[H_*(g_3)H^*(g_2)(H_*(\imath_X)(1_X))= H_*(g_1)(1_D)\cup H_*(\imath_X)(1_X).\]
Therefore it suffices to prove
\eq{ccmpbcd2}{H_*(g_1)(1_D)\cup H_*(\imath_X)(1_X)= n\cdot H_*(g_3\circ \imath_{D'})(1_{D'}).}
Let $\eta$ be the generic point of $D'$. Since  $H^{c+1}_Z(Y,\Omega^{c+1}_Y)=0$
for all closed subsets $Z\subset Y$ of codimension $\ge c+2$, by \cite[III, Prop. 3.3]{SGA2}, the restriction map
\begin{equation}\label{loc100}
H^{c+1}_{D'}(Y,\Omega_Y^{c+1}) \xr{} H^{c+1}_{\eta}(Y,\Omega_Y^{c+1})
\end{equation}
is injective. Thus it is sufficient to prove the equality \eqref{ccmpbcd2} in $H^{c+1}_{\eta}(Y,\Omega_Y^{c+1})$.

Since $X$ is smooth we may find a regular sequence  
$t_1,\dots,t_c\in \OO_{Y,\eta}$, which generates the ideal of $X$. 
If $D={\rm div}(f)$ around $\eta$ then
$$
(-1)\genfrac[]{0pt}{}{df}{f}\cup (-1)^c\genfrac[]{0pt}{}{dt_1\wedge \dots \wedge dt_c}{t_1,\dots,t_c}
=(-1)^{c+1}\genfrac[]{0pt}{}{df\wedge dt_1\wedge \dots \wedge dt_c}{f,t_1,\dots,t_c} 
$$
is the image of $H_*(g_1)(1_D)\cup H_*(\imath_X)(1_X)$ in $H^{c+1}_{\eta}(Y,\Omega_Y^{c+1})$.

Let $\pi\in \OO_{Y,\eta}$ be a lift of a generator of the maximal ideal in $\OO_{X,\eta}$. 
By the explicit description of the cycle class in \ref{classexplicit} we get 
$$
H_*(g_3\circ \imath_{D'})(1_{D'})=(-1)^{c+1}\genfrac[]{0pt}{}{d\pi\wedge dt_1\wedge \dots \wedge dt_c}{\pi,t_1,\dots,t_c}.
$$
Obviously $f=a\pi^n$ in $\OO_{X,\eta}$ for a unit $a\in \OO_{X,\eta}^*$. Choose a lift 
$\tilde{a}\in \OO_{Y,\eta}^*$ of $a$, thus $f=\tilde{a}\pi^n$ modulo $(t_1,\dots,t_c)$, and we obtain
\[
\genfrac[]{0pt}{}{df\wedge dt_1\wedge \dots \wedge dt_c}{f,t_1,\dots,t_c}=
\genfrac[]{0pt}{}{n \tilde{a}\pi^{n-1}\cdot d\pi\wedge dt_1\wedge \dots \wedge dt_c}{\tilde{a}\pi^n,t_1,\dots,t_c}  
=n\cdot \genfrac[]{0pt}{}{d\pi\wedge dt_1\wedge \dots \wedge dt_c}{\pi,t_1,\dots,t_c},
\]
which proves \eqref{ccmpbcd2}.
\end{proof}
\end{lemma}

\begin{thm}\label{propositionCHtoHodge}
There exists a morphism $cl:\CH\xr{}H=(H_*, H^*, T, e)$ in $\T$.
\begin{proof}
Since there is a morphism $HP=(HP_*, HP^*, T,e)\xr{} H$ in $\T$, it suffices to prove the existence of $cl:\CH\xr{} HP$.
 
This follows from Theorem \ref{mainthm1}, since $HP$ satisfies all the conditions listed there:
We have $HP\in\T$ and it satisfies the semi-purity condition \ref{semipurity} by Corollary \ref{purehogdeistriple},
it satisfies \ref{mainthm1}(\ref{traceofone}) by Proposition \ref{properties-pfwithsupports},(3) and (\ref{mainthm1})(\ref{multiplicities}) by Lemma \ref{ccmpbcd1}.
Finally the element $cl_{(X,W)}$ from  (\ref{mainthm1})(\ref{class}) is the cycle class constructed in Proposition \ref{cycleclassmap} and (\ref{mainthm1})(\ref{rationaleq}) is obvious.
\end{proof}
\end{thm}


\subsection{Main theorems}

\subsubsection{}\label{pb-pf-via-corr}
Let $f:(X,\Phi)\xr{} (Y,\Psi)$ be a morphism in $V_*$ or $V^*$. 
In view of Theorem \ref{propositionCHtoHodge}, Section \ref{VtoCor} and Lemma \ref{lemmaCortau}, the morphism 
$$
H_*(f): H_*(X,\Phi) \xr{} H_*(Y,\Psi), \quad \text{resp.}\; H^*(f):H^*(Y,\Psi)\xr{} H^*(X,\Phi) 
$$  
is given by (we write $cl$ instead of $Cor(cl)$)
$$
H_*(f)=\rho_H\circ cl\circ \tau_*^{\CH}(f), \quad \text{resp.}\; H^*(f)=\rho_H\circ cl\circ \tau^*_{\CH}(f).
$$
Thus we may use composition of correspondences in $Cor_{\CH}$ to compute 
$H^*(f)\circ H_*(f), H_*(f)\circ H^*(f),$ etc.

\begin{proposition} \label{propositionvanishing}
Let $X,Y$ be smooth and connected, and let 
\[\alpha\in \Hom_{Cor_{\CH}}(X,Y)^0= CH^{\dim X}(X\times Y, P(\Phi_X,\Phi_Y)).\]
\begin{enumerate}
\item 
If the support of $\alpha$ projects to an $r$-codimensional subset in $Y$, then
the restriction of $\rho_H\circ cl (\alpha)$ to 
$
\bigoplus_{j<r,i} H^i(X,\Omega_X^j)
$ 
vanishes.
\item If the support of $\alpha$ projects to an $r$-codimensional subscheme in $X$, then
       the restriction of $\rho_H\circ cl(\alpha)$ to 
$
\bigoplus_{j\ge \dim X-r+1, i} H^i(X,\Omega_X^j)
$
vanishes.
\end{enumerate}
\begin{proof}
(1). We may assume $\alpha=[V]$ for $V\subset X\times Y$ a closed irreducible subset 
of dimension $\dim(Y)=:d_Y$, with $p_Y(V)\subset Y$ of codimension $r$. We set $d_X=\dim X$. 

By definition of $\rho_H$ (see \ref{VtoCor}) and Lemma \ref{Hodgepf-cartesianproduct}
it is sufficient to prove that for all $0\le q\le r-1$ the 
image of the class $cl(V)$ vanishes via the map 
$$
H^{d_X}_V(\Omega^{d_X}_{X\times Y})\xr{{\rm proj.}} H^{d_X}_V(\pr_1^*\Omega^{d_X-q}_X\otimes \pr_2^*\Omega^q_Y).
$$
To prove this we may also localize at the generic point $\eta$ of $V$ (by \cite[III, Prop. 3.3]{SGA2}).

We write $B=\OO_{X\times Y,\eta}$ and $A=\OO_{Y, p_{Y}(\eta)}$. Now $A$ is a regular local ring
of dimension $r$ and $B$ is formally smooth over $A$. Let $t_1,\ldots, t_r\in A$ be a regular system 
of parameters of $A$. Since 
$B/(1\otimes t_1, \ldots, 1\otimes t_r)$ is a local regular ring 
there exist elements $s_{r+1}, \ldots, s_{d_X}\in B$ such that
$1\otimes t_1, \ldots, 1\otimes t_r, s_{r+1}, \ldots, s_{d_X}$ is a system of regular parameters for $B$.
Thus by the explicit description of the cycle class in \ref{classexplicit} we obtain 
$$
cl(V)_\eta=(-1)^{d_X}
\genfrac{[}{]}{0pt}{}{d(1\otimes t_1)\wedge \cdots \wedge d(1\otimes t_r) \wedge ds_{r+1}\wedge \cdots \wedge ds_{d_X}}{1\otimes t_1,\ldots,  1\otimes t_r, s_{r+1},\ldots ,s_{d_X}}.
$$
This clearly implies the claim.

(2) Let $\alpha=[V]$ be as in (1) and suppose $p_X(V)$ has codimension $r$ in $X$. As above it suffices to prove that for all $0\le q\le r-1$ the 
image of the class $cl(V)$ vanishes under the projection map 
$$
H^{d_X}_V(\Omega^{d_X}_{X\times Y})\xr{{\rm proj.}} H^{d_X}_V(\pr_1^*\Omega^q_X\otimes \pr_2^*\Omega^{d_X-q}_Y).
$$
Write $C=\sO_{X,p_X(\eta)}$.  Then as in (1) we find $\tau_1,\ldots, \tau_r\in C$ and $\sigma_{r+1},\ldots, \sigma_{d_X}\in B$, such that
$\tau_1\otimes 1, \ldots, \tau_r\otimes 1, \sigma_{r+1},\ldots,\sigma_{d_X}$ is a system of regular parameters for $B$.
Thus
$$ 
cl(V)_\eta=(-1)^{d_X}
\genfrac{[}{]}{0pt}{}{d(\tau_1 \otimes 1)\wedge \cdots \wedge d(\tau_r\otimes 1) \wedge d\sigma_{r+1}\wedge \cdots \wedge d\sigma_{d_X}}
                                                                                                 {\tau_1\otimes 1,\ldots,  \tau_r\otimes 1, \sigma_{r+1},\ldots ,\sigma_{d_X}},
$$
which implies the claim.
\end{proof}
\end{proposition}

\subsubsection{}\label{localversion}
Let $S$ be a $k$-scheme and let $f: X\to S$ and $g:Y\to S$ be two integral $S$-schemes, which are smooth over $k$. Let $Z\subset X\times_S Y$ be a closed integral subscheme of
dimension equal to the dimension of $Y$ and such that ${\pr_2}_{|Z}:Z\to Y$ is proper. 
For an open subset $U\subset S$, we denote by $Z_U\subset f^{-1}(U)\times_U g^{-1}(U)$ the pullback of $Z$ over $U$. This gives a correspondence 
$[Z_U]\in \Hom_{Cor_{\CH}}(f^{-1}(U),g^{-1}(U))^0$, which induces a morphism of $k$-vector spaces 
\[\rho_H\circ cl([Z_U]): H^i(f^{-1}(U),\Omega^j_{f^{-1}(U)})\to H^i(g^{-1}(U),\Omega^j_{g^{-1}(U)}),\quad \text{for all } i,j. \]

\begin{proposition}\label{higherdirectimages}
In the above situation,  the collection $\{\rho_H\circ cl([Z_U])\,|\, U\subset Z \text{ open}\}$ induces a morphism of quasi-coherent $\sO_S$-modules
\[\rho_H(Z/S): R^i f_*\Omega^j_X\to R^i g_*\Omega^j_Y,\quad \text{for all }i,j. \]
\end{proposition}

\begin{proof}
We need to show the following statements: 
\begin{enumerate}
\item The maps $\rho_H\circ cl([Z_U])$ are compatible with restriction to open sets.
\item The maps $\rho_H\circ cl([Z_U])$ are $\sO(U)$-linear. 
\end{enumerate}

\emph{For (1):} Let us denote by 
\begin{equation*}
\begin{split}
&\pr_{1,U}: f^{-1}(U)\times g^{-1}(U) \to f^{-1}(U), \quad \text{$\pr_{1,U}\in V^*$,}\\
&\pr_{2,U}: (f^{-1}(U)\times g^{-1}(U),P(\Phi_{f^{-1}(U)}, \Phi_{g^{-1}(U)}) \to g^{-1}(U), \quad \text{$ \pr_{2,U}\in V_*$,}
\end{split}
\end{equation*}
the morphism induced by the projections (see \ref{Phiset} and \ref{definitionP} for the definition of $P(\Phi_{f^{-1}(U)}, \Phi_{g^{-1}(U)})$). 
Let $j: V\inj U$ be an open immersion and denote by 
\begin{equation*}
j_f: f^{-1}(V)\xr{} f^{-1}(U) , \quad j_g: g^{-1}(V)\xr{} g^{-1}(U) 
\end{equation*}
the morphisms in $V^*$ induced by $j$.

We have to show that for all $a\in H^i(f^{-1}(U), \Omega^j_{f^{-1}(U)})$ the following equality 
holds:
\ml{restriction}{H^*(j_g) H_*(\pr_{2,U})(H^*(\pr_{1,U})(a) \cup cl([Z_U]))  \\
                           = H_*(\pr_{2,V})(H^*(\pr_{1,V})(H^*(j_f)(a) \cup cl([Z_V])).  }
As a first step from the left hand side to the right hand side we observe that 
$$
H^*(j_g)H_*(\pr_{2,U})=H_*(\pr_{2,V}') H^*(id_{f^{-1}(U)}\times j_g)
$$
where 
$$
\pr_{2,V}':(f^{-1}(U)\times g^{-1}(V),\Phi) \xr{} g^{-1}(V)
$$
in $V_*$ is induced by the projection and $\Phi:=(id\times j_g)^{-1}P(\Phi_{f^{-1}(U)},\Phi_{g^{-1}(U)})$. Denoting $\pr_{1,U}':f^{-1}(U)\times g^{-1}(V) \xr{} f^{-1}(U)$ as a morphism in $V^*$ we obtain the equality
$$
H^*(id_{f^{-1}(U)}\times j_g)(H^*(\pr_{1,U})(a) \cup cl([Z_U]))= H^*(\pr_{1,U}')(a)\cup cl([Z_V])
$$
in $H(f^{-1}(U)\times g^{-1}(V), \Phi)$; here 
we consider $Z_V\in \Phi$ as a closed subset of $f^{-1}(U)\times g^{-1}(V)$. Next, consider the 
morphisms 
\begin{equation*}
\begin{split}
j_f\times id_{g^{-1}(V)}&: f^{-1}(V)\times g^{-1}(V) \xr{} f^{-1}(U) \times g^{-1}(V), \\
\tau&: (f^{-1}(V)\times g^{-1}(V),Z_V) \xr{} (f^{-1}(U)\times g^{-1}(V), \Phi),\\
id'&: (f^{-1}(V)\times g^{-1}(V),Z_V) \xr{} (f^{-1}(V)\times g^{-1}(V), P(\Phi_{f^{-1}(V)},\Phi_{g^{-1}(V)})),
\end{split}
\end{equation*}
with 
$$
j_f\times id_{g^{-1}(V)}\in V^*, \quad \tau \in V_*,\quad  id' \in V_*,
$$
and where $id'$ is induced by the identity. The projection formula yields
\begin{multline*}
H^*(\pr_{1,U}')(a)\cup cl([Z_V]) = H^*(\pr_{1,U}')(a)\cup cl(\CH_*(\tau)([Z_V]))= \\H_*(\tau)( H^*(j_f\times id_{g^{-1}(V)})H^*(\pr_{1,U}')(a)\cup cl([Z_V])). 
\end{multline*}
Now the equalities 
\begin{eqnarray*}
H_*(\pr_{2,V}')H_*(\tau) & = & H_*(\pr_{2,V})H_{*}(id'), \\
H^*(j_f\times id_{g^{-1}(V)})H^*(\pr_{1,U}') & = & H^*(\pr_{1,V}) H^*(j_f),
\end{eqnarray*}
imply the claim \ref{restriction}. 

\emph{For (2):} It suffices to consider the case $U=S=\Spec\, R$.  The ring homomorphisms  $g^*:R\xr{} H^0(X,\OO_X)$
and $f^*:R\xr{} H^0(Y,\OO_Y)$ induce $R$-module structures on $H(X)$ and $H(Y)$ respectively.
  
We have to prove the following equality for all $r\in R$ and $a\in H^i(X, \Omega^j_X)$:
$$g^*(r)\cup H_*(\pr_2)(H^*(\pr_1)(a)\cup cl([Z]))=H_*(\pr_2)(H^*(\pr_1)(f^*(r)\cup a)\cup cl([Z])).$$
For this, it is enough to show that
\eq{linear2}{H^*(\pr_2)(g^*(r))\cup cl([Z])= H^*(\pr_1)(f^*(r))\cup cl([Z])\quad \text{in } H^d_Z(X\times Y, \Omega^d_{X\times Y}).} 
Choose an open set $U\subset X\times Y$ such that $Z\cap U$ is nonempty and smooth. 
Since the natural map $H^d_Z(X\times Y, \Omega^d_{X\times Y})\to H^d_{Z\cap U}(U, \Omega^d_{U})$ is injective 
it suffices to check \eqref{linear2} on $H^d_{Z\cap U}(U, \Omega^d_{U})$. 
We write $\imath_1:Z\cap U\xr{} (U,Z\cap U)$ in $V_*$ and $\imath_2:Z\cap U \xr{} U$ in $V^*$ for the obvious morphisms. 
By using the projection formula and $cl([Z\cap U])=H_*(\imath_1)(1)$ we reduce to the statement 
$$
H^*(\imath_2)H^*(\pr_2)(g^*(r))=H^*(\imath_2)H^*(\pr_1)(f^*(r)).
$$
This follows from $g\circ \pr_2 \circ \imath_2 = f\circ \pr_1 \circ \imath_2$.
\end{proof}

\begin{definition}\label{propbir}
Two integral schemes $X,Y$ over a base $S$ are called properly birational over $S$ if there 
is an integral scheme $Z$ over $S$ and morphisms over $S$: 
$$
\xymatrix
{
&
Z \ar[dr]^{\text{proper,birational}} \ar[dl]_{\text{proper,birational}}
&
\\
X
&
&
Y.
}
$$
\end{definition}

\begin{thm} \label{thmoverS}
 Let $S$ be a scheme over a perfect field $k$. Let $f:X\xr{} S$ and $g:Y\xr{} S$ 
be integral $S$-schemes, which are smooth over $k$ and properly birational over $S$.
Let $Z$ be an integral scheme together with proper birational morphisms 
$Z\xr{} X, Z\xr{} Y$ such that 
\[\xymatrix@=15pt{   &  Z\ar[dl]\ar[dr]  &   \\
              X\ar[dr]_f     &               & Y\ar[dl]^g  \\
                 &      S  &     } \]
is commutative.
We denote by $Z_0$ be the image of $Z$ in $X\times_S Y$.
Then, for all $i$,  $\rho_H(Z_0/S)$ induces isomorphisms of $\OO_S$-modules ($d=\dim X=\dim Y$)
\[R^if_*\OO_X\xr{\simeq} R^ig_*\OO_Y,\quad R^if_*\Omega^d_X\xr{\simeq} R^ig_*\Omega^d_Y.\]
\end{thm}

\begin{proof}
Recall that $\rho_H(Z_0/S)$ is defined in Proposition \ref{higherdirectimages}
as the sheafication of the maps 
\begin{equation} \label{equation-to-show-thmoverS}
\rho_H\circ cl([(Z_{0,U}]): H^i(f^{-1}(U),\Omega^j_{f^{-1}(U)})\to H^i(g^{-1}(U),\Omega^j_{g^{-1}(U)}), 
\end{equation}
where $U$ runs over all open sets of $S$ and $Z_{0,U}$ denotes the restriction
of $Z_0$ to $f^{-1}(U)\times_U g^{-1}(U)$. By Proposition \ref{higherdirectimages},
$\rho_H(Z_0/S)$ is a morphism of $\OO_S$-modules.

Obviously, it is sufficient to prove that \eqref{equation-to-show-thmoverS},
for $j=0$ and $j=d$,
is an isomorphism for every open $U$. Thus we may suppose that $U=S$, 
$f^{-1}(U)=X$, $g^{-1}(U)=Y$, $Z_{0,U}=Z$, and we need to prove that 
\begin{align*}
\rho_H\circ cl([(Z_{0}])&:H^i(X,\OO_X)\xr{} H^i(Y\OO_Y), \\
\rho_H\circ cl([(Z_{0}])&:H^i(X,\Omega^d_X)\xr{} H^i(Y,\Omega^d_Y), 
\end{align*}
are isomorphisms for all $i$. In other words, we reduced to the case $S=\Spec(k)$.
 
Obviously, we may assume that $Z\subset X\times Y$. Let $Z'\subset Z, X'\subset X, Y'\subset Y$ be non-empty open subsets such that $\pr_1^{-1}(X')=Z'$, 
$\pr_2^{-1}(Y')=Z'$, and $\pr_1:Z'\xr{} X'$, $\pr_2:Z'\xr{} Y'$ are isomorphisms.

We obtain a correspondence $[Z]\in \Hom_{Cor_{\CH}}(X,Y)^0$, and we 
denote by $[Z^t]\in \Hom_{Cor_{\CH}}(Y,X)^0$ the correspondence defined by $Z$ considered
as subset of $Y\times X$.  

We claim that  
$$
[Z]\circ [Z^t] = \Delta_Y + E_1, \quad [Z^t]\circ [Z] = \Delta_X + E_2,
$$
with cycles $E_1$, resp.~$E_2$, supported in $(Y\backslash Y')\times (Y\backslash Y')$,
resp.~$(X\backslash X')\times (X\backslash X')$. 
Indeed, in view of Lemma \ref{remark:support}, $[Z^t]\circ [Z]$ is naturally supported in 
$$
\supp(Z,Z^t)=\{(x_1,x_2)\in X\times X \mid \exists y\in Y, (x_1,y)\in Z, (y,x_2)\in Z^t\}.
$$ 
By using Lemma \ref{lemma-composition-localization} for the open $X'\subset X$
we conclude that $[Z^t]\circ [Z]$ maps to $[\Delta_{X'}]$ via the localization 
map
$$
\CH(\supp(Z,Z^t))\xr{} \CH(\supp(Z,Z^t)\cap (X'\times X')).
$$
Thus 
$$
[Z^t]\circ [Z] = \Delta_X + E_2
$$
with $E_2$ supported in $\supp(Z,Z^t)\backslash (X'\times X')$. Finally, 
we observe that 
$$
\supp(Z,Z^t)\cap ((X'\times X)\cup (X\times X'))=\Delta_{X'}=\supp(Z,Z^t)\cap (X'\times X'),
$$ 
and thus $E_2$ has support in $(X\times X)\backslash ((X'\times X)\cup (X\times X'))=(X\backslash X')\times (X\backslash X')$.
The same argument works for $[Z]\circ [Z^t]$.
 
Now, Proposition \ref{propositionvanishing} implies that $\rho_H\circ cl([Z])$ induce isomorphisms 
$$
H^*(X,\OO_X) \xr{\simeq} H^*(Y,\OO_Y), \quad H^*(X,\Omega^d_X)\xr{\simeq} H^*(Y,\Omega^d_Y).
$$  

\end{proof}

\begin{corollary}\label{corollarystandartapplication}
Let $k$ be an arbitrary field and let $f:X\xr{} Y$ be a proper birational morphism
between smooth schemes $X,Y$. Then 
$$
Rf_*(\OO_X)=\OO_Y, \quad Rf_*(\omega_X)=\omega_Y. 
$$
\end{corollary}
\begin{proof}
By base change we may assume that $k$ is algebraically closed. The claim 
follows from Theorem \ref{thmoverS} for $S=Y, X=Z$.
\end{proof}

\subsubsection{} \label{criterionrational}
Consider a commutative diagram 
$$
\xymatrix
{
Y_a \arir[rr] \ar[drr]^{\text{gen.fin.}} 
&
&
Y \are[d]^{f}
\\
\ti{X} \ar[rr]^{\text{bir.}}_{\pi}
&
&
X.
}
$$
Here, all morphisms are proper and all schemes are integral, 
$\ti{X},Y$ are smooth of dimension $d_X, d_Y$, $\pi$ is birational, $f$ is surjective, 
$Y_a\xr{} X$ is generically finite and surjective, and  
$Y_a\inj Y$ is a closed immersion.
Let $\eta$ be the generic point of $X$ then $Y_a\times_X \eta$ is  finite 
over  $\Spec\, k(\eta)$ of degree $\deg(Y_a/X)$.

Choose a nonempty open set $U\subset X$ with $\pi:\pi^{-1}(U)\xr{\cong} U$ and 
such that $Y_a':=Y_a\cap f^{-1}(U)\subset f^{-1}(U) \xr{f} U$ is a finite morphism.
Set 
$$Z_a=\overline{Y_a'\times_U \pi^{-1}(U)}\subset Y_a\times_X \ti{X}$$ 
which gives a morphism $[Z_a]:Y\xr{} \tilde{X}$ in $Cor_{\CH}$. Furthermore, set 
$$
\Gamma=\overline{\pi^{-1}(U)\times_U f^{-1}(U)}\subset \ti{X}\times_X Y
$$ 
which defines an element $[\Gamma]: \ti{X}\xr{} Y$ in $Cor_{\CH}$. By using 
Lemma \ref{remark:support} and \ref{lemma-composition-localization}, we obtain
$$
[Z_a]\circ [\Gamma] = \deg(Y_a/X)\cdot  id_{\ti{X}} + E_1,
$$
where $E_1$ has support in $\pi^{-1}(X\setminus U)\times\pi^{-1}(X\setminus U)$; thus $\rho_H\circ cl (E_1)$ acts trivially on $H^*(\ti{X},\OO_{\ti{X}})\oplus H^*(\ti{X}, \Omega^{d_X}_{\ti{X}})$ by Proposition \ref{propositionvanishing}.
On the other hand, Lemma \ref{remark:support} and \ref{lemma-composition-localization} imply  
$$
[\Gamma] \circ [Z_a]  = \overline{Y_a'\times_U f^{-1}(U)} + E_2,
$$
where $E_2$ has support in $f^{-1}(X\backslash U)\times f^{-1}(X\backslash U)$; thus 
$\rho_H\circ cl(E_2)=0$ on $H^*(Y,\OO_{Y})\oplus H^*(Y,\Omega^{d_Y}_Y)$. 
Moreover, by using Lemma \ref{remark:support} and \ref{lemma-composition-localization} again, 
$$
[\Gamma] \circ [Z_a] \circ [\Gamma] \circ [Z_a]  = \deg(Y_a/X)\cdot 
\overline{Y_a'\times_U f^{-1}(U)} + E_3,
$$
with a cycle $E_3$ supported in $f^{-1}(X\backslash U)\times f^{-1}(X\backslash U)$, 
and therefore $\rho_H\circ cl(E_3)=0$ on $H^*(Y,\OO_{Y})\oplus H^*(Y,\Omega^{d_Y}_Y)$.

We obtain an endomorphism 
$$
P(Y_a):=\rho_H\circ cl([\Gamma] \circ [Z_a])\mid H^*(Y,\OO_{Y})\oplus H^*(Y,\Omega^{d_Y}_Y)
$$
of $H^*(Y,\OO_{Y})\oplus H^*(Y,\Omega^{d_Y}_Y)$ such that $P(Y_a)^2=\deg(Y_a/X)\cdot P(Y_a)$. 
Note that $P(Y_a)$ does not depend on $\tilde{X}$, because it is given by 
$$
P(Y_a)=\rho_H(cl([\overline{Y_a'\times_U f^{-1}U}])).
$$

\begin{proposition}\label{lemmaprojector}
If $\deg(Y_a/X)$ is invertible in $k$ then 
$$
\rho_H\circ cl(\Gamma):H^*(\ti{X},\OO_{\ti{X}})\oplus H^*(\ti{X},\Omega^{d_X}_{\ti{X}})\xr{} 
                                   {\rm image} \left(P(Y_a)\mid H^*(Y,\OO_{Y})\oplus H^*(Y,\Omega^{d_Y}_Y)\right)
$$
is a well-defined isomorphism.
\begin{proof}
Indeed $(\rho_H\circ cl(Z_a))\circ (\rho_H\circ cl(\Gamma)) \mid H^*(\ti{X},\OO_{\ti{X}})\oplus H^*(\ti{X},\Omega^{d_X}_{\ti{X}})$ 
is multiplication by $\deg(Y_a/X)$. It follows that $\rho_H\circ cl(\Gamma)$ is injective 
and the image is contained in the image of $P(Y_a)$. The opposite inclusion is obvious. 
\end{proof}
\end{proposition}

\begin{corollary} \label{maincorollary}
Let $Y,\ti{X}_{},X$ be as in \ref{criterionrational}. Let $a\in Y$ be a closed
point of the generic fiber $f^{-1}(\eta)$ with $\deg_{k(\eta)}(a)\in k^*$, we denote 
the corresponding closed subvariety by $Y_a$. For $i\geq 0$ the following are equivalent:
\begin{enumerate}
\item $R^i\pi_*(\OO_{\ti{X}})\oplus R^i\pi_*(\Omega^{d_X}_{\ti{X}})=0$ 
\item $P(Y_a\cap f^{-1}(X'))$ vanishes on 
\[ H^i(f^{-1}(X'),\OO_{f^{-1}(X')})\oplus H^i(f^{-1}(X'),\Omega^{d_X}_{f^{-1}(X')})\] 
for every affine open subset $X'\subset X$.
\end{enumerate}
\begin{proof}
In view of Proposition \ref{lemmaprojector} we get
$$
H^i(\pi^{-1}(X'),\OO_{\pi^{-1}(X')})
= {\rm image} \left(P(Y_a\cap f^{-1}(X'))\mid H^i(f^{-1}(X'),\OO_{f^{-1}(X')}) \right) 
$$
and
$$
H^i(\pi^{-1}(X'),\Omega^{d_X}_{\pi^{-1}(X')})
= {\rm image} \left(P(Y_a\cap f^{-1}(X'))\mid H^i(f^{-1}(X'),\Omega^{d_Y}_{f^{-1}(X')}) \right) 
$$
for every open subset  $X'\subset X$.
\end{proof}
\end{corollary}

\begin{thm}\label{thmquotients}
Let $k$ be an arbitrary field.
Consider  
$$
\xymatrix
{
&
Y\ar[d]^{f}
\\
\ti{X}
\ar[r]^{\pi}
&
X,
}
$$
where $Y,\ti{X}$ are smooth and connected, $X$ is integral and normal, $f$ is surjective and finite with $\deg(f)\in k^*$, and finally
$\pi$ is birational and proper.
Then X is Cohen-Macaulay and  
$$
R\pi_*\OO_{\ti{X}}=\OO_X, \quad R\pi_*\omega_{\ti{X}}= \omega_{X},
$$
where $\omega_X$ is the dualizing sheaf of $X$.
\begin{proof}
Choose an algebraic closure $\bar{k}$ of $k$.
We claim that $X$ is geometrically normal, i.e.~ 
\begin{equation}\label{equation-thmquotients-geom-normal}
X\times_k \bar{k}=\coprod_{i=1}^{r} X_i  \quad \text{(disjoint union)} 
\end{equation}
with $X_i$ integral and normal for all $i$. Indeed, since $X$ is normal we obtain 
$\OO_X\xr{=}\pi_*\OO_{\tilde{X}},$ and this isomorphism is stable under 
the base change to $\bar{k}$. Because $\tilde{X}$ is smooth, 
$\tilde{X}\times_k \bar{k}$ is a disjoint union of smooth schemes
$$
\tilde{X}\times_k \bar{k}=\coprod_{i=1}^{r} \tilde{X}_i.
$$
From $\OO_{X\times_k\bar{k}}\xr{=}\pi_*\OO_{\tilde{X}\times_k \bar{k}}$ we conclude that
$
\pi\times_k \bar{k}
$
has connected fibres; thus we obtain the equality \eqref{equation-thmquotients-geom-normal} 
with $X_i:=(\pi\times_k \bar{k})(\tilde{X_i})$. Of course, $\tilde{X}_i\xr{} X_i$
is birational, and $\OO_{X_i}\xr{=} (\pi\times_k \bar{k})_*\OO_{\tilde{X}_i}$ implies that $X_i$ is normal.

We denote by $Y_{\bar{k}}, X_{\bar{k}}, \ti{X}_{\bar{k}}$ the base change to $\bar{k}$, and 
by $\sigma_X:X_{\bar{k}}\xr{} X$ the obvious morphism.
Since 
$$
\sigma_X^*R^i\pi_* \OO_{\ti{X}} = R^i\pi_{\bar{k}*}\OO_{\ti{X}_{\bar{k}}}, \quad \sigma_X^*R^i\pi_* \omega_{\ti{X}} = R^i\pi_{\bar{k}*}\omega_{\ti{X}_{\bar{k}}}, 
$$
and $\sigma_X$ is faithfully flat it is sufficient to prove 
$$
R^i\pi_{\bar{k}*}\OO_{\ti{X}_{\bar{k}}}=0=R^i\pi_{\bar{k}*}\omega_{\ti{X}_{\bar{k}}} \quad \text{for all $i>0$.} 
$$

Now $\ti{X}\times_k \bar{k}=\coprod_{i=1}^{r} \ti{X}_i$ with $\ti{X}_i$ smooth and connected 
such that $\pi_{\bar{k}}\mid_{\ti{X}_i}:\ti{X}_i\xr{} X_i$ is birational. We define 
$Y_i:=Y_{\bar{k}}\cap f^{-1}_{\bar{k}}(X_i)$ and let $Y_i=\coprod_j Y_{i,j}$ be the decomposition
into connected (smooth) components. Since $\deg(Y_i/X_i)\in k^*$ is invertible there exists 
$j$ such that $\deg(Y_{i,j}/X_i)\in k^*$. Thus we are reduced to proving the claim for 
an algebraically closed field $k$.
 
Since $f$ is affine the statement $R\pi_*\OO_{\ti{X}}=\OO_X$  follows 
from Corollary \ref{maincorollary} and $X$ normal. Applying $D_X$ and shift by $[-d_X]$ (with $d_X=\dim X$)
we obtain $R\pi_*\Omega^{d_X}_{\ti{X}}= \pi_X^!k[-d_X]$ (with $\pi_X: X\to \Spec \, k$ the structure map).
Now again by Corollary \ref{maincorollary} we obtain $R^i\pi_*\Omega^{d_X}_{\ti{X}}=0$ for all $i\neq 0$.
Thus $\pi_X^!k[-d_X]=\omega_X$ is a sheaf and hence $X$ is Cohen-Macaulay. 
\end{proof}
\end{thm}

\section{Generalization to tame quotients}

The goal of this section is to generalize Theorem \ref{thmoverS} by replacing 
the assumption on the smoothness of $X,Y$ with the weaker assumption that 
$X,Y$ are tame quotients (see Definition \ref{tamequotsing}). We already proved in Theorem 
\ref{thmquotients} that the cohomology of the structure sheaf and 
the dualizing sheaf of a tame quotient behaves like
for smooth schemes. Therefore it is a natural question to extend Theorem \ref{thmoverS}
in order to include tame quotients.

\subsection{The action of finite correspondences}\label{coract} 
Let $X$ be a smooth scheme and $Z$ a closed integral subscheme of pure codimension $c$.
 Then $\sH^j_Z(\Omega^c_X)=0$  for $j<c$. Consequently, there is a natural morphism 
\eq{HZtoRGammaZ}{ \sH^c_Z(\Omega^c_X)[-c]\to R\ul{\Gamma}_Z\Omega^c_X, }
which induces an isomorphism $H^c_Z(X, \Omega^c_X)=H^0(X, \sH^c_Z(\Omega^c_X))$.

\begin{definition}\label{localpf}
Let $f: X\to Y$ be a morphism between smooth $k$-schemes of pure dimension $d_X$ and $d_Y$, respectively. Let $Z\subset X$ be a $c:=d_X-d_Y$ 
codimensional integral subscheme such that the restriction of $f$ to $Z$ is {\em finite}.
Then we define for $q\ge 0$ the local push-forward
\[f_{Z*}: f_*\sH^c_Z(\Omega^q_X)\to \Omega^{q-c}_Y\]
in the following way: Choose a compactification of $f$, i.e. a proper morphism $\bar{f}: \bar{X}\to Y$ and an open immersion $j :X\inj\bar{X} $ such that $f=\bar{f}\circ j$, and then define
$f_{Z*}$ as the composition in $D^+_{\rm qc}(Y)$ of the natural map $f_*\sH^c_Z(\Omega_X^q)\xr{\eqref{HZtoRGammaZ}} Rf_* (R\ul{\Gamma}_Z(\Omega^q_X)[c])$ with 
\ml{derivedlocalpf}{ Rf_* (R\ul{\Gamma}_Z(\Omega^q_X)[c])\xrightarrow[+\text{excision}]{\eqref{isoD(Omega)}}
                       R\bar{f}_*R\ul{\Gamma}_Z(D_{\bar{X}}(\Omega^{d_X-q}_{\bar{X}})[c-d_X])\\ 
           \xr{\text{forget support}} R\bar{f}_*D_{\bar{X}}(\Omega^{d_X-q}_{\bar{X}})[-d_Y]
          \xr{\bar{f}_*} D_Y(\Omega^{d_X-q}_Y)[-d_Y]\xr{\eqref{isoD(Omega)}} \Omega^{q-c}_Y, }
where $\bar{f}_*$ is the morphism from Definition \ref{defpush-forward}.

Applying $H^0(X, -)$ gives a morphism 
\[H^0(Y, f_*\sH^c_Z(\Omega_X^q))=H^0(X, \sH^c_Z(\Omega_X))=H^c_Z(X,\Omega^q_X)\to H^0(Y,\Omega^{q-c}_Y),\]
which by the very definition coincides with the cohomological degree zero part of the pushforward for $f:(X,Z)\to Y$, see Definition \ref{pfwithsupport}.
This also implies that $f_{Z*}$ is independent of the chosen compactification. 
\end{definition}

\begin{definition}\label{localcorr}
Let $S$ be a $k$-scheme and let $f: X\to S$ and $g:Y\to S$ be two integral $S$-schemes, which are smooth over $k$. Let $Z\subset X\times_S Y$ be a closed integral subscheme 
such that ${\pr_2}_{|Z}:Z\to Y$ is {\em finite and surjective}. In particular the codimension of $Z$ in $X\times Y$ equals $\dim X:=c$.
The projections from $X\times Y$ to $X$ and $Y$ are denoted by $\pr_1$ and $\pr_2$ respectively.
For all $q\ge 0$ we define a morphism  
\[\varphi_{Z}^q: f_*\Omega^q_X\to g_*\Omega_Y^q\]
as follows:
 Let $cl(Z)\in H^c_Z(X\times Y, \Omega^c_{X\times Y})=H^0(X\times Y, \sH^{c}_Z(\Omega^c_{X\times Y}))$ be the cycle class of $Z$. The cup product with $cl(Z)$ yields a morphism
\eq{localcup}{\cup cl(Z): \Omega^q_{X\times Y}\to \sH^c_Z(\Omega^{c+q}_{X\times Y})} 
and hence a morphism $f_*\pr_{1*}\Omega^q_{X\times Y}\to f_*\pr_{1*}\sH^c_Z(\Omega^{c+q}_{X\times Y})$.
We claim that it also induces a morphism of $\sO_S$-modules
\eq{twistedcup}{f_*\pr_{1*}\Omega^q_{X\times Y}\to g_*\pr_{2*}\sH^c_Z(\Omega^{q+c}_{X\times Y}).}
Indeed since $\sH^c_Z(\Omega^{q+c}_{X\times Y})$ has support in $Z\subset X\times_S Y$, the two abelian sheaves $g_*\pr_{2*}\sH^c_Z(\Omega^{q+c}_{X\times Y})$ and $f_*\pr_{1*}\sH^c_Z(\Omega^{q+c}_{X\times Y})$ are equal;
we denote this abelian sheaf on $S$ by $\sA$.
Now there are two  $\sO_S$-module structures on $\sA$: One is induced by $\sO_S\xr{g^*}g_*\sO_Y\xr{\pr_2^*}g_*\pr_{2*}\sO_{X\times Y}$ and the other is induced by
 $\sO_S\xr{f^*}f_*\sO_X\xr{\pr_1^*}f_*\pr_{1*}\sO_{X\times Y}$. The claim \eqref{twistedcup}  is now a consequence of the following equality in $\sA$:
\[\pr_2^*g^*(a)\cdot(\beta\cup cl(Z))=\pr_1^*f^*(a)\cdot (\beta\cup cl(Z)),\quad \text{for all } a\in \sO_S,\beta\in f_*\pr_{1*}\Omega^q_{X\times Y}, \]
which holds by \eqref{linear2}. We can therefore define the morphism $\varphi_Z^q$ as the composition
\begin{equation}\label{phi_Z}
f_*\Omega^q_X\xr{\pr_1^*} f_*\pr_{1*}\Omega^q_{X\times Y}\xr{\eqref{twistedcup}} g_*\pr_{2*}\sH^c_Z(\Omega^{c+q}_{X\times Y})\xr{\pr_{2, Z*}} g_*\Omega^q_Y.
\end{equation}
We write $\varphi_Z=\oplus_q \varphi_Z^q$.

Let $\alpha=\sum_i n_i [Z_i]$ be a formal sum of integral closed subschemes $Z_i$ of $X\times_S Y$, which are finite and surjective over $Y$, with coefficients $n_i$ in $\Z$.
Then we define
\eq{3.4.2}{\varphi_\alpha:= \sum_i n_i\varphi_{Z_i}: \bigoplus_q f_*\Omega^q_X \to \bigoplus_q g_*\Omega^q_Y.}
\end{definition}

\begin{lemma}\label{globallocalpf}
In the above situation, assume additionally that $f$ and $g$ are affine. Then for any cycle $\alpha=\sum_i n_i [Z_i]$, with $Z_i\subset X\times_S Y$ integral closed subschemes, which are finite and surjective over $Y$, and 
$n_i\in\Z$, we have the following  equality 
\[\oplus_i H^i(S,\varphi_\alpha)= \rho_H(cl(\bar{\alpha})): \bigoplus_{i,j} H^i(X,\Omega^j_X)\to \bigoplus_{i,j} H^i(Y,\Omega^j_Y), \]
where $\bar{\alpha}$ is the image of $\alpha$ in $\CH_{\dim Y}(X\times Y, P(\Phi_X,\Phi_Y))$ with $P(\Phi_X,\Phi_X)$ as in \eqref{definitionP}, $\rho_H$ is defined in \ref{VtoCor} and 
$cl$ is a shorthand notation for $Cor(cl)$ with $cl: CH\to H$ the morphism from \ref{propositionCHtoHodge}.
\end{lemma}

\begin{proof}
Let $\pi: S\to \Spec\,k $ be the structure map. We may assume $\alpha=[Z]$ with $Z\subset X\times_S Y$ an integral closed subscheme, which is finite and surjective over $Y$.
It is easy to see that $\rho_H(cl(\bar{\alpha}))$ is induced by taking the cohomology of the following composition
\begin{align}\label{derivedglobalcorr}
    R\pi_* f_*\Omega^q_X & \xr{\pr_1^*} R\pi_* R(f\pr_1)_*\Omega^q_{X\times Y}\\
                         & \xr{\cup cl(Z)} R\pi_* R(f\pr_1)_*\sH^c_Z(\Omega^{q+c}_{X\times Y}) \nonumber\\
    &\xr{Z\subset X\times_S Y} R\pi_* R(g\pr_2)_*\sH^c_Z(\Omega^{q+c}_{X\times Y}) \nonumber\\
&\xr{} R\pi_* (g\pr_2)_*\sH^c_Z(\Omega^{q+c}_{X\times Y}) \label{derivedglobalcorr-step1}\\
                         & \xr{\eqref{HZtoRGammaZ}} R\pi_*R(g\pr_2)_*R\ul{\Gamma}_Z(\Omega^{q+c}_{X\times Y})[c] \nonumber \\
                        &  \xr{\eqref{derivedlocalpf}} R\pi_*g_*\Omega_Y^q \label{derivedglobalcorr-step2}.
\end{align}
We used for the fourth arrow the isomorphism 
$$
(g\pr_2)_*\sH^c_Z(\Omega^{q+c}_{X\times Y}) \xr{\cong} R(g\pr_2)_*\sH^c_Z(\Omega^{q+c}_{X\times Y}),
$$ 
because $\sH^c_Z(\Omega^{q+c}_{X\times Y})$ is a quasicoherent $\OO_{X\times Y}$-module with support in $Z$, $Z\xr{} Y$ is finite, and $g:Y\xr{} S$ is affine. 
For the third arrow, notice that there is no
map  $ R(f \pr_1)_*\sH^c_Z \to  R(g \pr_2)_*\sH^c_Z $ in the derived category
of $\sO_S$-modules, but in the derived category of sheaves of $k$-vector  
spaces on $S$  these two complexes are isomorphic and that's all we need 
to define the third arrow.

We have to compare the morphism from \eqref{derivedglobalcorr} to \eqref{derivedglobalcorr-step2} with $R\pi_*\varphi_Z$, 
where $\varphi_Z$ is defined in \eqref{phi_Z}. 
Obviously, the morphism from \eqref{derivedglobalcorr} to \eqref{derivedglobalcorr-step1} is equal to $R\pi_*(\text{\eqref{twistedcup}}\circ \pr_1^*)$. 
The morphism from \eqref{derivedglobalcorr-step1} to \eqref{derivedglobalcorr-step2} equals $R\pi_*(\pr_{2, Z*})$, which proves the claim.

\end{proof}

\subsection{Tame quotients}
\subsubsection{}\label{2.1} Let $X$ be a $k$-scheme which is normal, Cohen-Macaulay (CM) and equidimensional of pure dimension $n$ and denote by $\pi: X\to \Spec\, k$ its structure map.

Then $H^i(\pi^!k)=0$ for all $i\neq -n$ (see \cite[Thm 3.5.1]{Co}). The {\em dualizing sheaf of $X$} is then by definition
            \[\omega_X:=H^{-n}(\pi^!k).\]
            We list some well-known properties:
\begin{enumerate}
     \item $\omega_X[n]$ is canonically isomorphic to $\pi^!k$ in $D^+_c(X)$.
     \item $\omega_X$ is a dualizing complex on $X$, i.e. $\omega_X$ is coherent, has finite injective dimension and the natural map 
           $\sO_X\to R\sHom(\omega_X,\omega_X)$ is an isomorphism. (Indeed by \cite[V, \S 10, 1., 2.]{Ha}, $\pi^!k$ is a dualizing complex.)
     \item $\omega_X$ is CM with respect to the codimension filtration on $X$, i.e. 
               \[{\rm depth_{\sO_{X,x}}} \omega_{X,x}=\dim \sO_{X,x},\quad \text{ for all } x\in X,\]
          (By \cite[V,  Prop. 7.3]{Ha} $\omega_X$ is Gorenstein, in particular CM, with respect to its associated filtration. Therefore we have to show that the associated codimension function 
             to $\omega_X$ (see \cite[V, \S 7]{Ha}) is the usual codimension function.
           By \cite[V,  Prop. 7.1]{Ha} it suffices to show that $\Ext^0_{\sO_{X,\eta}}(k(\eta),\omega_{X,\eta})\neq 0$ for all generic points $\eta\in X$.
But $X$ is normal and thus  $\omega_{X,\eta}\cong k(\eta)$. 
     \item In case $X$ is smooth $\omega_{X}$ is canonically isomorphic to $\Omega^n_X$, via the isomorphism 
               \[\omega_X\cong \pi^!k[-n]\xr{e_\pi, \, \simeq} \pi^\#k[-n]= \Omega^n_X,\]
                 where $e_\pi: \pi^!\cong \pi^\#$ is the isomorphism from \cite[(3.3.21)]{Co}.
     \item If $u: U\to X$ is \'etale, then $u^*\omega_X$ is canonically isomorphic to $\omega_U$ via the isomorphism
                \[u^*\omega_X= u^\#\omega_X\xr{e_u, \, \simeq} u^!\omega_X\cong u^!\pi^!k\xr{c_{u,\pi}^{-1}, \, \simeq} (\pi\circ u)^!k\cong \omega_U,\]
             where $c_{u,\pi}: (\pi\circ u)^!\cong u^!\circ \pi^!$ is the isomorphism from \cite[3.3.14]{Co}.
     \item Let $U$ be an open subscheme of $X$, which is smooth over $k$ and contains all $1$-codimensional points and denote by $j: U\inj X$ the corresponding open immersion.
           Then adjunction induces an isomorphism
           \[\omega_X\xr{\simeq} j_*j^*\omega_X\cong j_*\omega_U\cong j_*\Omega_U^n,\]
           where the two last isomorphism are induced by (4) and (5). 
(This follows from (3). 
Indeed,  let  $V\subset X$ be open then 
$\Gamma(V,\omega_X)\to \Gamma(V,j_*j^*\omega_X)=\Gamma(V\cap U,\omega_X)$ 
is the restriction. Since all points in the complement of $U$ have 
      codimension $\ge 2$, we obtain from (3) that ${\rm depth}(\omega_{X,x})\ge 2$ for all $x\in X\setminus U$. 
Therefore $\Gamma(V,\omega_X)\to \Gamma(V\cap U,\omega_X)$ is bijective by
       \cite[Exp. III, Cor. 3.5]{SGA2}.)
\end{enumerate}

\subsubsection{}\label{pb-pf-sm-to-CM-zero-level}
Let $X$ be smooth and $Y$ a normal CM scheme both of pure dimension $n$ and let $f: X\to Y$ be a finite and surjective morphism.
Then we have the usual pull-back on the structure sheaves $f^*: \sO_Y\to f_*\sO_X$ as well as a trace map $\tau^0_f: f_*\sO_X\to \sO_Y$,
which extends the the usual trace over the smooth locus of $Y$ (over which $f$ is flat). We define a pull-back and a trace between the dualizing sheaves as follows.
\begin{definition}\label{pb-pf-sm-to-CM-top-level}
Let $X$ be smooth and $Y$ a normal CM scheme both of pure dimension $n$ and let $f: X\to Y$ be a finite and surjective morphism.
\begin{enumerate}
\item We define a pullback morphism 
\[f^*: \omega_Y\to f_*\omega_X \]
as follows:
Choose $j: U\inj Y$ open and smooth over $k$ such that it contains all 1-codimensional points of $Y$, let $j': U'=X\times_Y U\to X$ and $f':U'\to U$ be the base changes of $j$ and $f$.
Then we define $f^*$ as the composition 
\[\omega_Y\simeq j_*\Omega^n_U\xr{{f'}^*} j_*f'_*\Omega^n_{U'}=f_*j'_*\omega_{U'}\cong f_*\omega_{X},\]
for the last isomorphism observe, that $U'$ contains all 1-codimensional points of $X$. It is straightforward to check, that this morphism is independent of the choice of $U$. (One only needs
the compatibility statements (VAR1) and (VAR3) of \cite[VII, Cor. 3.4, (a)]{Ha}.)
\item We define the trace
\[\tau_f^n: f_*\omega_X\to \omega_Y\]
as the composition in $D^+_c(Y)$
\[f_*\omega_X\cong f_*\pi_X^!k[-n]\xr{c_{f,\pi_Y}} f_*f^!\pi_Y^!k[-n]\xr{\Tr_f} \pi_Y^!k[-n]\cong\omega_Y,\]
where $\pi_X$ and $\pi_Y$ are the structure maps of $X$ and $Y$ and $\Tr_f$ is the trace morphism \cite[(3.3.2)]{Co}.
\end{enumerate}
We will also write $f^*:\sO_Y\oplus\omega_Y\to f_*(\sO_X\oplus \omega_X)$ for the sum of the usual pull-back with the pull-back defined in (1), and we will write
$\tau_f:=\tau_f^0\oplus\tau_f^n:f_*(\sO_X\oplus\omega_X)\to \sO_Y\oplus\omega_Y.$
\end{definition}
\begin{remark}
 By its very definition the $\tau_f$ constructed above equals, when restricted to the smooth locus of $Y$, the $\tau_f$ from Proposition \ref{finite-pf}.
\end{remark}

\begin{corollary}\label{pb-pf-equals-deg}
Let $X$, $Y$ and $f$ be as in  Definition \ref{pb-pf-sm-to-CM-top-level}. 
Suppose that $X$ is connected. Then the composition
\[\tau_f\circ f^*: \sO_Y\oplus\omega_Y\to \sO_Y\oplus\omega_Y\]
is equal to the multiplication with the degree of $f$.
\end{corollary}

\begin{proof}
We have to check, that the section $s= \tau_f\circ f^*-\deg f$ of 
\[H^0(Y, \sHom_Y(\sO_Y,\sO_Y))\oplus H^0(Y, \sHom_Y(\omega_Y,\omega_Y))\]
is zero. But $H^0(Y, \sHom_Y(\sO_Y,\sO_Y))=H^0(Y,\sO_Y)=H^0(Y, \sHom_Y(\omega_Y,\omega_Y))$ (for the last equality we need that $\omega_Y$ is a dualizing complex).
Therefore it is enough to check that $s$ is zero over an open and dense subset $U$ of $Y$. We may choose $U$ such that it is smooth and contains all $1$-codimensional points of $Y$.
Thus the statement follows from Proposition  \ref{finite-pf}, (3).
\end{proof}

\begin{definition}\label{tamequotsing}
Let $X$ be a $k$-scheme. 
We say that $X$ is a {\em  tame quotient} if $X$ is integral, normal and there exists a smooth and integral scheme $X'$ with a finite and surjective morphism $f:X'\to X$ whose degree is invertible in $k$.
\end{definition}

\begin{remark}
Assume $X$ is a tame quotient. Then $X$ is CM (see \cite[Prop. 5.7, (1)]{KM}).
\end{remark}

We may describe the cohomology of the structure sheaf and of the dualizing sheaf of a tame quotient
as a direct summand of the corresponding cohomology of a smooth scheme as follows.

\begin{proposition}\label{direct-summand}
Let $f: X\to Y$ be a finite and surjective morphism between integral schemes.
Assume $X$ is smooth and $Y$ is normal. Furthermore, we assume that $\deg f$ is invertible in $k$.
Set 
\[\alpha:= [X\times_Y X]\quad \text{in } \CH_{\dim X}(X\times X, P(\Phi_X,\Phi_X))=:\Hom_{Cor_{\CH}}(X,X)^0,\]
(cf. \eqref{definitionP} for the definition of $P(\Phi_X,\Phi_X)$). 
Then, for all $i$, the pull-back morphism
\[f^*: (H^i(Y,\sO_Y)\oplus H^i(Y,\omega_Y))\to (H^i(X,\sO_X)\oplus H^i(X,\omega_X))\]
induces an isomorphism
\[(H^i(Y,\sO_Y)\oplus H^i(Y,\omega_Y))\cong \rho_H(cl(\alpha))(H^i(X,\sO_X)\oplus H^i(X,\omega_X)).\]
(The functor $\rho_H$ is defined in \ref{VtoCor} and $cl$ is a shorthand 
notation for $Cor(cl)$ with $cl: CH\to H$ the morphism from Theorem \ref{propositionCHtoHodge}.)
\end{proposition}

\begin{proof}
Write $\alpha=[X\times_Y X]=\sum_T n_T [T]$, where the sum is over all irreducible components $T$ of $X\times_Y X$. Notice that all $T$'s
have dimension equal to $\dim X$ and project (via both projections) finitely and surjectively to $X$.
Therefore
\[\varphi_\alpha: f_*(\sO_X\oplus\omega_X)\to f_*(\sO_X\oplus\omega_X)\]
is defined, where $\varphi_\alpha$ is the morphism from Definition \ref{localcorr}.
By Lemma \ref{globallocalpf}  we have for all $i$:
\eq{globalsections}{H^i(Y,\varphi_\alpha)=\rho_H\circ cl(\alpha): H^i(X, \sO_X\oplus\omega_X)\to H^i(X,\sO_X\oplus\omega_X).}
We claim 
\eq{3.6.1}{\varphi_\alpha=f^*\circ \tau_f: f_*(\sO_X\oplus\omega_X)\to f_*(\sO_X\oplus\omega_X).}
Let $U\subset Y$ be a non-empty smooth open subscheme which contains all 1- codimensional points of $Y$. Then $f^{-1}(U)$ is smooth and contains all $1$-codimensional points of $X$.
Hence for any open $V\subset Y$ we have that the restriction map 
\[H^0(f^{-1}(V), \sO_X\oplus\omega_X)\to H^0(f^{-1}(V)\cap f^{-1}(U), \sO_X\oplus\omega_X)\]
is an isomorphism (see \cite[Exp. III, Cor. 3.5]{SGA2}). Since both maps in \eqref{3.6.1} are compatible with restriction to open subsets of $Y$, we may therefore assume that $Y$ is smooth.
In particular $f$ is flat and thus $\alpha$ equals $[\Gamma_f^t]\circ [\Gamma_f]$, where $\Gamma_f$ is the graph of $f$ and $\Gamma_f^t$ its transposed. Now
the identity \eqref{3.6.1} follows from \eqref{globalsections} (in the case $i=0$), Proposition \ref{properties-pfwithsupports}, (3) and \ref{pb-pf-via-corr}.
Thus applying again \eqref{globalsections} we obtain
\mlnl{\rho_H(cl(\alpha))(H^i(X,\sO_X)\oplus H^i(X,\omega_X))\\
= {\rm Image}(f^*\circ \tau_f: H^i(X, \sO_X\oplus\omega_X))\to H^i(X, \sO_X\oplus\omega_X)).}
Since $\tau_f\circ f^* : (H^i(Y,\sO_Y)\oplus H^i(Y,\omega_Y))\to (H^i(Y,\sO_Y)\oplus H^i(Y,\omega_Y))$ is multiplication with the degree of $f$ (by Corollary \eqref{pb-pf-equals-deg})
the statement of the proposition follows.
\end{proof}

\subsection{Main theorem for tame quotients}

\begin{thm}\label{thm-tame-quotients}
Let $S$ be a scheme over a perfect  field $k$. Let  $\pi_X:X\to S$ and $\pi_Y:Y\to S$ be two integral $S$-schemes, 
which are tame quotients (cf.~Definition \ref{tamequotsing}). Furthermore, we assume that $X$ and $Y$ are properly birational equivalent.
Then any $Z$ as in Definition \ref{propbir} induces isomorphisms of $\sO_S$-modules 
\[R^i \pi_{X*}\sO_X\cong R^i\pi_{Y*}\sO_Y,\quad R^i\pi_{X*}\omega_X\cong R^i\pi_{Y*}\omega_Y,\quad\text{for all } i\ge 0.\]
These isomorphisms depend only on the $\sO_{S,\eta}$-isomorphism $k(X)\cong k(Y)$ induced by $Z$, where $\eta=\pi_X(\text{generic point of }X)=\pi_Y(\text{generic point of }Y)$.
\end{thm}

\begin{proof} We first prove 

{\em  Claim 1: There are isomorphisms as in the statement in the case $S=\Spec\, k$.}

Choose integral and smooth schemes $X'$ and $Y'$ with finite and surjective morphisms $f: X'\to X$ and $g: Y'\to Y$ whose degree is invertible in $k$.
Choose $Z$ as in Definition \ref{propbir}. We may assume that $Z\subset X\times Y$ is a closed integral subscheme. We define $Z_{X'}$, $Z_{Y'}$ and $Z'$ by the cartesian diagram
\[\xymatrix{       &        &     Z'\ar[dl]\ar[dr]  &      &       \\
                   & Z_{X'}\ar[dl]_{\sim}\ar[dr] &    & Z_{Y'}\ar[dr]^{\sim}\ar[dl] &   \\
            X'\ar[dr]_{f}  &     &  Z\ar[dl]_\sim\ar[dr]^\sim  &    &  Y'\ar[dl]^g\\
                          & X     &                     & Y. &      }\]
Here the arrows with an $\sim$ are proper and birational morphisms between integral schemes and all other morphisms are finite and surjective.
Notice that we may identify $Z'$ with a closed subscheme of $X'\times Y'$ whose irreducible components are proper and surjective over both $X'$ and $Y'$, 
and all irreducible components have the same dimension equal 
to $d:=\dim X=\dim X'=\dim Y=\dim Y'$ (since $f$ and $g$ are finite and universally equidimensional).
Therefore $Z'$ and its transpose define cycles $[Z']\in \CH^d(X'\times Y', P(\Phi_{X'}, \Phi_{Y'}))$ and $[Z']^t\in \CH^d(Y'\times X', P(\Phi_{Y'}, \Phi_{X'}))$. 
Now choose non-empty smooth open subschemes $X_o$, $Y_o$ of $X,Y,$ such that the morphisms $Z\to X, Z\to Y$ induce 
isomorphisms $Z_o\xr{\simeq} X_o, Z_o\xr{\simeq} Y_o$ with 
$$
Z_o:=X_o\times_X Z= Z\times_Y Y_o.
$$
Set $X_o'=f^{-1}(X_o)$ and $Y_o'= g^{-1}(Y_o)$ and denote by $f_o$ and $g_o$ the restrictions
of $f$ and $g$ to $X_o'$ and $Y_o'$, respectively. We define $Z_{X_o'}$, $Z_{Y_o'}$  and $Z'_o$ by the cartesian diagram
\[\xymatrix{       &        &     Z_o'\ar[dl]\ar[dr]  &      &       \\
                   & Z_{X_o'}\ar[dl]_{\simeq}\ar[dr] &    & Z_{Y_o'}\ar[dr]^{\simeq}\ar[dl] &   \\
            X_o'\ar[dr]_{f_o}  &     &  Z_o\ar[dl]_\simeq\ar[dr]^\simeq  &    &  Y_o'\ar[dl]^{g_o}\\
                          & X_o     &                     & Y_o. &      }\]
Here the arrows with an $\simeq$ are isomorphisms, all other arrows are finite and surjective.
We set  $X'_c=X'\setminus X'_o$ and  $Y'_c=Y'\setminus Y'_o$; these are closed subsets  of codimension $\ge 1$. 
Now we define 
\begin{align*}
\alpha&:=[X'\times_X X']\in \CH^d(X'\times X', P(\Phi_{X'}, \Phi_{X'})),\\
\beta&:=[Y'\times_Y Y']\in \CH^d(Y'\times Y', P(\Phi_{Y'}, \Phi_{Y'})).
\end{align*}
We claim 
\eq{3.9.2}{\deg g\cdot ([Z'] \circ \alpha)- \deg f\cdot (\beta\circ [Z'])\in {\rm image}(\CH_*(X'_c\times Y'_c)), }
\eq{3.9.3}{\deg f\cdot ([Z']^t\circ \beta)-\deg g\cdot (\alpha\circ [Z']^t)\in {\rm image}(\CH_*(Y'_c\times X'_c)),}
\eq{3.9.4}{([Z']^t\circ [Z']\circ \alpha) -\deg f\deg g\cdot\alpha \in {\rm image}(\CH_*(X'_c\times X'_c)),}
\eq{3.9.5}{([Z']\circ [Z']^t\circ \beta)- \deg f\deg g \cdot \beta \in {\rm image}(\CH_*(Y'_c\times Y'_c)).}

By symmetry, it suffices to prove \eqref{3.9.2} and \eqref{3.9.4}.
Let us prove \eqref{3.9.2}. By using Lemma \ref{remark:support} we can 
consider 
$$
\alpha\in \CH(X'\times_X X'), \quad \beta\in \CH(Y'\times_Y Y'), \quad [Z']\in \CH(Z'),
$$ 
and see that
$[Z'] \circ \alpha$ and $\beta\circ [Z']$ are naturally supported in $\CH(Z')$.

Since $Z'\cap ((X_o'\times Y')\cup (X'\times Y'_o))=Z'\cap (X_o'\times Y'_o)$, 
Lemma \ref{lemma-composition-localization} and the localization sequence for Chow groups 
implies the claim provided that the 
equality 
\begin{equation}\label{equation-1-thm-3.9}
\deg(g) \cdot [Z']_{\mid X'\times Y'_o}\circ \alpha_{\mid X'_o\times X'} = \deg(f)\cdot \beta_{\mid Y'\times Y'_o}\circ [Z']_{\mid X'_o\times Y'}
\end{equation}
holds in $\CH(Z'\cap (X_o'\times Y'_o))$. Here we have 
\begin{align*}
 \alpha_{\mid X'_o\times X'}&\in \CH(X'_o\times_{X} X')= \CH(X_o'\times_{X_o} X_o')\\
\beta_{\mid Y'\times Y'_o} &\in  \CH(Y'\times_{Y} Y'_o) = \CH(Y'_o\times_{Y_o} Y'_o) \\
[Z']_{\mid X_o'\times Y'}&\in \CH(Z'\cap (X_o'\times Y'))=\CH(Z'\cap (X_o'\times Y'_o))\\
[Z']_{\mid X'\times Y_o'}&\in \CH(Z'\cap (X'\times Y_o'))= \CH(Z'\cap (X_o'\times Y'_o)).
\end{align*}
Obviously,
\begin{align}\label{align-thm-3.9-on-o}
\alpha_{\mid X'_o\times X'}&=[X_o'\times_{X_o} X_o']=[\Gamma_{f_o}^{t}]\circ 
[\Gamma_{f_o}] \\
\beta_{\mid Y'\times Y'_o}&=[Y_o'\times_{Y_o} Y_o']=[\Gamma_{g_o}^{t}]\circ [\Gamma_{g_o}] \nonumber\\
[Z']_{\mid X_o'\times Y'}&=[Z'_o]=[\Gamma_{g_o}^t]\circ [Z_o]\circ [\Gamma_{f_o}] \nonumber\\
[Z']_{\mid X'\times Y_o'}&=[Z'_o]. \nonumber
\end{align}
Thus \eqref{equation-1-thm-3.9} follows from 
\begin{align}\label{align-thm-3.9-Gamma-Gammat}
[\Gamma_{f_o}]\circ [\Gamma_{f_o}^{t}]&= \deg(f)[\Delta_{X_o}]\\
[\Gamma_{g_o}]\circ [\Gamma_{g_o}^{t}]&= \deg(g)[\Delta_{Y_o}].\nonumber
\end{align}
This finishes the proof of \eqref{3.9.2}. The  proof of \eqref{3.9.4} is similar. 
The cycles $[Z']^t\circ [Z']\circ \alpha$ and  $\alpha$ are supported in 
$$
B=\{(x_1',x_2')\in X'\times X' \mid \exists y\in Y:\; (f(x_1'),y)\in Z, (f(x_2'),y)\in Z \}.
$$ 
We see that $B\cap((X_o'\times X')\cup (X'\times X'_o))=B\cap(X_o'\times X'_o)$,
and by using Lemma \ref{lemma-composition-localization} it is 
sufficient to prove 
$$
[Z'_o]^t\circ [Z_o']\circ [X'_o\times_{X_o} X'_o] =\deg f\deg g\cdot[X'_o\times_{X_o} X'_o].
$$
In view of \eqref{align-thm-3.9-on-o} this follows immediately from \eqref{align-thm-3.9-Gamma-Gammat}.

Since $\deg f$ and $\deg g$ are invertible in $k$, 
it follows from Proposition \ref{propositionvanishing} and \eqref{3.9.2} that $\frac{1}{\deg f}\rho_H\circ cl([Z'])$ induces a morphism
\[(\rho_H\circ cl(\alpha))H^*(X', \sO_{X'}\oplus\omega_{X'})\to (\rho_H\circ cl(\beta))H^*(Y', \sO_{Y'}\oplus\omega_{Y'})\]
and by \eqref{3.9.3} $\frac{1}{\deg g}\rho_H\circ cl([Z']^t)$ induces a morphism  
\[(\rho_H\circ cl(\beta))H^*(Y', \sO_{Y'}\oplus\omega_{Y'})\to (\rho_H \circ cl(\alpha))H^*(X', \sO_{X'}\oplus\omega_{X'}).\]
By \eqref{3.9.4} and \eqref{3.9.5} these two morphism are inverse to each other. 
Thus Proposition \ref{direct-summand} yields isomorphisms
\[H^i(X,\sO_X)\cong H^i(Y,\sO_Y),\quad H^i(X,\omega_X)\cong H^i(Y,\omega_Y),\quad \text{for all }i\ge 0.\]
This proves Claim 1.

{\em Claim 2: The isomorphisms constructed in Claim 1 depends only on the isomorphism $k(X)\cong k(Y)$ induced by a $Z$. }

We use the shorthand notation $H^i(X)= H^i(X,\sO_X\oplus\omega_X)$. 
Choose $Z$ as in Definition \ref{propbir}. Denote  by $Z_0$ the image of $Z$ in $X\times Y$.
Choose $f_1: X_1\to X$ and $g_1: Y_1\to Y$ finite and surjective, with $X_1$, $Y_1$ smooth and integral and $\deg f_1$, $\deg g_1\in k^*$.
Define 
\[\ul{\alpha}_1:= \rho_H\circ cl([X_1\times_X X_1]), \quad \ul{\beta}_1:=\rho_H\circ cl([Y_1\times_Y Y_1]),\]
\[\ul{\gamma}_1(Z):=\frac{1}{\deg f_1}\rho_H\circ cl([X_1\times_X Z_0\times_Y Y_1]).\] 
Then we saw in the proof of Claim 1 above, that we obtain isomorphisms
\[H^i(X)\xr{f_1^*,\,\simeq} \ul{\alpha}_1 H^i(X_1)\xr{\ul{\gamma}_1(Z),\,\simeq} \ul{\beta}_1H^i(Y_1)\xl{g_1^*,\,\simeq} H^i(Y) .\]
Now choose two different $Z$'s as in Definition \ref{propbir}, say $Z_1$, $Z_2$, which induce the same isomorphism $k(X)\cong k(Y)$.
Then we can find smooth open subschemes $X_o$, $Y_o$, $Z_{1,o}$, $Z_{2,o}$ of $X$, $Y$, $Z_1$, $Z_2$, such that for $i=1,2$ we have
\[Z_{i,o}= Z_i\times_X X_o= Z_i\times_Y Y_o, \]
the projections $Z_{i,o}\xr{\simeq} X_o$, $Z_{i,o}\xr{\simeq} Y_o$ are isomorphisms and the induced isomorphisms $h_i: X_o\xr{\simeq} Z_{o,i}\xr{\simeq} Y_o$, $i=1,2$, are {\em equal}.
Proposition \ref{propositionvanishing} implies 
$$
\ul{\gamma}_1(Z_1)=\ul{\gamma}_1(Z_2)
$$ on $H^i(X_1)$. Therefore $\ul{\gamma}_1(Z)$ depends only on the 
isomorphism $k(X)\cong k(Y)$, which $Z$ induces. From now on we fix such an isomorphism and simply write $\ul{\gamma}_1$.

Now choose $f_2: X_2\to X$ and $g_2: Y_2\to Y$ finite and surjective, with $X_2$, $Y_2$ smooth and integral and $\deg f_2$, $\deg g_2\in k^*$.
Define $\ul{\alpha}_2$, $\ul{\beta}_2$ and $\ul{\gamma}_2$ as above (in the above formulas replace 1 by 2) and set
\[\ul{\alpha}_{12}:=\frac{1}{\deg f_1}\rho_H\circ cl([X_1\times_X X_2]),\quad \ul{\alpha}_{21}:=\frac{1}{ \deg f_2}\rho_H\circ cl([X_2\times_X X_1])\]
\[\ul{\beta}_{12}:=\frac{1}{  \deg g_1}\rho_H\circ cl([Y_1\times_Y Y_2]), \quad \ul{\beta}_{21}:=\frac{1}{\deg g_2}\rho_H\circ cl([Y_2\times_Y Y_1]).\]
Then one checks as in the proof of Claim 1 that $\ul{\alpha}_{12}: H^i(X_1)\to H^i(X_2)$ induces an isomorphism $\ul{\alpha}_1H^i(X_1)\xr{\simeq} \ul{\alpha}_2H^i(X_2)$
with inverse $\ul{\alpha}_{21}$ and $\ul{\beta}_{12}: H^i(Y_1)\to H^i(Y_2)$ induces an isomorphism $\ul{\beta}_1H^i(Y_1)\xr{\simeq}\ul{\beta}_2 H^i(Y_2)$ with inverse
$\ul{\beta}_{21}$. Further one checks that $\ul{\beta}_{12}\circ \ul{\gamma_1}\circ\ul{\alpha}_1= \ul{\gamma}_2\circ\ul{\alpha}_{12}\circ\ul{\alpha}_1$.
Thus we obtain the following commutative diagram
\eq{isoindep}{\xymatrix{           &     \ul{\alpha}_1H^i(X_1)\ar[r]^{\ul{\gamma}_1,\,\simeq}\ar[dd]^{\ul{\alpha}_{12}}_{\simeq} & \ul{\beta}_1H^i(Y_1)\ar[dd]^{\ul{\beta}_{12}}_{\simeq}   &           \\
             H^i(X)\ar[ur]^{f_1^*,\,\simeq}\ar[dr]_{f_2^*,\,\simeq} &                &       & H^i(Y).\ar[ul]_{g_1^*,\,\simeq}\ar[dl]^{g_2^*,\, \simeq} \\
                 &     \ul{\alpha}_2H^i(X_2)\ar[r]^{\ul{\gamma}_2,\,\simeq}                            & \ul{\beta}_2H^i(Y_2)   & 
                }}
Therefore the isomorphisms of Claim 1 do not depend on the choice of $f_1,g_1$. This proves Claim 2 and also the theorem in the case $S=\Spec\, k$.

{\em Finally we consider the case of a general basis $S$.}
Choose $Z$ as in Definition \ref{propbir} and choose integral, smooth schemes $X',Y'$ with 
finite, surjective morphisms $f: X'\to X, g: Y'\to Y$ whose degree is invertible in $k$.
For $U\subset S$ open denote by $X_U$, $f_U$, etc.~the pull-backs over $U$.
By Proposition \ref{direct-summand} the pull-back $f_U^*$ realizes $H^i(X_U, \sO_{X_U}\oplus \omega_{X_U})$ as a direct summand of $H^i(X'_U, \sO_{X'_U}\oplus \omega_{X'_U})$.
This is clearly compatible with restrictions along opens $V\subset U\subset S$ and thus the pull-back $f^*$ realizes
$R^i\pi_{X*}(\sO_X\oplus \omega_{X})$ as a direct summand of the $\sO_S$-module $R^i\pi_{X*}f_*(\sO_{X'}\oplus \omega_{X'})$.
In the same way $g^*$ realizes $R^i\pi_{Y*}(\sO_Y\oplus \omega_{Y})$ as a direct summand of the $\sO_S$-module $R^i\pi_{Y*}g_*(\sO_{Y'}\oplus \omega_{Y'})$.
Further by the case $S=\Spec\,k$ considered above, the map 
\[\frac{1}{\deg f}\rho_H\circ cl([Z'_U]):H^*(X'_U, \sO_{X'_U}\oplus\omega_{X'_U})\to H^*(Y'_U, \sO_{Y'_U}\oplus\omega_{Y'_U})\]
induces an isomorphism between $H^*(X_U, \sO_{X_U}\oplus\omega_{X_U})$ and $H^*(Y_U, \sO_{Y_U}\oplus\omega_{Y_U})$.
Write $[Z']=\sum_T n_T [T] $, where the sum is over the irreducible  components of $Z'$. Then the collection
$\{\frac{1}{\deg f}\rho_H\circ cl([Z'_U])\,|\, U\subset S\}$ induces a morphism of $\sO_S$-modules (by Proposition \ref{higherdirectimages})
\[\rho_H(Z'/S):=\sum_T n_T\, \rho_H(T/S): R^i\pi_{X*}f_*(\sO_{X'}\oplus \omega_{X'})\to R^i\pi_{Y*}g_*(\sO_{Y'}\oplus \omega_{Y'}),\]
which by the above induces an isomorphism
\eq{relisomainthm}{R^i\pi_{X*}(\sO_{X}\oplus \omega_{X})\xr{\simeq}R^i\pi_{Y*}(\sO_{Y}\oplus \omega_{Y}).}
Claim 2 implies that \eqref{relisomainthm} depends only on the 
$\sO_{S,\eta}$-isomorphism $k(X)\cong k(Y)$ 
induced by $Z$.
\end{proof}

\begin{remark} 
Theorem \ref{thm-tame-quotients} implies Theorem \ref{thmoverS} and Theorem 
\ref{thmquotients}.
\end{remark}

\begin{corollary} 
Let $\pi: X\to Y$ be a birational and proper morphism between integral schemes over a perfect field $k$.
Assume $X$ and $Y$ are tame quotients.
Then $\pi^*$ induces isomorphisms
\[R\pi_*\sO_X\cong\sO_Y,\quad R\pi_*\omega_X\cong \omega_Y.\]
\end{corollary}
\begin{proof}
In Theorem \ref{thm-tame-quotients} take $S=Y$, $\pi_X=\pi$ and $\pi_Y=\id_Y$.
\end{proof}

\subsection{Open questions}

Questions in ${\rm char}(k)=p$.    
 

\subsubsection{}
Do the statements in Corollary \ref{thmoverS} and Theorem \ref{thmquotients} hold when $k$ is 
not perfect and  smooth is replaced by regular? 

\subsubsection{}
Let $f:Y\xr{} X$ be a surjective projective morphism with connected fibres between smooth varieties $Y,X$.
Is $R^{\dim(Y)-\dim(X)}f_*\omega_Y=\omega_X$?
In ${\rm char}(k)=0$ this holds by \cite[Proposition~7.6]{Ko1}.

\subsubsection{}
Let $f:Y\xr{} X$ be a surjective projective morphism with connected fibres between smooth varieties $Y,X$. 
Is $R^{e}f_*\omega_Y=0$ for $e>\dim(Y)-\dim(X)$?
In ${\rm char}(k)=0$ this holds by \cite[Theorem~2.1(ii)]{Ko1}.

\appendix
\section{ }

{\em All schemes in this Appendix are assumed to be finite dimensional and noetherian.}

\subsection{Local Cohomology}\label{A1}
Let $Y=\Spec \,B$ be an affine scheme and $X\subset Y$ a closed subscheme of pure codimension $c$, defined by the ideal $I\subset B$.
We assume that there exists a $B$-regular sequence $t=t_1,\ldots,t_c\in I$ with $\sqrt{(t)}=\sqrt{I}$, where $(t)$ denotes the ideal $(t_1,\ldots, t_c)\subset B$.
We denote by $K^\bullet(t)$ the  Koszul complex of the sequence $t$,
i.e. $K^{-q}(t)=K_q(t)= \bigwedge^q B^c$, $q=0,\ldots, c$, and if $\{e_1,\ldots, e_c\}$ is the standard basis of $B^c$ and
$e_{i_1,\ldots, i_q}:= e_{i_1}\wedge\ldots\wedge e_{i_q}$, then the differential is given by
\[d^{-q}_{K^\bullet}(e_{i_1,\ldots, i_q})=d^{K_\bullet}_{q}(e_{i_1,\ldots, i_q})=
                        \sum_{j=1}^{q}(-1)^{j+1} t_{i_j}e_{i_1,\ldots,\widehat{i_j},\ldots i_q}.\]
For any $B$-module $M$ we define the complex 
\[K^\bullet(t,M):= \Hom_B(K^{-\bullet}(t), M),\]
and denote its $n$-th cohomology by $H^n(t, M)$. The map 
\[\Hom_B(\bigwedge^c B^c, M)\to M/(t)M,\quad \varphi\mapsto \varphi(e_{1,\ldots,c})\] induces a canonical isomorphism
$H^c(t, M)\simeq M/(t)M$. 

If $t$ and $t'$  are two sequences as above with $(t')\subset (t)$, then there exists a $c\times c$-matrix $T$ with coefficients in $B$ such that
$t'= T t$ and $T$ induces a morphism of complexes $K^\bullet(t')\to K^\bullet(t)$, which is the unique (up to homotopy) morphism lifting
the natural map $B/(t')\to B/(t)$. Furthermore we observe that, for any pair of sequences $t$, $t'$
as above there exists an $N\ge 0$ such that $(t^N)\subset (t')$, where $t^N$ denotes the sequence $t_1^N,\ldots, t_c^N$. Thus 
the sequences $t$ form a directed set and $H^c(t, M)\to H^c(t', M)$, $(t')\subset (t)$, becomes a direct system. 
It follows from \cite[Exp. II, Prop. 5]{SGA2}, that we have an isomorphism
\[\varinjlim_t M/(t)M=\varinjlim_t H^c(t, M) \cong H^c_X(Y, \widetilde{M}),\]
where the limit is over all $B$-regular sequences $t=t_1,\ldots, t_c$ in $B$ with $V((t))=X$ and $\widetilde{M}$ is the sheaf associated to $M$.
We denote by
\[\genfrac{[}{]}{0pt}{}{m}{t}\]
the image of $m\in M$ under the composition
\[M\to M/(t)M\to H^c(t, M)\to H_X^c(Y, \widetilde{M}).\]
It is a consequence of the above explanations that we have the following properties:
\begin{enumerate}
 \item Let $t$ and $t'$ be two sequences as above with $(t')\subset (t)$. Let $T$ be a $c\times c$-matrix with $t'=Tt$, then
      \[\genfrac{[}{]}{0pt}{}{\det(T)\,m}{t'}=\genfrac{[}{]}{0pt}{}{m}{t}.\]
\item\[\genfrac{[}{]}{0pt}{}{m+m'}{t}=\genfrac{[}{]}{0pt}{}{m}{t}+ \genfrac{[}{]}{0pt}{}{m'}{t}, \quad 
                 \genfrac{[}{]}{0pt}{}{t_i m}{t}=0 \quad\text{all }i.  \]
\item If $M$ is any  $B$-module of finite rank, then
        \[H_X^c(Y,\OO_Y)\otimes_B M\xr{\simeq} H^c_X(Y, \widetilde{M}), \quad 
                                              \genfrac{[}{]}{0pt}{}{b}{t}\otimes m\mapsto \genfrac{[}{]}{0pt}{}{bm}{t}\]
is an isomorphism.
\end{enumerate}
 
\begin{remark}\label{extsign}
Since for a $B$-regular sequence  $t$ as above $K^\bullet(t)\to B/(t)$ is a free resolution, we have an isomorphism for all $n$,
\[\Ext^n(B/(t), M)\simeq H^n(\Hom^\bullet_B(K^\bullet(t), M)).\] 
 Notice that we also have an isomorphism 
\[\Hom^\bullet_B(K^\bullet(t), M)\simeq K^\bullet(t,M),\]
which is given by multiplication with $(-1)^{n(n+1)/2}$ in degree $n$. 
We obtain an isomorphism
\[\psi_{t,M}: \Ext^c(B/(t),M)\xr{\simeq} H^c(t,M)= M/(t)M, \]
which has the sign $(-1)^{c(c+1)/2}$ in it. 
In particular under the composition
\[\Ext^c(B/(t), M)\xr{\psi_{t,M}} M/(t)M \to H^c_X(Y,\widetilde{M})\]
the class of a map $\varphi\in \Hom_B(\bigwedge^c B^c,M)$ is sent to 
\[(-1)^{c(c+1)/2}\genfrac{[}{]}{0pt}{}{\varphi(e_{1,\ldots,c})}{t}.\]
\end{remark}

\begin{lemma}\label{Cousindiff}
Let $Y=\Spec\, B$ be as above,  $\mc{M}$ a quasi-coherent sheaf on $Y$,  $c\ge 0$ and $t_1,\ldots, t_{c+1}$ a $B$-regular sequence. 
We set  $X':=V(t_1,\ldots, t_{c+1})\subset X:= V(t_1,\ldots, t_c)$. 
Let $\partial: H^c_{X\setminus X'}(Y\setminus X', \mc{M})\to H^{c+1}_{X'}(Y,\mc{M})$ be the boundary map of the localisation long exact sequence.
Then
\[\partial\genfrac{[}{]}{0pt}{}{m/t_{c+1}}{t_1,\ldots,t_c}= \genfrac{[}{]}{0pt}{}{m}{t_1,\ldots,t_c, t_{c+1}}.\]
 
\end{lemma}
\begin{proof}
 Let $M$ be the $B$-module of global sections of $\mc{M}$. By \cite[Exp. II, Cor. 4]{SGA2} and $\check{C}$ech computations, we may identify 
\[H^c_{X\setminus X'}(Y\setminus X', \mc{M})=\frac{M_{t_1\cdots t_c t_{c+1}}}{\sum_{i=1}^c M_{t_1\cdots\widehat{t_i}\cdots t_c t_{c+1}}},\quad
 H^{c+1}_{X'}(Y, \mc{M})=\frac{M_{t_1\cdots t_{c+1}}}{\sum_{i=1}^{c+1}M_{t_1\cdots\widehat{t_i}\cdots t_{c+1}}}\]
and $\partial$ is the natural map from left to right. Under this identifications the map
$M_{t_{c+1}}/(t_1,\ldots, t_c)= H^c(K^\bullet(t,M))\to H^c_{X\setminus X'}(Y, \mc{M})$ sends the class of $m/t_{c+1}$ ($m\in M$)
to the class of $\frac{m/t_{c+1}}{t_1\cdots t_c}$ and similar for $M/(t_1,\ldots, t_{c+1})\to H^{c+1}_{X'}(Y, \mc{M})$.
This proves the claim.
\end{proof}

\subsection{The trace for a regular closed embedding}
In this section we give an explicit description of the trace morphism for a regular closed embedding. This is well-known and appears 
in various incarnations in the literature, see e.g. \cite{Lipman}, \cite[Section 4]{HS}. But in all the articles we are aware of, more elementary version of duality
theory are used (e.g. no derived categories appear). Since the compatibility of this theories with the one we are working with, namely the one developed in \cite{Ha} and \cite{Co}, is not evident to the authors, and also
to be sure about the signs, we recall the description of the trace in this situation.

Let $i: X\inj Y$ be a closed immersion of pure codimension $c$ between two Gorenstein schemes and 
assume that the ideal sheaf $\mc{I}$ of $X$ is generated  by a sequence  $t=(t_1, \ldots, t_c)$ of global sections of $\sO_Y$.
Then the image of $t$ in any local ring of $Y$ is automatically a regular sequence. We denote by $K^\bullet(t)$ the sheafified
Koszul complex of $t$ and set
\[\omega_{X/Y}:= \bigwedge^c \sHom_{\OO_X}(\mc{I}/\mc{I}^2, \OO_X).\]
The fundamental local isomorphism (see e.g. \cite[2.5]{Co}) gives an isomorphism in $D^b_c(Y)$
\eq{FLI}{\eta_i: i_*\omega_{X/Y}[-c]\xl{\simeq}\sHom^\bullet_Y(K^\bullet(t), \OO_Y)\cong R\sHom_Y(i_*\OO_X, \OO_Y)\cong i_*i^!\OO_Y. }
The first map is induced by
\ml{FLI-explicit}{\sHom(\bigwedge^c\OO_Y^c,\OO_Y)=\sHom^c(K^\bullet(t), \OO_Y)\to i_*\omega_{X/Y}, \\ \varphi\mapsto (-1)^{\frac{c(c+1)}{2}}\varphi(e_{1,\dots,c}) t_1^\vee\wedge\ldots\wedge t_c^\vee. }
(The reason for the sign is Remark \ref{extsign}.)
Composing the morphism $\eta_i$ with the trace $\Tr_i: i_*i^!\OO_Y\to \OO_Y$ (see e.g. \cite[3.4]{Co}) we obtain a morphism in $D^b_{\rm c}(Y)$
\eq{lci-trace}{i_*\omega_{X/Y}[-c]\xr{\eta_i} i_*i^!\OO_Y\xr{\Tr_i} \OO_Y}
which factors in $D_{\rm qc}^b(Y)$ as
\eq{lci-trace2}{i_*\omega_{X/Y}[-c]\xr{\eta_i} i_*i^!\OO_Y\xr{\Tr_i} R\ul{\Gamma}_X\OO_Y.}

\begin{lemma}\label{trace-for-lci-appendix}
In the above situation there is a natural isomorphism 
\eq{coh-suuport-derived}{R\ul{\Gamma}_X\OO_Y\cong \sH^c_X(\sO_Y)[-c],\quad \text{in } D^b_{\rm qc}(Y)}
 and $\sH^c(\eqref{lci-trace2})$ is given by
\eq{lci-trace-with-support}{i_*\omega_{X/Y}\lra \sH^c_X(\OO_Y),\quad  at_1^\vee\wedge\ldots\wedge t_c^\vee\mapsto (-1)^{\frac{c(c+1)}{2}}\genfrac{[}{]}{0pt}{}{\tilde{a}}{t_1,\ldots,t_c},}
where $\tilde{a}\in\sO_Y$ is any lift of $a\in\sO_X$.
\end{lemma}
\begin{proof}
The first statement is equivalent to $\sH^i_X(\sO_Y)=0$, for $i\neq c$, and hence we may assume that $Y$ is affine. We have the vanishing for $i<c$ since $Y$ is CM  (by \cite[Exp. III, Prop. 3.3]{SGA2})
and for $i>c$ since the ideal of $X$ in $Y$ is generated by $c$ elements, which by a $\check{C}$ech-argument implies that $H^i(Y\setminus X, \sO_Y)=0$ for $i>c$.

We denote by $E^\bullet=E^\bullet(\OO_Y)$ the Cousin complex of $\OO_Y$ (see e.g. \cite[IV,\S 2]{Ha}). In particular $E^\bullet$ is an injective resolution of $\OO_Y$ (since $Y$ is Gorenstein) and 
if $Y^{(c)}$ denotes the set of points of codimension $c$ in $Y$, then
\[E^q= \bigoplus_{y\in Y^{(c)}}i_{y*}H^q_y(Y,\OO_Y),\]
where $i_y: y\to Y$ is the inclusion and $H^q_y(Y,\OO_Y)=\varinjlim_{y\in U} H^q_{\overline{y}\cap U}(Y,\OO_Y)$, the limit being over all open subsets $U\subset Y$ which contain $y$. Further we denote $K^\bullet:=K^\bullet(t)$. 

The trace $\Tr_i: i_*i^!\sO_Y\to \sO_Y$ is now induced by the ``evaluation at 1'' morphism $\sHom^\bullet(i_*\sO_X, E^\bullet)\to E^\bullet$.
Furthermore the augmentation morphisms $K^\bullet\to i_*\sO_X$ and $\sO_Y\to E^\bullet$ induce quasi-isomorphisms
\eq{compute-Rhom}{\sHom^\bullet(K^\bullet, \OO_Y)\xr{\sim} \sHom^\bullet(K^\bullet, E^\bullet)\xl{\sim} \sHom^\bullet(i_*\OO_X, E^\bullet).}
To prove the second statement, we may assume $a=1\in \sO_X$.
We define $\alpha\in \sHom^c(K^\bullet, \OO_Y)=\sHom(\wedge^c\OO_Y^c,\OO_Y)$ by
\[\alpha(e_{1,\ldots,c})=1\]
and $\beta\in \sHom^c(i_*\OO_X, E^\bullet)= \sHom(i_*\OO_X, E^c)$ by
\[\beta(1)= (\beta_y)\in E^c=\bigoplus_{y\in Y^{(c)}}i_{y*}H^c_y(Y,\OO_Y)\]
           with
           \[\beta_y=\begin{cases} \genfrac{[}{]}{0pt}{}{1}{t_1,\ldots, t_c}, & \text{if } y \text{ is a generic point of }X,\\
                                 0, & \text{else.} \end{cases}\]
Then 
\[Tr_i(\bar{\beta})= \genfrac{[}{]}{0pt}{}{1}{t_1,\ldots,t_c} \in \sH^c_X(\sO_Y), \]
where $\bar{\beta}$ is the residue class of $\beta$ in $\sH^c(\sHom^\bullet(i_*\sO_X, E^\bullet))$ and 
\[\eta_i(\alpha)= (-1)^{\frac{c(c+1)}{2}}t_1^\vee\wedge\ldots\wedge t_c^\vee \in i_*\omega_{X/Y}.\]
Thus the second statement of the lemma follows if we can show that the images of $\alpha$ and $\beta$ in $\sHom^c(K^\bullet, E^\bullet)$ differ
by an element in $d_{\sHom^\bullet}^{c-1}(\sHom^{c-1}(K^\bullet, E^\bullet))$.

For $j=0,\ldots, c-1$ we define 
\[\gamma^{c-1-j}=(\gamma_y^{c-1-j})\in E^{c-1-j}=\bigoplus_{y\in Y^{(c-1-j)}} i_{y*} H^{c-1-j}_y(Y,\OO_Y)\] 
by
\[\gamma_{y}^{c-1-j}:=
\begin{cases} \genfrac{[}{]}{0pt}{}{1/t_{c-j}}{t_1,\ldots, t_{c-j-1}}, & \text{ if } y\in Y^{(c-j-1)}\cap V(t_1,\ldots,t_{c-1-j}),\\
              0, & \text{ else.}\end{cases}\]
In particular $\gamma^0=\frac{1}{t_1}\in \bigoplus_i H^0_{\eta_i}(Y,\OO_Y)=\oplus_i k(\eta_i)$, with $\eta_i$ the generic points of $Y$.
Notice that (by Lemma \ref{Cousindiff})
\begin{eqnarray}
\label{gamma-rel-1} d_E \gamma^{c-1-j} & = & t_{c-j+1} \gamma^{c-j},\quad   \forall\, j\ge 1.\\
\label{gamma-rel-2} t_i\gamma^{c-1-j} & = & 0,  \quad \forall \,j\ge 0,\, i\in\{1,\ldots,c-1-j\}. 
\end{eqnarray}
Further define
\[\psi=(\psi_0,\ldots,\psi_{c-1})\in\sHom^{c-1}(K^\bullet, E^\bullet)=\bigoplus_{j=0}^{c-1}\sHom(K_j, E^{c-1-j})\]
by
\[\psi_j(e_{i_1,\ldots,i_j})=\begin{cases}(-1)^{j(c+j)}\gamma^{c-1-j}, & \text{ if } (i_1,\ldots, i_j)=(c+1-j,\ldots,c),\\
                                             0 , &\text{ else.} \end{cases}\]
By definition and \eqref{gamma-rel-2} we have
\eq{psi-vanishing}{t_{i_q}\psi_{j-1}(e_{i_1,\ldots,\widehat{i_q},\ldots, i_j})\neq 0 \text{ if } q\neq 1
                                                                                      \text{ or } (i_1,\ldots,i_j)\neq (c+1-j,\ldots, c).  }
Now we calculate the boundary of $\psi$, 
\[d^{c-1}_{\sHom}\psi=(d_E^{c-1}\circ\psi_0,\ldots, d_E^{c-1-j}\circ\psi_j+(-1)^c\psi_{j-1}\circ d^K_j,\ldots, (-1)^c\psi_{c-1}\circ d^K_c).\]

{\em 1st Case: $j=0$. } 
\[d_E^{c-1}\circ\psi_0(1)=d_E^{c-1}\gamma^{c-1}= \beta(1),\quad \text{by } \eqref{gamma-rel-1}.\]

{\em 2nd Case: $1\le j\le c-2$. }
By \eqref{psi-vanishing} and the definition of $\psi$ we have
\[(d_E^{c-1-j}\circ\psi_j+(-1)^c\psi_{j-1}\circ d^K_j)(e_{i_1,\ldots,i_j})=0\] 
if $(i_1,\ldots,i_j)\neq(c+1-j,\ldots, c)$. Else 
\mlnl{(d_E^{c-1-j}\circ\psi_j+(-1)^c\psi_{j-1}\circ d^K_j)(e_{c+1-j,\ldots,c})\\
      = (-1)^{j(c+j)}d_E^{c-1-j}\gamma^{c-1-j}+ (-1)^c t_{c+1-j}\psi_{j-1}(e_{c+2-j,\ldots,c})\\
    = (-1)^{j(c+j)}(d_E^{c-1-j}\gamma^{c-1-j}-t_{c-j+1}\gamma^{c-j})=0, \quad \text{by \eqref{gamma-rel-1}.}}

{\em 3rd Case: $j=c-1$. } By \eqref{psi-vanishing} we have
\[(-1)^c\psi_{c-1}(d_c^K(e_{1,\ldots,c}))= (-1)^ct_1\psi_{c-1}(e_{2,\ldots, c})=-1=-\alpha(e_{1,\ldots,c}).\]

All in all we obtain
\[d^{c-1}_\sHom( \psi)=(\beta, 0,\ldots,0,-\alpha)\]
and this proves the lemma. 
\end{proof}

\begin{lemma}\label{1.7.11}
Let $S$ be a Gorenstein scheme and $i: X\inj Y$ a closed immersion between smooth, separated and equidimensional $S$-schemes with structure maps $\pi_X: X\to S$, $\pi_Y: Y\to S$ 
and denote by $d_{X/S}$ and $d_{Y/S}$ their relative dimensions. We set $\omega_{X/S}:=\Omega^{d_{X/S}}_{X/S}$, $\omega_{Y/S}:=\Omega^{d_{Y/S}}_{Y/S}$ and $c=d_{Y/S}-d_{X/S}$.
Assume that the ideal sheaf of $X$ in $Y$ is generated by a sequence $t=t_1, \ldots, t_c$ of global sections of $\sO_Y$.
Define a morphism $\imath_X$ by
\[\imath_X: i_*\omega_{X/S}\to \sH^c_X(\omega_{Y/S}), \quad \alpha\mapsto (-1)^c\genfrac{[}{]}{0pt}{}{dt\tilde{\alpha}}{t},\]
with $\tilde{\alpha}\in \Omega^{d_{X/S}}_{Y/S}$ any lift of $\alpha$ and $dt= dt_1\wedge\ldots\wedge dt_c$. 
Then the following diagram in $D^b_{\rm qc}(\sO_Y)$ is commutative
\eq{1.7.11.0}{\xymatrix{ i_* \pi_X^! \sO_S\ar[r]^-{c_{i, \pi_Y}}& i_* i^! \pi_Y^! \sO_S \ar[r]^-{\Tr_i} &  \pi_Y^! \sO_S\\ 
            i_*\omega_{X/S}[d_{X/S}]\ar[u]^\simeq\ar[r]^-{\imath_X} & \sH^c_X(\omega_{Y/S}) [d_{X/S}] \ar[r]^{\simeq} & R\ul{\Gamma}_X(\omega_{Y/S}) [d_{Y/S}]\ar[u]         } }
where the vertical map on the left is the well-known canonical isomorphism (see \cite[(3.3.21)]{Co}), the vertical map on the right
is the composition of the forget supports map $R\ul{\Gamma}_X(\omega_{Y/S}) [d_{Y/S}]\to \omega_{Y/S}[d_{Y/S}]$ with the canonical isomorphism $\omega_{Y/S}[d_{Y/S}]\cong \pi_Y^!\sO_S$
and $c_{i,\pi_Y}: \pi_X^!\cong i^!\pi_Y^!$ is the canonical isomorphism \cite[(3.3.14)]{Co}.

\end{lemma}

\begin{proof}
Let $\sI\subset \sO_Y$ be the ideal sheaf of $X$. As above we write 
\[\omega_{X/Y}:= \bigwedge^c \sHom_{\sO_X}(\sI/\sI^2, \sO_X).\]
Further let $\tau_Y: \pi_Y^!\sO_S\cong\omega_{Y/S}[d_{Y/S}]$ and $\tau_X:\pi_X^!\sO_S\cong \omega_{X/S}[d_{X/S}]$ be the canonical isomorphisms and
$\eta_i: \omega_{X/Y}[-c]\cong i^!\sO_Y$ the fundamental local isomorphism \eqref{FLI}.
Consider the following diagram in $D^b_{\rm c}(\sO_Y)$
\eq{1.7.11.1}{\xymatrix{i_*\pi_X^!\sO_S\ar@{=}[d]\ar[r]^-{c_{i,\pi_Y},\,\simeq} & i_*i^!\pi_Y^!\sO_S\ar[r]^{\Tr_i}\ar[d]^\simeq & \pi_Y^!\sO_S\ar@{=}[d]\\
            i_*\pi_X^!\sO_S\ar[d]^{\tau_X}_\simeq\ar[r]^-\simeq &  i_*i^!\sO_Y\otimes^L \pi_Y^!\sO_S\ar[d]^{\eta_i^{-1}\otimes \tau_Y}_\simeq\ar[r]^-{\Tr_i\otimes\id } & \pi_Y^!\sO_S\ar[d]^{\tau_Y}_\simeq\\
            i_*\omega_{X/S}[d_{X/S}]\ar[r]^-\simeq        & i_*\omega_{X/Y}[-c]\otimes \omega_{Y/S}[d_{Y/S}]\ar[r]^-{\Tr'_i\otimes \id} & \omega_{Y/S}[d_{Y/S}]. } }
Here some explanations: The middle horizontal arrow on the left is defined such that the upper left square commutes.  We have a canonical identification
$i_*i^!(-)=R\sHom_{\sO_Y}(i_*\sO_X, (-))$ and since $\pi_Y^!\sO_S$ is isomorphic to a shifted locally free $\sO_Y$-module we have 
$R\sHom_{\sO_Y}(i_*\sO_X, \pi_Y^!\sO_S)=R\sHom_{\sO_Y}(i_*\sO_X, \sO_Y) \otimes^L\pi_Y^!\sO_S$; this defines the upper vertical arrow in the middle. Furthermore $\Tr_i: i_*i^!(-)\to (-)$ may be identified
with $R\sHom_{\sO_Y}(i_*\sO_X, -)\to (-)$ given by the evaluation at 1. This shows, that in the above diagram the upper square on the right commutes.
The map $\Tr'_i: i_*\omega_{X/Y}[-c]\to \sO_Y$ on the right bottom is the composition \eqref{lci-trace} and thus the lower square on the right commutes by definition.
The horizontal isomorphism on the lower left is given by (see \cite[p.29, 30, (c), (2.2.6)]{Co}, note the ordering)
\ml{1.7.11.2}{i_*\omega_{X/S}[d_{X/S}]\to i_*\omega_{X/Y}[-c]\otimes \omega_{Y/S}[d_{Y/S}], \\ \alpha\mapsto (t_1^\vee\wedge\ldots\wedge t_c^\vee)\otimes dt_c\wedge\ldots\wedge dt_1\wedge \tilde\alpha,}
with $\tilde\alpha\in \Omega^{d_{X/S}}_{Y/S}$ any lift of $\alpha$. That the square on the lower left commutes, follows from
\cite[Thm 3.3.1, (3.3.27)]{Co} and \cite[Lem 3.5.3]{Co}. (Notice that by \cite[p.5, pp.160-164]{CoErr} the statement of \cite[Lem 3.5.3]{Co} should be ``(3.5.8) is equal to (3.5.7)'' instead of
 ``(3.5.8) is equal to $(-1)^{n(N-n)}$ times (3.5.7)''.) Thus the whole diagram commutes. The upper line equals the upper line in \eqref{1.7.11.0} 
and the lower line factors as the composition 
\ml{1.7.11.3}{i_*\omega_{X/S}[d_{X/S}]\to i_*\omega_{X/Y}[-c]\otimes \omega_{Y/S}[d_{Y/S}]\\ \xr{{\rm tr}_i\otimes\id} \sH^c_X(\sO_Y)[-c]\otimes\omega_{Y/S}[d_{Y/S}]\simeq \sH^c_X(\omega_{Y/S})[d_{X/S}]}
with the natural map 
\[\sH^c_X(\omega_{Y/S})[d_{X/S}] \simeq R\ul{\Gamma}_X(\omega_{Y/S})[d_{Y/S}]\to \omega_{Y/S}[d_{Y/S}].\]
Here ${\rm tr}_i$ denotes the composition of $R\ul{\Gamma}_X(\Tr_i')$ with the isomorphism $R\ul{\Gamma}_X(\sO_Y)\cong \sH^c_X(\sO_Y)[-c]$.
Thus the lemma is proved once we know that \eqref{1.7.11.3} equals $\imath_X$. But by Lemma \ref{trace-for-lci-appendix} the map ${\rm tr_i}$  is given by
\[i_*\omega_{X/Y}\to \sH^c_X(\sO_Y), \quad t_1^\vee\wedge\ldots\wedge t_c^\vee \mapsto (-1)^{\frac{c(c+1)}{2}}\genfrac{[}{]}{0pt}{}{1}{t}.\]
Together with \eqref{1.7.11.2} we obtain, that \eqref{1.7.11.3} is given by
\mlnl{\alpha\mapsto (t_1^\vee\wedge\ldots\wedge t_c^\vee)\otimes dt_c\wedge\ldots\wedge dt_1\wedge \tilde\alpha \\ 
            \mapsto (-1)^{\frac{c(c+1)}{2}}\genfrac{[}{]}{0pt}{}{1}{t}\otimes dt_c\wedge\ldots\wedge dt_1\wedge \tilde\alpha =(-1)^c\genfrac{[}{]}{0pt}{}{dt\tilde\alpha}{t},}
which by definition equals $\imath_X$. This proves the lemma.
\end{proof}

\subsection{The trace for a finite and surjective morphism}

\subsubsection{}\label{A3.1}
Let $S$ be a Gorenstein scheme and $f:X\to Y$ a finite and surjective morphism of smooth, separated and equidimensional  $S$-schemes, which both have relative dimension $n$. 
We denote by $\pi_X: X\to S$ and $\pi_Y: Y\to S$ the respective structure maps. 
Then we define  the trace map
\eq{A3.1.1}{\tau_f^n: f_*\omega_{X/S}\to \omega_{Y/S}}
to be the composition 
\[f_*\omega_{X/S}\cong Rf_*\pi_X^!\sO_S[-n]\xr{\simeq\, c_{f,\pi_Y}} Rf_*f^!\pi_Y^!\sO_S[-n]\xr{\Tr_f} \pi_Y^!\sO_S[-n]\cong \omega_{Y/S}. \]

In the lemma below we give a well-known explicit description of this trace map, for which we could not find an appropriate reference. 
There are well studied ad hoc definitions of this trace map not using the machinery of duality theory (see e.g. \cite[\S 16]{Ku}), but it is a priori not clear that these construction coincide with the one above.

\begin{lemma}\label{A3.2}
Let $f: X\to Y$ be as above and assume it factors as
\[\xymatrix{X\ar@{^(->}[r]^i\ar[d]_f & P\ar[dl]^\pi\\
             Y,}\]
where $\pi$ is smooth and separated of pure relative dimension $d$ and $i$ is a closed immersion whose ideal sheaf $\sI\subset \sO_P$ is generated by global sections $t_1,\ldots, t_d\in \Gamma(P,\sO_P)$.
Then for any local section $\alpha\in f_*\omega_{X/S}$ we have the following formula for $\tau^n_f(\alpha)$:
 Let $\tilde{\alpha}\in \Omega^n_{P/S}$ be any lift of $\alpha$ and write in $i^*\Omega_{P/S}^{n+d}= i^*\omega_{P/S}= i^*(\omega_{P/Y})\otimes f^*\omega_{Y/S}$
\[i^*(dt_d\wedge\ldots\wedge dt_1\wedge \tilde{\alpha})= \sum_j i^*\gamma_j\otimes f^*\beta_j,\quad \gamma_j\in \omega_{P/Y}, \, \beta_j\in\omega_{Y/S}.\]
Then  
\[\tau^n_f(\alpha) = (-1)^{d(d-1)/2} \sum_j  \Res_{P/Y}\genfrac{[}{]}{0pt}{}{\gamma_j}{t_1,\ldots, t_d} \beta_j \in \omega_{Y/S}, \]
where $\Res_{P/Y}\genfrac{[}{]}{0pt}{}{\gamma_j}{t_1,\ldots, t_d}\in \sO_Y$ is the residue symbol defined in \cite[(A.1.4)]{Co}.
\end{lemma}

\begin{proof}
The proof is a collection of compatibility statements from \cite{Co}. First we collect some notations.
\begin{enumerate}
 \item $\zeta'_{i,\pi_P}: \omega_{X/S}\to \omega_{X/P}\otimes i^*\omega_{P/S}$ is defined in \cite[p.29/30, (c)]{Co} and sends $\alpha$ to 
            $(t_1^\vee\wedge\ldots\wedge t_d^\vee)\otimes i^*(dt_d\wedge\ldots\wedge dt_1\wedge\tilde{\alpha})$, where we identify $\omega_{X/P}=\bigwedge^d (\sI/\sI^2)^\vee$.
\item $\eta_i: \sExt^d_P(i_*\sO_X,-) \xr{\simeq} \omega_{X/P}\otimes i^*(-)$ is the fundamental local isomorphism \cite[(2.5.1)]{Co}.
\item For a smooth and separated morphism of pure relative dimension $n$ between two schemes $g: V\to W$, $e_g: g^!\xr{\simeq} g^\#=\omega_{V/W}[n]\otimes^{\bm L}(-)$ denotes the natural transformation \cite[(3.3.21)]{Co}.
\item In case $g$ as above factors as $g=h\circ i$ with $i:V\to Z$ finite and $h: Z\to W$ smooth, then $\psi_{i,h}: g^\#\xr{\simeq}i^\flat h^\#$ is the isomorphism defined in \cite[(2.7.5)]{Co}, where
        $i^\flat(-)= i^{-1}R\sHom_Z(i_*\sO_V,-)\otimes_{i^{-1}i_*\sO_V} \sO_V$ is defined in \cite[(2.2.8)]{Co}.
\item $d_f: f^!\xr{\simeq} f^\flat$ is the isomorphism defined in \cite[(3.3.19)]{Co}.
\item $\Tr_f: f_*f^!\to \id$ is the trace morphism defined in \cite[3.4]{Co}, and ${\rm Trf}_f: f_*f^\flat\to \id$ is the finite trace morphism defined in \cite[(2.2.9)]{Co} and which is induced by evaluation at 1.
\end{enumerate}
Consider the following diagram:
\eq{A3.2.1}{\mbox{\tiny $ \xymatrix{ f_*\omega_{X/S}\ar@{=}[r]\ar[d]_{\simeq}\ar@{}[dr]|*[o]+[F]{1} & f_*\omega_{X/S}\ar[d]^\simeq\ar[rr]^{\zeta'_{i,\pi_P}}\ar@{}[drr]|*[o]+[F]{2} &
                                                                                                                               &       f_*(\omega_{X/P}\otimes i^*\omega_{P/S})\ar[d]^{\eta_i^{-1}}  \\
             f_*\pi_X^!\sO_S[-n]\ar[d]_{c_{f,\pi_Y}}\ar[r]^{e_{\pi_X}}\ar@{}[dr]|*[o]+[F]{3} & f_*\pi^\#_X\sO_S[-n]\ar[r]^{\psi_{i,\pi_P}}\ar[d]^{\psi_{f,\pi_Y}}\ar@{}[dr]|*[o]+[F]{4} & 
                                                               f_*(i^\flat\pi^\#_P\sO_S[-n])\ar[d]^\simeq\ar[r]^\simeq\ar@{}[ddr]|*[o]+[F]{5} &     f_*\sExt^d_P(i_*\sO_X,\omega_{P/S})\ar[dd]^\simeq  \\
             f_*f^!\pi_Y^!\sO_S[-n]\ar[d]_\simeq\ar[r]^{d_f}_{e_{\pi_Y}}\ar@{}[dr]|*[o]+[F]{6} & f_*f^\flat\pi^\#_Y\sO_S[-n]\ar[d]^\simeq\ar[r]^{\psi_{i,\pi}}\ar@{}[dr]|*[o]+[F]{7} & 
                                           f_*(i^\flat\pi^\#\pi_Y^\#\sO_S[-n])\ar[d]^\simeq &  \\
             f_*f^!\omega_{Y/S}\ar[r]^{d_f}\ar[d]_{\Tr_f}\ar@{}[dr]|(.3)*[o]+[F]{8} & f_*f^\flat\omega_{Y/S}\ar[d]^\simeq\ar[dl]|{{\rm Trf}_f}\ar[r]^{\psi_{i,\pi}}\ar@{}[dr]|*[o]+[F]{10} & 
                                         f_* i^\flat\pi^\#\omega_{Y/S}\ar[d]^\simeq\ar[r]^-\simeq\ar@{}[dr]|*[o]+[F]{11} &  f_*\sExt^d_P(i_*\sO_X, \omega_{P/Y}\otimes\pi^*\omega_{Y/S})\ar[d]^\simeq  \\
             \omega_{Y/S}   &  f_*f^\flat(\sO_Y)\otimes\omega_{Y/S}\ar[l]^-{{\rm Trf}\otimes\id}\ar[r]^{\psi_{i,\pi}}\ar@{}[ul]|(.3)*[o]+[F]{9} & 
                                                                  f_*(i^\flat\pi^\#\sO_Y)\otimes \omega_{Y/S}\ar[r]^-\simeq &  f_*\sExt^d_P(i_*\sO_X,\omega_{P/Y})\otimes \omega_{Y/S}. 
                                   } $}  }
Let us describe the different squares and triangles in this diagram:
\begin{enumerate}
 \item[1] The vertical isomorphism on the right in square 1 is immediate from the definition of $\pi_Y^\#$, the left vertical isomorphism is defined such that the square commutes.
\item[2] See \cite[Thm. 3.5.1, Cor. 3.5.2]{Co} for the isomorphism in the lower right of square 2. The square commutes by \cite[Lem. 3.5.3]{Co}.
         (Notice by \cite[comment to pp. 160-164]{CoErr} the last statement of \cite[Lem. 3.5.3]{Co} should be {\em ..., then (3.5.8) is equal to (3.5.7)}.)
\item[3] Square 3 commutes by \cite[(3.3.27)]{Co}.
\item[4] The vertical isomorphism on the right of square 4 is induced by the natural isomorphism $(\pi\pi_Y)^\#\cong \pi_Y^\#\pi^\#$. For the commutativity of the square see
         the discussion in \cite[p. 83/84]{Co} (our case is point three).
\item[5] The isomorphism on the right of square 5 is induced by the natural isomorphism $\omega_{P/S}\cong \omega_{P/Y}\otimes \pi^*\omega_{Y/S}$. The square commutes 
        by the functoriality of the horizontal isomorphisms, which are just induced by taking the $0$-th cohomology (the other cohomology groups being zero).
\item[6] The vertical isomorphism on the right of square 6 is induced by $\omega_{Y/S}\cong \pi_Y^!\sO_S[-n]\cong \pi_Y^\#\sO_S[-n]$. Thus the square commutes by the functoriality of $d_f$.
\item[7]The vertical isomorphism on the right of square 7 is defined as above. Thus the square commutes by the functoriality of $\psi_{i,\pi}$.
\item[8] Triangle 8 commutes by \cite[Lem 3.4.3, (TRA2)]{Co}.
\item[9] By \cite[III, proof of Prop. 6.5]{Ha} we may identify $f_*f^\flat\omega_{Y/S}$ with the sheaf  $\sHom_Y(f_*\sO_X,\omega_{Y/S})$ (since $\omega_{Y/S}$ is locally free) and ${\rm Trf}_f$ is given
         by evaluation at 1. The vertical isomorphism on the right in triangle 9 is defined by the isomorphism $\sHom_Y(f_*\sO_X,\omega_{Y/S})\cong \sHom_Y(f_*\sO_X,\sO_Y)\otimes \omega_{Y/S}$.
         The triangle thus obviously commutes.
\item[10] By \cite[(2.8.3) and the paragraph after this, p.100/101]{Co} we have a commutative square
            \[\xymatrix{f^\flat\omega_{Y/S}\ar[r]^{\psi_{i,\pi}}\ar[d]^\simeq & i^\flat\pi^\#\omega_{Y/S}\ar[d]^\simeq\\
                        f^\flat\sO_Y\otimes f^*\omega_{Y/S}\ar[r]^{\psi_{i,\pi}} & i^\flat\pi^\#\sO_Y\otimes f^*\omega_{Y/S}. }\]
          Applying $f_*$ to this diagram and using projection formula defines the commutative square 10.
\item[11] The horizontal maps in square 11 are induced by taking the $0$-th cohomology, the vertical maps are the natural isomorphisms ($\omega_{Y/S}$ is locally free \& projection formula).
           The commutativity of the diagram is clear.
\item[12] The isomorphism in the upper right of triangle 12 is induced by the isomorphism $\omega_{P/S}\cong\omega_{P/Y}\otimes \pi^*\omega_{Y/S}$ and the projection formula.
          The triangle commutes by \cite[Thm. 2.5.2, 1.]{Co}.
\end{enumerate}
Thus diagram \eqref{A3.2.1} is commutative. The composition of the vertical maps along the left outer edge of the diagram equals $\tau_f^n$ by definition.
The composition of the vertical maps along the right outer edge of the diagram is by \cite[Thm. 2.5.2, 1.]{Co} equal to the composition
\mlnl{f_*(\omega_{X/P}\otimes i^*\omega_{P/S})\xr{\simeq} f_*(\omega_{X/P}\otimes i^*\omega_{P/Y})\otimes \omega_{Y/S}\\
                                      \xr{\eta^{-1}_i\otimes\id} f_*\sExt^d_P(i_*\sO_X,\omega_{P/Y})\otimes \omega_{Y/S}.}
All together we see that $\tau_f^n$ equals the composition 
\mlnl{f_*\omega_{X/S}\xr{\zeta'_{i,\pi_P}} f_*(\omega_{X/P}\otimes i^*\omega_P)\cong f_*(\omega_{X/P}\otimes i^*\omega_{P/Y})\otimes \omega_{Y/S}\xr{\eta_i^{-1}\otimes \id } \\f_*(i^\flat\pi^\#\sO_Y)\otimes \omega_{Y/S}
       \xr{\psi_{i,\pi}^{-1}} f_*f^\flat\sO_Y\otimes \omega_Y\cong\sHom_Y(f_*\sO_X,\sO_Y)\otimes \omega_{Y/S}\xr{\text{eval. at } 1} \omega_{Y/S}. }
Hence the claim follows from the definition of $\zeta'_{i,\pi_P}$ (see (1) above) and the definition of the residue symbol in \cite[(A.1.4)]{Co}.
\end{proof}

\end{document}